%% file: article.tex
\documentclass[11pt]{amsart}
\usepackage[utf8]{inputenc}
\usepackage{natbib}
\usepackage{graphicx}
\usepackage{amsmath}
\usepackage{amsthm}
\usepackage{amssymb}
\usepackage{lmodern}
\usepackage{microtype}
\usepackage{enumerate}
\usepackage{tikz}
\usepackage{tikz-cd}
\usepackage{hyperref}
\usepackage{cleveref}
\usepackage{xcolor}
\usepackage{thmtools}
\usepackage{marvosym}
\usepackage{subfiles}
\usepackage{textcmds}
\usepackage{verbatim}
\usepackage{mathtools}
\usepackage[a4paper, left=3cm, right=3cm, top=2.5cm, bottom=2.5cm]{geometry}

\declaretheoremstyle[
spaceabove=15pt,
spacebelow=15pt,
headfont=\bfseries,
notefont=\bfseries, notebraces={[}{]},
bodyfont=\selectfont,
headpunct={},
headformat={\NUMBER~\NAME~\NOTE},
]{sub_part_style}

\declaretheoremstyle[
spaceabove=15pt,
spacebelow=15pt,
headfont=\bfseries,
notefont=\bfseries, notebraces={[}{]},
bodyfont=\selectfont,
headpunct={},
headformat={\NAME~\NUMBER.\NOTE},
]{sub_point_style}

\declaretheoremstyle[
spaceabove=15pt,
spacebelow=15pt,
headfont=\bfseries, notebraces={[}{]},
notefont=\bfseries,
bodyfont=\selectfont,
headpunct={},
headformat={\NAME.\NOTE},
qed=\qedsymbol
]{proof_style}

\declaretheorem[numbered=yes, title=Theorem, style=sub_point_style, numberwithin=section, refname={Theorem,Theorems}]{theorem}
\declaretheorem[numbered=yes, title=Lemma, sibling=theorem,style=sub_point_style, refname={Lemma,Lemmata}]{lemma}
\declaretheorem[numbered=yes, title=Proposition, sibling=theorem, style=sub_point_style, refname={Proposition,Propositions}]{proposition}
\declaretheorem[numbered=yes, title=Corollary, sibling=theorem, style=sub_point_style, refname={Corollary,Corolarrys}]{corollary}
\declaretheorem[numbered=yes, title=Definition, sibling=theorem, style=sub_point_style, refname={Definition,Definitions}]{definition}

\declaretheorem[numbered=no, title=Remark, style=sub_part_style]{remark*}
\declaretheorem[numbered=no, title=Proof, style=proof_style]{proof*}
\crefname{section}{Section}{Sections}

\title{Hecke Algebras for Minimal Equicontinuous Cantor Actions}
\author{Pablo Lummerzheim}
\address{Mathematical Institute, Leiden University, P.O. Box 9512, 2300 RA Leiden,
	The Netherlands}
\email{p.lummerzheim@math.leidenuniv.nl}

\begin{document}
\maketitle

\begin{abstract}
	In this work we extend the concept of Hecke algebras to a pair of a separable infinite profinite group and a closed core-free subgroup. Motivated by applications to minimal equicontinuous Cantor actions, which are characterized by such pairs, we characterize isometric \mbox{*-algebra} isomorphisms between such Hecke algebras by showing that they always arise as Koopman operators of homeomorphisms between the double coset spaces combined with the multiplication by some character. We then prove that two minimal equicontinuous Cantor actions with isometric isomorphic Hecke algebras have the same Steinitz order. We provide examples including the infinite dihedral group and the Grigorchuk group.
\end{abstract}

\input{ch-Introduction.tex}
\input{ch-Cantor_Actions.tex}

\input{ch-Hecke_Algebras.tex}

\input{ch-Hecke_Algebra_as_Invariant.tex}

\input{ch-Examples.tex}

\input{ch-Appendix.tex}

\section*{Acknowledgment}
I like to express my gratitude to Jaime Gómez Ortiz, Max Schmidt, Dean Wardell and to my supervisor Olga Lukina for the many helpful discussions and useful hints.

\bibliographystyle{alpha}
\bibliography{generalized_Hecke_pairs.bib}

\end{document}

%% file: ch-Introduction.tex
\section{Introduction}\label{sec-Introduction}

In this article we study Hecke algebras for pairs \((\widehat{G}, D)\), where \(\widehat{G}\) is a separable infinite profinite group and \(D\) a closed core-free subgroup, where core-free means that \(D\) has no non-trivial normal subgroup. The motivation and applications come from the study of minimal equicontinuous Cantor actions. This connection is explained further below. Hecke algebras have applications to representation theory (\citep{hall_hecke_1999}) and number theory (\citep{bost_hecke_1995}).

In the literature Hecke algebras are mostly studied for pairs \((\widehat{G}, D)\), where either \(\widehat{G}\) is a countable group with a subgroup \(D\) or \(\widehat{G}\) is a totally disconnected group and \(D\) a compact open subgroup (see \citep{hall_hecke_1999,krieg_hecke_1990,glockner_topologization_2001}). Additionally it is required that the left \(D\)-orbit of each coset in \(\widehat{G}/D\) is finite. Such pairs \((\widehat{G}, D)\) are called \textbf{Hecke pairs}.

A pair \((\widehat{G}, D)\) consisting of a separable infinite profinite group \(\widehat{G}\) and a closed subgroup \(D\) only fits into this setting if \(D\) is open, in which case \(\widehat{G}/D\) is discrete. If \(D\) is not open, then \(\widehat{G}/D\) is not discrete, so we can in general not expect that the \(D\)-orbits of \(\widehat{G}/D\) are finite. However, by virtue of \(D\) being closed, thus compact, it follows that the \(D\)-orbits in \(\widehat{G}/D\) are compact. It turns out that this suffices to define a Hecke algebra for a pair of a separable infinite profinite group \(\widehat{G}\) and a closed subgroup \(D\). We additionally restrict to \(D\) being core-free. We explain in the beginning of \cref{sec-Hecke_Algebras}, why we do not lose any generality by making this restriction.

Now fix such a pair \((\widehat{G}, D)\). Let \(\mu\) be the Haar probability measure on \(\widehat{G}\) and \(\nu\) the pushforward measure of the map \(\hat{g} \mapsto D\hat{g}D\) on the double coset space \(D\backslash\widehat{G}/D\). The Hecke algebra \(\mathcal{H}(\widehat{G}, D)\) consists of all \(\operatorname{L}^1\)-functions on \(D\backslash\widehat{G}/D\) as in the definition below. If \(D\) is trivial then the Hecke algebra coincides with the \(\operatorname{L}^1\)-group algebra of \(\widehat{G}\). This is in particular the case if \(\widehat{G}\) is abelian.

\begin{definition}
	The \textbf{Hecke algebra} of \((\widehat{G}, D)\) is the Banach-*-algebra \(\mathcal{H}(\widehat{G}, D) := \operatorname{L}^1(\nu)\) with convolution defined by
	\[
		(\varphi * \psi)(D\hat{g}D) := \int_{\widehat{G}} \varphi(D\hat{f}D) \cdot \psi(D\hat{f}^{-1}\hat{g}D)\ \mathrm{d}\mu(\hat{f})
	\]
	and involution defined by
	\[
		\varphi^*(D\hat{g}D) := \overline{\varphi(D\hat{g}^{-1}D)}
	\]
	for all \(\varphi, \psi \in \mathcal{H}(\widehat{G}, D)\) and \(\nu\)-almost all \(D\hat{g}D \in D\backslash\widehat{G}/D\).
\end{definition}

In \cref{the-*-Algebra_Structure_H} we prove that \(\mathcal{H}(\widehat{G}, D)\) indeed satisfies all the properties of a Banach-*-algebra and prove that the convolution is well-defined. While these properties can be proven by direct computation, we use an approximation by finite dimensional Hecke algebras of so-called group chains. Group chains are fundamental to the study of minimal equicontinuous Cantor actions and are explained in \cref{sec-Cantor_Actions}.

In the main theorem we characterize isometric *-algebras isomorphisms between such Hecke algebras as Koopman operators of measure preserving homeomorphisms combined with the multiplication by a character, i.e. a continuous group homomorphism into the unit circle \(\mathbb{S}^1 \subseteq \mathbb{C}\). This theorem is a Hecke algebra analog to a theorem by Wendel on group algebras (\citep[Theorem 1.]{wendel_isometric_1951}).

\begin{theorem}[\cref{the-Isometric_Isomorphisms_Hecke_Algebras}]\label{the-Isometric_Isomorphic_Hecke_algebras-INTRO}
	Let \(\widehat{G}\), \(\widehat{H}\) be two separable infinite profinite groups with closed core-free subgroups \(D \leq \widehat{G}\) and \(E \leq \widehat{H}\). Assume that \(\mathcal{T}: \mathcal{H}(\widehat{G}, D) \to \mathcal{H}(\widehat{H}, E)\) is an isometric *-algebra isomorphism. Then there exist a unique continuous group homomorphism \(\chi: \widehat{H} \to \mathbb{S}^1\) with \(E \leq \ker(\chi)\) and a unique measure preserving homeomorphism \(\theta: E\backslash\widehat{H}/E \to D\backslash\widehat{G}/D\) such that
	\[
	\mathcal{T}\varphi(E\hat{h}E) = \chi(\hat{h}) \cdot \varphi(\theta(E\hat{h}E))
	\]
	for all \(\varphi \in \mathcal{H}(\widehat{G}, D)\) and \(\nu_H\)-almost all \(E\hat{h}E \in E\backslash\widehat{H}/E\).
\end{theorem}

Now we turn to minimal equicontinuous Cantor actions. A Cantor set \(X\) is a non-empty compact totally disconnected perfect (i.e. no isolated points) metrizable space. A Cantor action \((X, G, \Phi)\) is given by a countable infinite group \(G\) and a group homomorphism \(\Phi: G \to \operatorname{Homeo}(X)\). We abbreviate \(\Phi(g)x\) by \(gx\). Similarly for a subset \(U \subseteq X\) we write \(gU\) for the translate of \(U\) by \(g\).

\begin{definition}\label{def-Equicontinuous_Minimal_Action}
	Let \(d_X\) be any metric on \(X\) that is compatible with its topology. The action \((X, G, \Phi)\) is \textbf{equicontinuous} if for any \(\varepsilon > 0\) there exists some \(\delta > 0\) such that \(d_X(x, y) < \delta\) implies \(d_X(gx, gy) < \varepsilon\) for all \(g \in G\) and all \(x, y \in X\). The action \((X, G, \Phi)\) is \textbf{minimal} if the orbits \(\{gx: g \in G\}\) are dense for all \(x \in X\).
\end{definition}

The definition of equicontinuity is independent of the choice of a compatible metric as \(X\) is compact. Throughout this paper we assume all actions to be minimal equicontinuous Cantor actions. Such actions are studied in many works, see \citep{cortez_g_2008,dyer_discriminant_2016,giordano_mathbbzd-odometers_2019} and references therein. Each minimal equicontinuous Cantor action can be represented as a group action on a rooted tree, see \citep{groger_measures_2021} for the construction. More details on groups acting on rooted trees can be found in \citep{nekrashevych_self-similar_2005,grigorchuk_topics_2011}.

For the action \((X, G, \Phi)\) we denote by \(\widehat{G}\) is the closure of \(\Phi(G)\) in the topology of uniform convergence in \(\operatorname{Homeo}(X)\). Since our action is equicontinuous \(\widehat{g}\) is a compact group and we call it the \textbf{Ellis group}. In the case of a Cantor action the Ellis group \(\widehat{G}\) is topologically a Cantor set, in particular it is a separable infinite profinite group (see \citep[Chapter 3.]{auslander_minimal_1988} and \citep[Chapter 2.]{dyer_discriminant_2016}). Fix any \(x \in X\) and denote by \(D\) the stabilizer of \(x\) in \(\widehat{G}\), i.e.
\[
	D := \{\hat{g} \in \widehat{G}: \hat{g}x = x\}.
\]
It is a closed core-free subgroup of the Ellis group \(\widehat{G}\) (\citep[Proposition 5.3.]{dyer_discriminant_2016}). We refer to \(D\) as the \textbf{discriminant group}. The idea to consider the pair \((\widehat{G}, D)\) as a classification invariant for the minimal equicontinuous Cantor actions was introduced in \citep{dyer_discriminant_2016} and developed in the further works of the authors. The Ellis group \(\widehat{G}\) acts naturally on the homogeneous space \(\widehat{G}/D\), and since \(\Phi\) is a group homomorphism from \(G\) into \(\widehat{G}\), we also obtain a \(G\)-action on \(\widehat{G}/D\). We have the following structure theorem.

\begin{theorem}[{\citealp[Theorem 3.6]{auslander_minimal_1988}}]\label{the-Structure_Theorem_Equicontinuous_Actions}
	The actions of \(G\) on \(X\) and on \(\widehat{G}/D\) are conjugate.
\end{theorem}

Such a conjugacy of the action \((X, G, \Phi)\) is called a \textbf{homogeneous model}. It follows from the construction in \citep[Lemma 3.7.]{fokkink_homogeneous_2002} and \citep[Proposition 10.1.]{dyer_molino_2017} attributed to Lenstra that every pair of a separable infinite profinite group \(\widehat{G}\) with a closed core-free subgroup can be identified with the Ellis group and discriminant group of some minimal equicontinuous Cantor action, for some countable dense subgroup \(G\) of \(\widehat{G}\).

We define the Hecke algebra of a minimal equicontinuous Cantor action \((X, G, \Phi)\) as the Hecke algebra \(\mathcal{H}(\widehat{G}, D)\), where \(\widehat{G}\) is the Ellis group and \(D\) the discriminant group. The subgroup \(D\) is the stabilizer of a point \(x\), so it depends on the choice of this basepoint \(x \in X\), however since \(\widehat{G}\) acts transitively on \(X\) the stabilizers of different base points are conjugate. More precisely, if \((\widehat{H}, E)\) is another homogeneous model of the same action \((X, G, \Phi)\), then there exists a topological group isomorphism \(\theta: \widehat{H} \to \widehat{G}\) with \(\theta(E) = D\) (see \cref{lem-G_U_G_V_isomorphic}). For such isomorphisms we have the theorem below.

\begin{theorem}[\cref{the-Isomorphic_Hecke_Pairs}]
	Let \(\widehat{G}\) and \(\widehat{H}\) be two separable infinite profinite groups with closed core-free subgroups \(D \leq \widehat{G}\) and \(E \leq \widehat{H}\). Assume that there exists a topological group isomorphism \(\theta: \widehat{H} \to \widehat{G}\) with \(\theta(E) = D\). Then there exists an isometric *-algebra isomorphism \(\mathcal{T}: \mathcal{H}(\widehat{G}, D) \to \mathcal{H}(\widehat{H}, E)\).
\end{theorem}

Thus we may speak of \textbf{the Hecke algebra of a minimal equicontinuous Cantor action \((X, G, \Phi)\)}, which we denote by \(\mathcal{H}(\Phi) := \mathcal{H}(\widehat{G}, D)\). This makes the Hecke algebra an invariant for minimal equicontinuous Cantor actions.

We compare the Hecke algebra to another invariant of minimal equicontinuous Cantor actions, the Steinitz order. We prove that the Hecke algebra is a stronger invariant than the Steinitz order.

A \textbf{Steinitz number} or \textbf{supernatural number} is a formal infinite product
\[
	\xi = \prod_{p \in \Pi} p^{\chi_\xi(p)}
\]
 over the set of all primes \(\Pi\), where \(\chi_\xi: \Pi \to \mathbb{N}_0 \cup \{\infty\}\), i.e. \(\chi_\xi(p)\) denotes the multiplicity of the prime factor \(p\). In contrast to natural numbers a Steinitz number may have infinitely many prime factors and a prime factor may have infinite multiplicity. For Steinitz numbers \((\xi_n)_{n \in \mathbb{N}}\) we can define the least common multiple \(\xi = \operatorname{lcm} \{\xi_n\}\) by
\[
	\chi_\xi(p) := \sup\{\chi_{\xi_n}(p): n \in \mathbb{N}\}
\]
for all \(p \in \Pi\). This allows for a generalization of the concept of the index of subgroup for profinite groups. Namely for a profinite group \(\widehat{G}\) and a closed subgroup \(D\) we define the \textbf{profinite index} by
\[
	[\widehat{G}:D] := \operatorname{lcm}\big\{|\widehat{G}/N|: N \text{ open normal subgroup of } \widehat{G},\ D \subseteq N\big\}.
\]
For more details on this see \citep[Section 2.3]{ribes_profinite_2000}. Now we can define the Steinitz order of an equicontinuous Cantor action \((X, G, \Phi)\).

\begin{definition}[{\citealp[Proposition 3.7.]{hurder_type_2026}}]\label{def-Steinitz-Order}
	The \textbf{Steinitz order} of the minimal equicontinuous Cantor action \((X, G, \Phi)\) is defined by
	\[
		\xi(X, G, \Phi) := [\widehat{G}:D],
	\]
	where \((\widehat{G}, D)\) is a homogeneous model of \((X, G, \Phi)\).
\end{definition}

This definition is independent of the choice of a homogeneous model (\citep[Theorem 1.5.]{hurder_type_2026}). In \citep[Theorem 1.5]{hurder_type_2026} it is shown that the order of return equivalent actions differs only by multiplication by an integer. In \citep[Theorem 1.5.]{giordano_mathbbzd-odometers_2019} it is proven that two \(\mathbb{Z}^d\)-actions are orbit equivalent if and only if they have the same Steinitz order. We show that the Steinitz order is preserved under isometric *-algebra isomorphisms.

\begin{theorem}[\cref{cor-Steinitz_Order}]\label{the-Steinitz_Order}
	Let \((X, G, \Phi)\) and \((Y, H, \Psi)\) be two minimal equicontinuous Cantor actions with Steinitz orders \(\xi(X, G, \Phi)\) and \(\xi(Y, H, \Psi)\) respectively. Assume that there exists an isometric *-algebra isomorphism \(\mathcal{T}: \mathcal{H}(\Phi) \to \mathcal{H}(\Psi)\) between their Hecke algebras. Then
	\[
		\xi(X, G, \Phi) = \xi(Y, H, \Psi).
	\]
\end{theorem}

This result gives immediate access to a wide range of example. In \citep[Section 5B. and 5C.]{hurder_prime_2023} the authors construct minimal equicontinuous Cantor actions of the discrete Heisenberg group with any Steinitz order. Combining this result with \cref{the-Steinitz_Order} we obtain the following corollary.

\begin{corollary}\label{cor-Hecke_Algebra_Heisenberg_Group}
	Let \(G\) be the discrete Heisenberg group. There are uncountably many minimal equicontinuous Cantor actions of \(G\) that have non-isometric isomorphic Hecke algebras.
\end{corollary}

In the last section of our paper we consider three examples. All of them are given by groups acting on the binary tree. The first one is an action of the infinite dihedral group \(H\). For it we calculate finite dimensional Hecke algebras that approximate the Hecke algebra of the action (\cref{the-Hecke_Algebra_Dihedral_Group}).

The other two examples are the Grigorchuk group \(G\) and another group \(\widetilde{G}\) constructed and discussed in \citep{bartholdi_spectrum_2000, bartholdi_parabolic_2002}. Both are defined as groups of automorphism of the binary tree and are branch groups. In \citep{bartholdi_parabolic_2002} finite dimensional Hecke algebras that correspond to their actions were already discussed. We show in \cref{pro-Double_Coset_Grigorchuk} that for both groups the double coset space \(D\backslash\widehat{G}/D\) of the Ellis group \(\widehat{G}\) and the discriminant group \(D\) is homeomorphic to the one point compactification \(\mathbb{N}_0 \cup \{\infty\}\) of the natural numbers. Using this we describe the Hecke algebras of their minimal equicontinuous Cantor actions, i.e. the Hecke algebras associated to \((\widehat{G}, D)\), explicitly in \cref{the-Hecke_Algebra_Grigorchuk's_Group}.

We then see that the actions of \(G\) and \(\widetilde{G}\) give rise to isometric isomorphic Hecke algebras, while the Hecke algebra of the action of the dihedral group \(H\) is not isometric isomorphic to the Hecke algebras of the actions of \(G\) and \(\widetilde{G}\). However, the actions of \(G\) and \(H\) have the same Steinitz order, namely \(2\) is the only prime factor, which has infinite multiplicity. The example of \(G\) and \(\widetilde{G}\) shows that the Hecke algebra is not a complete invariant of minimal equicontinuous Cantor actions, since two non-conjugate actions can have isometric *-isomorphic Hecke algebras. On the other hand they are a finer invariant than the Steinitz order, since the actions of \(H\) and \(G\) have the same Steinitz order but non-isometric *-isomorphic Hecke algebras. This shows that \cref{the-Steinitz_Order} has no converse.

The structure of the rest of the paper is the following: In \cref{sec-Cantor_Actions} we give necessary construction and statements for minimal equicontinuous Cantor action. \cref{sec-Cantor_Actions} consists of three parts, the first part covers group chains, the second double coset spaces and in the third part probability measures on various groups and double coset spaces are introduced. Hecke algebras are defined in \cref{sec-Hecke_Algebras}. In the first part of that section finite dimensional Hecke algebras are discussed. In the second part the Hecke algebra for minimal equicontinuous Cantor actions and for pairs consisting of a separable infinite profinite group with a closed core-free subgroup is defined in \cref{the-*-Algebra_Structure_H}. In \cref{sec-Invariant} we explore how the Hecke algebra can be used as an invariant, proving our main result in \cref{the-Isometric_Isomorphisms_Hecke_Algebras}. In the second part of \cref{sec-Invariant} we consider the applications of this theorem to minimal equicontinuous Cantor actions. In \cref{sec-Examples} we discuss the examples, namely the dihedral group in the first part and Grigorchuk group in the second part. In the end there is an appendix, where some theorems used in this work are discussed.

%% file: ch-Cantor_Actions.tex
\section{Cantor Actions}\label{sec-Cantor_Actions}

In this section we give the necessary definitions and statements concerning minimal equicontinuous Cantor actions introduced in \cref{def-Equicontinuous_Minimal_Action}. Throughout this section \((X, G, \Phi)\) always denotes a minimal equicontinuous Cantor action. In the first part we discuss how such an action can be described using group chains. In the second part we turn to double coset spaces for group chains. Those spaces are essential for defining Hecke algebras in the next section. In the final part of this section we add a measure theoretic perspective by equipping various groups and spaces with probability measures.

\subsection{Group Chains}\label{sub-Group_Chains}

In the case of minimal equicontinuous Cantor actions the Ellis group \(\widehat{G}\) and the discriminant group \(D\) can be more explicitly described using group chains. Also group chains provide a convenient way to construct a minimal equicontinuous Cantor action for a given group \(G\). More details on group chains can be found in \citep{dyer_discriminant_2016}.

\begin{definition}
	Let \(G\) be a countable group. A \textbf{group chain} is a strictly decreasing sequence \((G_n)_{n \in \mathbb{N}_0}\) of finite index subgroups of \(G\) with \(G_0 = G\).
\end{definition}

With such a group chain given we can define a Cantor space as the inverse limit of finite coset spaces
\[
	\varprojlim_{n \in \mathbb{N}_0} G/G_n = \bigg\{(g_kG_k)_{k \in \mathbb{N}_0} \in \prod_{k \in \mathbb{N}_0} G/G_k: g_lG_k = g_kG_k\ \forall l \geq k \geq 0\bigg\}.
\]
A metric can be naturally defined on this space by letting the distance between two sequences be \(2^{-k}\), where \(k\) is the greatest integer on which the sequences agree. For \(g \in G\) we define its action on \((g_kG_k)_{k \in \mathbb{N}_0} \in \varprojlim_{n \in \mathbb{N}_0} G/G_n\) by
\begin{equation}\label{eq-G-action_of_Group_Chain}
 	g(g_kG_k)_{k \in \mathbb{N}_0} := (gg_kG_k)_{k \in \mathbb{N}_0}.
\end{equation}
It is immediate to check that this is a well defined \(G\)-action. The group action is isometric with respect to the metric defined above, in particular it is equicontinuous. Also the \(G\)-action on \(G/G_n\) is minimal for all \(n \in \mathbb{N}_0\), from which it follows that also the \(G\)-action on \(\varprojlim_{n \in \mathbb{N}_0} G/G_n\) is minimal. We call the action in \eqref{eq-G-action_of_Group_Chain} the minimal equicontinuous Cantor action associated to the group chain \((G_n)_{n \in \mathbb{N}_0}\). Conversely if a minimal equicontinuous action \((X, G, \Phi)\) is given we can construct a group chain associated to it using adapted sets.

\begin{definition}\label{def-Adapted_Set}
	A clopen subset \(U \subseteq X\) is called \textbf{adapted} if \(\{gU: g \in G\}\) is a finite partition of \(X\).
\end{definition}

\begin{proposition}[{\citealp[Proposition A.1.]{dyer_discriminant_2016}}]\label{pro-Existence_Group_Chain}
	Let \((X, G, \Phi)\) be a minimal equicontinuous Cantor action. Fix \(x \in X\). Then there exists a strictly decreasing clopen neighborhood basis \((U_n)_{n \in \mathbb{N}_0}\) of \(x\) that consists of adapted sets with \(U_0 = X\). Then
	\[
		G_n := \{g \in G: gU_n = U_n\}
	\]
	is a finite index subgroup for all \(n \in \mathbb{N}_0\) and \((G_n)_{n \in \mathbb{N}_0}\) forms a group chain of \(G\). The actions of \(G\) on \(X\) and on \(\varprojlim_{n \in \mathbb{N}_0} G/G_n\) (see \eqref{eq-G-action_of_Group_Chain}) are conjugate, i.e. there exists a homeomorphism \(\theta: X \to \varprojlim_{n \in \mathbb{N}_0} G/G_n\) such that
	\[
		\theta(gy) = g\theta(y)
	\]
	for all \(g \in G\) and \(y \in X\) and \(\theta(x) = (eG_k)_{k \in \mathbb{N}_0}\).
\end{proposition}

We call the point \(x\) the \textbf{base point} of the group chain. For a minimal equicontinuous Cantor action \((X, G, \Phi)\) fix a base point \(x\) and a neighborhood basis \((U_n)_{n \in \mathbb{N}_0}\) of \(x\) consisting of adapted sets and let \((G_n)_{n \in \mathbb{N}_0}\) the group chain arising from this neighborhood basis. Denote by \(\widehat{G}\) its Ellis group and by \(D\) its discriminant group, i.e. the stabilizer of the base point \(x\). For \(n \in \mathbb{N}_0\) we define
\[
	C_n := \operatorname{core}_G(G_n) := \bigcap_{g \in G} g^{-1}G_ng,
\]
which are normal subgroups of \(G\) with finite index. For \(m \geq n\geq 0\) we have \(C_m \leq C_n\) since \(G_m \leq G_n\), so we can define the coset inclusion \(p_{n, m}:G/C_m \to G/C_n\) by
\[
	p_{n, m}(gC_m) := gC_n
\]
for all \(gC_m \in G/C_m\). These maps are surjective group homomorphisms. They form an inverse system, i.e. \(p_{n, k} = p_{n, m} \circ p_{m, k}\) for \(k \geq m \geq n \geq 0\) (see \citep[Section 1.1]{ribes_profinite_2000}). Thus we can define the inverse limit
\[
	\varprojlim_{n \in \mathbb{N}_0} G/C_n = \bigg\{(g_kC_k)_{k \in \mathbb{N}_0} \in \prod_{k \in \mathbb{N}_0} G/C_k: p_{k, l}(g_lC_l) = g_kC_k\ \forall l \geq k \geq 0\bigg\}.
\]
This is a profinite group with multiplication given by
\[
	(g_kC_k)_{k \in \mathbb{N}_0}(h_kC_k)_{k \in \mathbb{N}_0} = (g_kh_kC_k)_{k \in \mathbb{N}_0}.
\]
There is a natural group homomorphism \(G \to \varprojlim_{n \in \mathbb{N}_0} G/C_n\) given by \(g \mapsto (gC_k)_{k \in \mathbb{N}_0}\). It has dense image in \(\varprojlim_{n \in \mathbb{N}_0} G/C_n\). It is injective if and only if the action \((X, G, \Phi)\) is effective, i.e. if \(g \in G\) fulfills \(gy = y\) for all \(y \in X\) if and only if \(g\) is the identity element.

For \(m \geq n \geq 0\) we have \(p_{n, m}(G_m/C_m) \subseteq G_n/C_n\) (see \citep[Lemma 4.2]{dyer_discriminant_2016}), so we can also define the group
\[
	\varprojlim_{n \in \mathbb{N}_0} G_n/C_n = \bigg\{(g_kC_k)_{k \in \mathbb{N}_0} \in \prod_{k \in \mathbb{N}_0} G_k/C_k: p_{k, l}(g_lC_l) = g_kC_k\ \forall l \geq k \geq 0\bigg\},
\]
which is a subgroup of \(\varprojlim_{n \in \mathbb{N}_0} G/C_m\).

\begin{theorem}[{\citealp[Theorem 1.8]{dyer_discriminant_2016}}]\label{the-Inverse_Limit_Ellis_Group}
	There exists a topological group isomorphism \(\theta: \widehat{G} \to \varprojlim_{n \in \mathbb{N}_0} G/C_n\) such that \(\theta(D) = \varprojlim_{n \in \mathbb{N}_0} G_n/C_n\).
\end{theorem}

An immediate consequence of this is that the topological group \(\widehat{G}\) is topologically a Cantor set. From now on we identify \(\widehat{G}\) with \(\varprojlim_{n \in \mathbb{N}_0} G/C_n\) and the discriminant group \(D\) with \(\varprojlim_{n \in \mathbb{N}_0} G_n/C_n\). We summarize all previously defined maps in the following commuting diagram of topological groups for \(m \geq n \geq 0\):
\begin{equation}\label{dia-Inverse_Limit_Ellis_Group}
	\begin{tikzcd}
		D \arrow[rrr, hook] \arrow[rd, "p_m"] \arrow[rdd, bend right, "p_n"] & & &\widehat{G} \arrow[ld, ->>, "p_m"] \arrow[ldd, bend left, ->>, "p_n"]\\
		&G_m/C_m \arrow[r, hook] \arrow[d, "p_{n, m}"] &G/C_m \arrow[d, ->>, "p_{n, m}"]\\
		&G_n/C_n \arrow[r, hook] &G/C_n
	\end{tikzcd}
\end{equation}
Here the unlabeled hooked arrows are set inclusions and the double tipped arrows correspond to surjections. The group homomorphisms \(p_n: \widehat{G} \to G/C_n\) are defined by
\[
	p_n((g_kC_k)_{k \in \mathbb{N}_0}) := g_nC_n
\]
for all \((g_kC_k)_{k \in \mathbb{N}_0} \in \widehat{G}\) and \(n \in \mathbb{N}_0\).

A group chain associated to \((X, G, \Phi)\) depends on a choice of basepoint \(x\) and a choice of adapted neighborhood basis and so it is not unique. The theorem below gives a classification of all group chains that give rise to conjugate Cantor actions.

\begin{theorem}[{\citealp[Theorem 1.4.]{dyer_discriminant_2016}}]\label{the-Equivalence_Group_Chains}
	Let \((G_n)_{n \in \mathbb{N}_0}\) and \((G_n')_{n \in \mathbb{N}_0}\) be two group chains in \(G\) with Ellis groups \(\widehat{G}\) and \(\widehat{G}'\) as well as their discriminant groups \(D\) and \(D'\). Then
	\begin{enumerate}[(a)]
		\item there exists a conjugacy \(\phi: \widehat{G}/D \to \widehat{G}'/D\) with \(\phi(eD) = eD'\) if and only if the group chains \((G_n)_{n \in \mathbb{N}_0}\) and \((G_n')_{n \in \mathbb{N}_0}\) are \textbf{equivalent}, i.e. there are strictly increasing sequences \((n_k)_{k \in \mathbb{N}_0}\) and \((m_k)_{k \in \mathbb{N}_0}\) in \(\mathbb{N}_0\) such that \(G_{n_k} > G_{m_k}'\) and \(G_{m_k}' > G_{n_{k+1}}\) for all \(k \in \mathbb{N}_0\)
		\item there exists a conjugacy \(\phi: \widehat{G}/D \to \widehat{G}'/D\) if and only if the group chains \((G_n)_{n \in \mathbb{N}_0}\) and \((G_n')_{n \in \mathbb{N}_0}\) are \textbf{conjugate equivalent}, i.e. there are \((g_n)_{n \in \mathbb{N}_0}\) in \(G\) with \(g_nG_n = g_mG_n\) for all \(m \geq n \geq 0\) such that the group chains \((g_nG_ng_n^{-1})_{n \in \mathbb{N}_0}\) and \((G_n')_{n \in \mathbb{N}_0}\) are equivalent.
	\end{enumerate}
\end{theorem}

In the introduction we introduced Steinitz numbers and defined the Steinitz order of a minimal equicontinuous Cantor action in \cref{def-Steinitz-Order}. It can be easily computed using group chains and does not depend on the choice of group chain.

\begin{proposition}[{\citealp[Theorem 1.5.]{hurder_type_2026}}]\label{pro-Steinitz_Order}
	Let \((X, G, \Phi)\) be a minimal equicontinuous Cantor action with associated group chain \((G_n)_{n \in \mathbb{N}_0}\). Then
	\[
		\xi(X, G, \Phi) = \operatorname{lcm}\{|G/G_n|: n \in \mathbb{N}_0\}.
	\]
\end{proposition}

\subsection{Double Coset Spaces}\label{sub-Double_Coset_Spaces}

Hecke algebras are defined on double cosets spaces, so we take a look at those. Fix a minimal equicontinuous Cantor action \((X, G, \Phi)\) with associated group chain \((G_n)_{n \in \mathbb{N}_0}\), Ellis group \(\widehat{G}\) and discriminant group \(D\). For \(g \in G\) its \(G_n\)-double coset is defined by
\[
	G_ngG_n = \{hgh': h, h' \in G_n\} \subseteq G
\]
and we denote by \(G_n\backslash G/G_n\) the space of all double cosets of \(G_n\). There is a one-to-one correspondence between double cosets and left \(G_n\)-orbits in the coset space \(G/G_n\). We do not strictly distinguish between these two notions, we will however occasionally write
\[
	G_n(gG_n) := \{hgG_n: h \in G_n\} \subseteq G/G_n
\]
if we want to emphasize that we consider a double coset as a set of single cosets. The same notation is also used for \(D\backslash \widehat{G}/D\). We equip \(D\backslash\widehat{G}/D\) with the quotient topology arising from the quotient map \(\hat{g} \mapsto D\hat{g}D\). Since \(D\) is closed it follows that each double coset \(D\hat{g}D\) is a closed subset of \(\widehat{G}\), thus by \citep[10.4 Proposition 8.]{bourbaki_general_1995} we obtain that \(D\backslash\widehat{G}/D\) is compact Hausdorff.

For \(n \in \mathbb{N}_0\) we define the inclusion map \(r_n: G/C_n \to G_n\backslash G/G_n\) by
\[
	r_n(gC_n) := G_ngG_n
\]
for all \(gC_n \in G/C_n\), and for \(m \geq n \geq 0\) the double coset inclusion \(q_{n, m}: G_m\backslash G/G_m \to G_n\backslash G/G_n\) by
\[
	q_{n, m}(G_mgG_m) := G_ngG_n
\]
for all \(G_mgG_m \in G_m\backslash G/G_m\). Both maps are surjective. Also we define the surjective map \(r: \widehat{G} \to D\backslash\widehat{G}/D\) by
\[
	r(\hat{g}) := D\hat{g}D
\]
for all \(\hat{g} \in \widehat{G}\). At last we define surjective maps \(q_n: D\backslash\widehat{G}/D \to G_n\backslash G/G_n\) by
\[
	q_n(D(g_kC_k)_{k \in \mathbb{N}_0}D) := G_ng_nG_n
\]
for all \((g_kC_k)_{k \in \mathbb{N}_0} \in \widehat{G}\) and \(n \in \mathbb{N}_0\). All these maps are well defined. For \(r_n\) and \(q_{n, m}\) this follows from \(C_n \leq G_n\) and \(G_m \leq G_n\). For \(q_n\) we note that if \(D(g_kC_k)_{k \in \mathbb{N}_0}D = D(h_kC_k)_{k \in \mathbb{N}_0}D\) then there are \((d_kC_k)_{k \in \mathbb{N}_0}\) and \((d_k'C_k)_{k \in \mathbb{N}_0}\) in \(D\) such that
\[
	(d_kg_kd_k'C_k)_{k \in \mathbb{N}_0} = (d_kC_k)_{k \in \mathbb{N}_0}(g_kC_k)_{k \in \mathbb{N}_0}(d_k'C_k)_{k \in \mathbb{N}_0} = (h_kC_k)_{k \in \mathbb{N}_0}.
\]
Thus \(d_ng_nd_n'C_n = h_nC_n\) and \(d_n, d_n' \in G_n\) since \((d_kC_k)_{k \in \mathbb{N}_0}, (d_k'C_k)_{k \in \mathbb{N}_0} \in D\) implies with \(C_n \leq G_n\) that
\[
	G_ng_nG_n = G_nh_nG_n,
\]
so \(q_n\) is indeed well defined. It is immediate from the definitions that the following diagram commutes for all \(m \geq n \geq 0\):
\begin{equation}\label{dia-Inverse_Limit_Double_Cosets}
	\begin{tikzcd}
		\widehat{G} \arrow[rrr, ->>, "r"] \arrow[rd, ->>, "p_m"] \arrow[rdd, bend right, ->>, "p_n"] & & &D\backslash \widehat{G}/D \arrow[ld, ->>, "q_m"] \arrow[ldd, bend left, ->>, "q_n"]\\
		&G/C_m \arrow[r, ->>, "r_m"] \arrow[d, ->>, "p_{n, m}"] &G_m\backslash G/G_m \arrow[d, ->>, "q_{n, m}"]\\
		&G/C_n \arrow[r, ->>, "r_n"] &G_n\backslash G/G_n
	\end{tikzcd}
\end{equation}

All the maps in the diagram are continuous. This is clear for \(p_{n, m}, q_{n, m}\) and \(r_n\) since these are maps between finite spaces. The group \(\widehat{G}\), when seen as the inverse limit of the \(G/C_n\), is equipped with the product topology and the maps \(p_n\) are the coordinate projection, hence they are continuous. The space \(D\backslash\widehat{G}/D\) is equipped with the quotient topology arising from the quotient map \(r\), so \(r\) is continuous and from \(q_n \circ r = r_n \circ p_n\) the fundamental property of the quotient topology implies that \(q_n\) is continuous since \(r_n \circ p_n\) is continuous.

Since the maps \(q_{n, m}\) satisfy \(q_{n, k} = q_{n, m} \circ q_{m, k}\) for \(k \geq m \geq n \geq 0\) we see that double coset spaces \(G_n\backslash G/G_n\) with the maps \(q_{n, m}\) form an inverse system of compact topological spaces (\citep[Section 1.1]{ribes_profinite_2000}). In fact its inverse limit is homeomorphic to \(D\backslash\widehat{G}/D\), which we prove in the proposition below. This in particular shows that \(D\backslash\widehat{G}/D\) is a totally disconnected Polish space (\citep[Proposition 1.1.3]{ribes_profinite_2000}). However, it is not necessarily perfect and might thus not be a Cantor set. An example of this is the Grigorchuk group, see \cref{pro-Double_Coset_Grigorchuk}. This happens if for some \(n \in \mathbb{N}_0\) there is a double coset \(G_ngG_n\) such that the preimage \(q_{n, m}^{-1}(\{G_ngG_n\})\) consists of a single double coset element for all \(m \geq n\). This yields an isolated point in the inverse limit \(\varprojlim_{n \in \mathbb{N}_0} G_n\backslash G/G_n\). This can not happen for the inverse limits \(\varprojlim_{n \in \mathbb{N}_0} G/C_n\) and \(\varprojlim_{n \in \mathbb{N}_0} G/G_n\), since their bonding maps are constant-to-one.

\begin{proposition}\label{pro-Inverse_Limit_Double_Cosets}
	The compact topological spaces \(D\backslash\widehat{G}/D\) and \(\varprojlim_{n \in \mathbb{N}_0} G_n\backslash G/G_n\) are homeomorphic.
\end{proposition}
\begin{proof*}
	Denote
	\[
		Z := \varprojlim_{n \in \mathbb{N}_0} G_n\backslash G/G_n = \{(G_kg_kG_k)_{k \in \mathbb{N}_0}: q_{k, l}(G_lg_lG_l) = G_kg_kG_k\ \forall l \geq k \geq 0\}.
	\]
	We define a map \(\theta: D\backslash\widehat{G}/D \to Z\) by
	\[
		\theta(D\hat{g}D) := (q_k(D\hat{g}D))_{k \in \mathbb{N}_0}
	\]
	for all \(D\hat{g}D \in D\backslash\widehat{G}/D\). This is well defined since
	\[
		q_{k, l}(q_l(D\hat{g}D)) \overset{\eqref{dia-Inverse_Limit_Double_Cosets}}{=} q_k(D\hat{g}D)
	\]
	for all \(l \geq k \geq 0\), so \(\theta(D\hat{g}D)\) is indeed in \(Z\). It is continuous since \(q_k\) is continuous for all \(k \in \mathbb{N}_0\). Since all maps \(q_k\) are surjective and \(D\backslash\widehat{G}/D\) is compact Hausdorff it follows by \citep[Corollary 1.1.6]{ribes_profinite_2000} that \(\theta\) is also surjective.
	
	Next we prove that \(\theta\) is injective. Therefore fix \(D\hat{g}D, D\hat{f}D \in D\backslash\widehat{G}/D\) with \(\hat{g} = (g_kC_k)_{k \in \mathbb{N}_0}\) and \(\hat{f} = (f_kC_k)_{k \in \mathbb{N}_0}\), such that they have the same image under \(\theta\). Thus
	\[
		G_ng_nG_n = q_n(D\hat{g}D) = q_n(D\hat{f}D) = G_nf_nG_n
	\]
	for all \(n \in \mathbb{N}_0\), Therefore there are \(d_n, d_n' \in G_n\) such that \(d_ng_nd_n' = f_n\). Then \(\hat{d}_n := (d_nC_k)_{k \in \mathbb{N}}\) and \(\hat{d}_n' := (d_n'C_k)_{k \in \mathbb{N}_0}\) are elements in \(\widehat{G}\). For each \(m \geq n \geq 0\) we have
	\begin{alignat}{3}\label{eq-Conjugacy_m_geq_n}
		p_n(\hat{d}_m\hat{g}\hat{d}_m') &\overset{\eqref{dia-Inverse_Limit_Ellis_Group}}{=} p_{n, m}(p_m(\hat{d}_m\hat{g}\hat{d}_m'))\\
		&= p_{n, m}(d_mg_md_m'C_m)\nonumber\\
		&= p_{n, m}(f_mC_m)\nonumber\\
		&= f_nC_n &= p_n(\hat{f}).\nonumber
	\end{alignat}
	Since \(\widehat{G}\) is compact, we find some subsequence \((n_m)_{m \in \mathbb{N}_0}\) and elements \(\hat{d}, \hat{d}' \in \widehat{G}\) such that
	\[
		\hat{d} = \lim_{m \to \infty} \hat{d}_{n_m} \quad\text{and}\quad \hat{d}' = \lim_{m \to \infty} \hat{d}_{n_m}'.
	\]
	We can use the same subsequence for both sequences by first obtaining a subsequence on which \((\hat{d}_n)_{n \in \mathbb{N}_0}\) converges and then taking a subsequence of this subsequence on which \((\hat{d}_n')_{n \in \mathbb{N}_0}\) converges. For all \(n \in \mathbb{N}_0\) we have for \(n_m \geq n\) that \(p_n(\hat{d}_{n_m}) = d_{n_m}C_n \in G_n/C_n\), since \(d_{n_m} \in G_{n_m} \leq G_n\). This implies, since \(p_n\) is continuous and \(G_n/C_n\) closed in \(G/C_n\), that
	\[
		p_n(\hat{d}) = \lim_{m \to \infty} p_n(\hat{d}_{n_m}) \in G_n/C_n
	\]
	for all \(n \in \mathbb{N}_0\), thus \(\hat{d} \in D\). Similarly \(\hat{d}' \in D\). Since \(p_n\) is continuous we have
	\[
		p_n(\hat{d}\hat{g}\hat{d}') = \lim_{m \to \infty} p_n(\hat{d}_{n_m}\hat{g}\hat{d}_{n_m}') \overset{\eqref{eq-Conjugacy_m_geq_n}}{=} p_n(\hat{f}),
	\]
	which implies \(\hat{d}\hat{g}\hat{d}' = \hat{f}\). Hence \(D\hat{g}D = D\hat{f}D\), which proves that \(\theta\) is also injective, hence bijective. Thus \(\theta\) is a homeomorphism since \(Z\) and \(D\backslash\widehat{G}/D\) are compact Hausdorff.
\end{proof*}

\subsection{Probability Measures on Group Chains}\label{sub-Probability_Measures_Group_Chains}

In this section we equip the various groups and spaces defined above with suitable probability measures. The groups \(G/C_n\) for \(n \in \mathbb{N}_0\) and \(\widehat{G}\) are compact Hausdorff, so we can equip them with Haar probability measures \(\mu_n\) and \(\mu\). These measures are left and right invariant and the taking inverses also preserves the measure. As the groups \(G/C_n\) are finite the measure \(\mu_n\) is given by the counting measure rescaled by the index of \(C_n\). i.e.
\[
	\mu_n(\{gC_n\}) = \frac{1}{|G/C_n|}
\]
for all \(gC_n \in G/C_n\). On \(G_n\backslash G/G_n\) we use the the pushforward measure \(\nu_n := r_{n*}\mu_n\), i.e.
\[
	\nu_n(\{G_ngG_n\}) := \mu_n(r_n^{-1}(\{G_ngG_n\})).
\]
for all \(G_ngG_n \in G_n\backslash G/G_n\). We note that \(r_n(fC_n)\) is equal to \(G_ngG_n\) if and only if \(fG_n\) is in the \(G_n\)-orbit \(G_n(gG_n)\) of \(gG_n\). Also each \(fG_n \in G_n(gG_n)\) contains exactly \(\frac{|G/C_n|}{|G/G_n|}\) cosets of \(C_n\). Therefore
\[
	|r_n^{-1}(\{G_ngG_n\})| = |\{fC_n \in G/C_n: G_nfG_n = G_ngG_n\}| = |G_n(gG_n)| \cdot \frac{|G/C_n|}{|G/G_n|}
\]
and hence by the definition of \(\mu_n\)
\begin{equation}\label{eq-nu_n}
	\nu_n(\{G_ngG_n\}) := \mu_n(r_n^{-1}(\{G_ngG_n\})) = \frac{|G_n(gG_n)|}{|G/G_n|}.
\end{equation}
Similarly we equip \(D\backslash\widehat{G}/D\) with the probability measure
\[
	\nu := r_*\mu.
\]
With these choices of measures it turns out that all the maps in \eqref{dia-Inverse_Limit_Double_Cosets} are measure preserving. This is the main reason why it is advantageous for us to equip \(G/C_n\) with the normalized over the non-normalized counting measure.

\begin{proposition}\label{pro-Maps_Measure_Preserving}
	For \(m \geq n \geq 0\) we have that all the maps \(p_{n, m}\), \(p_n\), \(q_{n, m}\), \(q_n\), \(r_n\) and \(r\) (see \eqref{dia-Inverse_Limit_Double_Cosets}) are measure preserving.
\end{proposition}
\begin{proof*}
	For \(r_n\) and \(r\) this is true by definition of \(\nu_n\) and \(\nu\). For \(p_{n, m}\) we note that since \(C_m \leq C_n\) every coset \(gC_n \in G/C_n\) consists of exactly \(\frac{|G/C_m|}{|G/C_n|}\) cosets in \(G/C_m\), thus
	\begin{alignat*}{3}
		\mu_m(p_{n, m}^{-1}(\{gC_n\})) &= \mu_m(\{fC_m \in G/C_m: fC_m \subseteq gC_n\})\\
		&= \frac{1}{|G/C_m|} \cdot \frac{|G/C_m|}{|G/C_n|}\\
		&= \frac{1}{|G/C_n|}\\
		&= \mu_n(\{gC_n\})
	\end{alignat*}
	for all \(gC_n \in G/C_n\), so \(p_{n, m}\) is indeed measure preserving.
	
	Next we prove that \(p_n\) is measure preserving. Let \(i := |G/C_n|\) be the index of \(C_n\) and \(g_1, \dots, g_i \in G\) representatives of the cosets of \(C_n\). Then \(\{p_n^{-1}(\{g_1C_n\}), \dots, p_n^{-1}(\{g_iC_n\})\}\) is a clopen partition of \(\widehat{G}\). Since \(p_n\) is surjective, there are \(\hat{g}_1, \dots, \hat{g}_i \in \widehat{G}\) such that \(p_n(\hat{g}_j) = g_jC_n\) for all \(j \in \{1, \dots, i\}\). Thus by the invariance of \(\mu\) and since \(p_n\) a group homomorphism we have for \(j, k \in \{1, \dots, i\}\) that
	\[
		\mu(p_n^{-1}(\{g_jC_n\})) = \mu(\hat{g}_k\hat{g}_j^{-1}p_n^{-1}(\{g_jC_n\})) = \mu(p_n^{-1}(g_kg_j^{-1}\{g_jC_n\})) = \mu(p_n^{-1}(\{g_kC_n\})),
	\]
	so all the sets in the partition \(\{p_n^{-1}(\{g_1C_n\}), \dots, p_n^{-1}(g_iC_n)\}\) have the same measure. Thus
	\[
		\mu(p_n^{-1}(\{g_jC_n\})) = \frac{1}{|G/C_n|} = \mu_n(\{g_jC_n\})
	\]
	for all \(j \in \{1, \dots, i\}\), so \(p_n\) is indeed measure preserving.
	
	The fact that \(q_{n, m}\) is measure preserving follows from \(q_{n, m} \circ r_m = r_n \circ p_{n, m}\) (see \eqref{dia-Inverse_Limit_Double_Cosets}) and the fact that \(r_m\), \(r_n\) and \(p_{n, m}\) are measure preserving. Similarly \(q_n\) is measure preserving since \(r\), \(r_n\) and \(p_n\) are and \(q_n \circ r = r_n \circ p_n\).
\end{proof*}

%% file: ch-Hecke_Algebras.tex
\section{Hecke Algebras}\label{sec-Hecke_Algebras}

In this section we define the Hecke algebra for pairs \((\widehat{G}, D)\) consisting of a separable infinite profinite group \(\widehat{G}\) with a closed core-free subgroup \(D\). We note that if \(D\) has a non-trivial core \(C\), then we can consider the pair \((\widehat{G}/C, D/C)\) instead. For this pair we have for the double cosets space that
\[
	(D/C)\backslash(\widehat{G}/C)/(D/C) \cong D\backslash\widehat{G}/D
\]
holds. Thus, as the Hecke algebra consists of functions on the double coset space, we do not lose any generality by restricting to core-free subgroups \(D\).

In the introduction we discussed that any such pair \((\widehat{G}, D)\) with \(D\) core-free arises from a minimal equicontinuous Cantor action \((X, G, \Phi)\), where \(G\) is a countable dense subgroup of \(\widehat{G}\). Thus we assume that \(\widehat{G}\) is the Ellis group of some minimal equicontinuous Cantor action \((X, G, \Phi)\) and \(D\) its discriminant group. Further we fix a group chain \((G_n)_{n \in \mathbb{N}_0}\) that corresponds to this action (see \cref{sec-Cantor_Actions}).

In the first part of this section we look at the finite dimensional Hecke algebras \(\mathcal{H}(G, G_n)\) of the pairs \((G, G_n)\). These pairs fit into the setting of Hecke algebras in the literature (\citep{krieg_hecke_1990,hall_hecke_1999}), however, we modify the definition by adding a normalizing factor. The motivation and consequences of this modification are discussed in \cref{sub-Finite_Dimensional_Hecke_algebras} directly below. In the second part we define the infinite dimensional Hecke algebra \(\mathcal{H}(\widehat{G}, D)\). For this we use an approximation by the finite dimensional Hecke algebras \(\mathcal{H}(G, G_n)\). We see that the Hecke algebra is independent of the choice of a minimal equicontinuous Cantor action associated to a pair of a separable infinite profinite group with closed core-free subgroup and independent of the choice of base point and group chain.

\subsection{Finite Dimensional Hecke Algebras}\label{sub-Finite_Dimensional_Hecke_algebras}

The pair \((G, G_n)\) is a Hecke pair, since \(G_n\) is a finite index subgroup of the countable group \(G\), so we can consider its Hecke algebra. In the literature the Hecke algebra \(\mathcal{H}(G, G_n)\) is defined in two different ways, either as the vector space of all functions \(G \to \mathbb{C}\) that are \(G_n\) left and right invariant (see \citep[Section 2.2.2]{hall_hecke_1999}), or as functions on double cosets(see \citep[Section 1.2]{krieg_hecke_1990}). It is clear that these two definitions are equivalent, since the underlying vector spaces are isomorphic. We use the second version, so for us the Hecke algebra \(\mathcal{H}(G, G_n)\) consists of functions from the double coset space \(G_n\backslash G/G_n\) into \(\mathbb{C}\). However, it might sometimes be helpful to think of them as functions on \(G\).

\begin{definition}\label{def-Hecke_Alegbra_Finite_Dimensional}
	The \textbf{Hecke algebra} \(\mathcal{H}_n := \mathcal{H}(G, G_n)\) is the vector space of all functions \(G_n\backslash G/G_n \to \mathbb{C}\). It is equipped with the convolution
	\[
		(\varphi * \psi)(G_ngG_n) := \frac{1}{|G/G_n|} \sum_{fG_n \in G/G_n} \varphi(G_nfG_n) \cdot \psi(G_nf^{-1}gG_n)
	\]
	for all \(\varphi, \psi \in \mathcal{H}_n\) and \(G_ngG_n \in G_n\backslash G/G_n\) and the involution
	\[
		\varphi^*(G_ngG_n) := \overline{\varphi(G_ng^{-1}G_n)}
	\]
	for all \(\varphi \in \mathcal{H}_n\) and \(G_ngG_n \in G_n\backslash G/G_n\).
\end{definition}

The convolution is associative and bilinear, the involution is antilinear, self-inverse and we have \((\varphi * \psi)^* = \psi^* * \varphi^*\) (\citep[Lemmata 2.14. and 2.17]{hall_hecke_1999}). The Hecke algebra has a unit element given by \(|G/G_n| \cdot 1_{\{G_neG_n\}}\). The convolution can alternatively be written as
\begin{align}\label{eq-Convolution_Finite_Dimensional}
	(\varphi * \psi)(G_ngG_n) &= \frac{1}{|G/C_n|} \sum_{fC_n \in G/C_n} \varphi(G_nfG_n) \cdot \psi(G_nf^{-1}gG_n)\\
	&= \int_{G/C_n} \varphi(G_nfG_n) \cdot \psi(G_nf^{-1}gG_n)\ \mathrm{d}\mu_n(fC_n).\nonumber
\end{align}
The first term differs from \cref{def-Hecke_Alegbra_Finite_Dimensional} by taking the sum over \(G/C_n\) instead of \(G/G_n\). This is equivalent to \cref{def-Hecke_Alegbra_Finite_Dimensional}, since every left coset \(fG_n\) contains exactly \(\frac{|G/C_n|}{|G/G_n|}\) left cosets \(f'C_n\). The second term in \eqref{eq-Convolution_Finite_Dimensional} is equivalent to the first term in \eqref{eq-Convolution_Finite_Dimensional} since \(\mu_n\) is the normalized counting measure on \(G/C_n\). Note that this is an integral over \(G/C_n\), but \(\varphi\) and \(\psi\) are functions on \(G_n\backslash G/G_n\), so we implicitly compose them with the map \(r_n: G/C_n \to G_n\backslash G/G_n\), we could also explicitly write \(\varphi(r_n(fC_n))\) and \(\psi(r_n(f^{-1}gG_n))\) instead. In light of this integral formula for the convolution we could also consider elements in the Hecke algebra \(\mathcal{H}(G, G_n)\) to be function on \(G/C_n\) that are \(G_n/C_n\) left and right invariant.

In the literature the convolution is defined without the normalizing factor \(\frac{1}{|G/G_n|}\) (see \citep[Formula (2.3)]{hall_hecke_1999}). While the addition of the factor does not change any of the algebraic properties, it allows us to rewrite the sum as an integral, which is useful when considering the measures \(\mu_n\) and \(\nu_n\). In particular this modification is necessary to extend the concept of Hecke algebras to pairs of separable infinite profinite groups and closed core-free subgroups in the next part of this section.

We can equip the Hecke algebra \(\mathcal{H}(G, G_n)\) with the \(\operatorname{L}^1(\nu_n)\)-norm, where \(\nu_n = r_{n*}\mu_n\) is the probability measure on \(G_n\backslash G/G_n\) from \eqref{eq-nu_n}. 
We prove below that \(\mathcal{H}(G, G_n)\) is a Banach-*-algebra when equipped with this norm. Similar statements can be found in \citep[Definition 1.9.]{palma_enveloping_2013} and \citep[Section 3.]{tzanev_hecke_2003}. In these references the convolution is defined without the normalizing factor \(\frac{1}{|G/G_n|}\) and the measure of some double coset \(G_ngG_n \in G_n\backslash G/G_n\) is given by the number \(|G_n(gG_n)|\) of left cosets in the \(G_n\)-orbit of \(G_n\) instead of by \(\nu_n(\{G_ngG_n\}) = \frac{|G_n(gG_n)|}{|G/G_n|}\). This rescaling by \(|G/G_n|\) in both the convolution and the measure ensures that the sub-multiplicativity still holds with our definitions.

\begin{proposition}\label{pro-H_n_Banach-*-algebra}
	The Hecke algebra \(\mathcal{H}_n\) equipped with the \(\operatorname{L}^1(\nu_n)\)-norm is a Banach-*-algebra with isometric involution, i.e.
	\[
		\|\varphi * \psi\|_{\operatorname{L}^1(\nu_n)} \leq \|\varphi\|_{\operatorname{L}^1(\nu_n)} \cdot \|\psi\|_{\operatorname{L}^1(\nu_n)} \quad\text{and}\quad \|\varphi^*\|_{\operatorname{L}^1(\nu_n)} = \|\varphi\|_{\operatorname{L}^1(\nu_n)}
	\]
	for all \(\varphi, \psi \in \mathcal{H}_n\).
\end{proposition}
\begin{proof*}
	The space \(\mathcal{H}_n\) is complete as it is finite dimensional. The other properties are straightforward calculation:
	\begin{align*}
		\|\varphi * \psi\|_{\operatorname{L}^1(\nu_n)} &= \int_{G_n\backslash G/G_n} |(\varphi * \psi)(G_ngG_n)|\ \mathrm{d}\nu_n(G_ngG_n)\\
		&= \int_{G/C_n} |(\varphi * \psi)(G_ngG_n)|\ \mathrm{d}\mu_n(gG_n),\\
		\intertext{where we used the substitution given by the measure preserving map \(r_n\). We continue by using the integral definition of the convolution in \eqref{eq-Convolution_Finite_Dimensional}:}
		&= \int_{G/C_n} \bigg| \int_{G/C_n} \varphi(G_nfG_n) \cdot \psi(G_nf^{-1}gG_n)\ \mathrm{d}\mu_n(fC_n) \bigg|\ \mathrm{d}\mu_n(gC_n)\\
		&\leq \int_{G/C_n} |\varphi(G_nfG_n)| \int_{G/C_n} |\psi(G_nf^{-1}gG_n)|\ \mathrm{d}\mu_n(gC_n)\ \mathrm{d}\mu_n(fC_n)\\
		\intertext{where in the last step we use Tonelli's theorem. Here we apply in the inner integral the measure preserving transformation \(gC_n \mapsto fgC_n\) to obtain}
		&= \int_{G/C_n} |\varphi(G_nfG_n)| \int_{G/C_n} |\psi(G_ngG_n)|\ \mathrm{d}\mu_n(gC_n)\ \mathrm{d}\mu_n(fC_n)\\
		&= \bigg( \int_{G/C_n} |\varphi(G_nfG_n)|\ \mathrm{d}\mu_n(fC_n) \bigg) \cdot \bigg( \int_{G/C_n} |\psi(G_ngG_n)|\ \mathrm{d}\mu_n(gC_n) \bigg)\\
		&= \|\varphi\|_{\operatorname{L}^1(\nu_n)} \cdot \|\psi\|_{\operatorname{L}^1(\nu_n)},
	\end{align*}
	where for the last step we use again the substitution with \(r_n\).
	
	Since \(\mu_n\) is the Haar measure on the compact group \(G/C_n\) the inversion \(gC_n \mapsto g^{-1}C_n\) is measure preserving. Therefore we obtain by rewriting the integral over \(G_n\backslash G/G_n\) as an integral over \(G/C_n\) that
	\begin{alignat*}{3}
		\|\varphi^*\|_{\operatorname{L}^1(\nu_n)} &= \int_{G/C_n} |\overline{\varphi(G_ng^{-1}G_n)}|\ \mathrm{d}\mu_n(gC_n)\\
		&= \int_{G/C_n} |\varphi(G_ngG_n)|\ \mathrm{d}\mu_n(gC_n) &= \|\varphi\|_{\operatorname{L}^1(\nu_n)}.
	\end{alignat*}
\end{proof*}

In this proposition we see that we need the rescaling by \(|G/G_n|\) in the definition of the convolution when using the normalized measure \(\mu_n\), to ensure that \(\mathcal{H}(G, G_n)\) is a Banach algebra. The use of the normalized measure \(\nu_n\) was also relevant in \cref{pro-Maps_Measure_Preserving} to ensure that the inclusion map \(q_{n, m}: G_m\backslash G/G_m \to G_n\backslash G/G_n\) for \(m \geq n \geq 0\) is measure preserving. We use this fact now, since it implies that the Koopman operator \(\mathcal{Q}_{n, m}: \mathcal{H}_n \to \mathcal{H}_m\) defined by
\[
	\mathcal{Q}_{n, m} \varphi := \varphi \circ q_{n, m}
\]
for all \(\varphi \in \mathcal{H}_n\) is an \(\operatorname{L}^1\)-isometry (see \citep[Chapter 7]{eisner_operator_2015}). Also \(\mathcal{Q}_{n, m}\) is a *-algebra homomorphism.

\begin{proposition}\label{pro-Q_n_m_*-Algebra_Homomorphism}
	For \(m \geq n \geq 0\) we have that \(\mathcal{Q}_{n, m}\) is an isometric *-algebra homomorphism, i.e.
	\[
		\mathcal{Q}_{n, m}(\varphi * \psi) = \mathcal{Q}_{n, m}\varphi * \mathcal{Q}_{n, m}\psi \quad\text{and}\quad (\mathcal{Q}_{n, m}\varphi)^* = \mathcal{Q}_{n, m}\varphi^*
	\]
	for all \(\varphi, \psi \in \mathcal{H}_n\). Furthermore we have
	\[
		\mathcal{Q}_{n, k} = \mathcal{Q}_{m, k} \mathcal{Q}_{n, m}
	\]
	for all \(k \geq m \geq n \geq 0\).
\end{proposition}

To summarize the reasons for our modification from the definitions in the literature: we equip the spaces \(G_n\backslash G/G_n\) with probability measures to ensure that the *-algebra homomorphisms \(\mathcal{Q}_{n, m}\) between the Hecke algebras \(\mathcal{H}_n\) and \(\mathcal{H}_m\) are isometries. As a consequence we have to add the normalization factor \(\frac{1}{|G/G_n|}\) in the definition of the convolution, so that that the Hecke algebra \(\mathcal{H}_n\) is a Banach-*-algebra when equipped with the \(\operatorname{L}^1\)-norm.

\begin{proof*}[\cref{pro-Q_n_m_*-Algebra_Homomorphism}]
	The fact that \(\mathcal{Q}_{n, m}\) is isometric was already discussed. The equality \(\mathcal{Q}_{n, k} = \mathcal{Q}_{m, k} \mathcal{Q}_{n, m}\) follows directly from \(q_{n, k} = q_{n, m} \circ q_{m, k}\). For the convolution we calculate for \(m \geq n \geq 0\) and all \(\varphi, \psi \in \mathcal{H}_n\) and \(G_mgG_m \in G_m\backslash G/G_m\) that
	\begin{align*}
		\mathcal{Q}_{n, m}(\varphi * \psi)(G_mgG_m) &= (\varphi * \psi)(q_{n, m}(G_mgG_m))\\
		&= (\varphi * \psi)(G_ngG_n)\\
		&= \int_{G/C_n} \varphi(G_nfG_n) \cdot \psi(G_nf^{-1}gG_n)\ \mathrm{d}\mu_n(fC_n).\\
		\intertext{Now we transform the integral using the measure preserving map \(fC_m \mapsto fC_n\) (see \cref{pro-Maps_Measure_Preserving}). Note that we have \(C_m \leq C_n\leq G_n\).}
		&= \int_{G/C_m} \varphi(G_nfG_n) \cdot \psi(G_nf^{-1}gG_n)\ \mathrm{d}\mu_m(fC_m)\\
		&= \int_{G/C_m} \mathcal{Q}_{n, m}\varphi(G_mfG_m) \cdot \mathcal{Q}_{n, m}\psi(G_mf^{-1}gG_m)\ \mathrm{d}\mu_m(fC_m)\\
		&= (\mathcal{Q}_{n, m}\varphi * \mathcal{Q}_{n, m}\psi)(G_mgG_m).
	\end{align*}
	For the involution we calculate for all \(G_mgG_m \in G_m\backslash G/G_m\) that
	\begin{alignat*}{3}
		(\mathcal{Q}_{n, m}\varphi)^*(G_mgG_m) &= \overline{\mathcal{Q}_{n, m}\varphi(G_mg^{-1}G_m)}\\
		&= \overline{\varphi(G_ng^{-1}G_n)}\\
		&= \varphi^*(G_ngG_n) &= \mathcal{Q}_{n, m}\varphi^*(G_mgG_m).
	\end{alignat*}
\end{proof*}

\subsection{The Infinite Dimensional Hecke Algebra}\label{sub-Hecke_Algebra_G,D}

Now we define infinite dimensional Hecke algebras. We first define the Hecke algebra \(\mathcal{H}(\widehat{G}, D)\) in \cref{the-*-Algebra_Structure_H} for pairs \((\widehat{G}, D)\) consisting of the Ellis group \(\widehat{G}\) and the discriminant group \(D\) of a minimal equicontinuous Cantor action \((X, G, \Phi)\) with group chain \((G_n)_{n \in \mathbb{N}_0}\). For the construction we utilize the finite dimensional Hecke algebras \(\mathcal{H}_n = \mathcal{H}_n(G, G_n)\) introduced in \cref{sub-Finite_Dimensional_Hecke_algebras}. In \cref{the-*-Algebra_Structure_H} we will immediately see that \(\mathcal{H}(\widehat{G}, D)\) does not depend on the choice of group chain. In \cref{cor-Hecke_Algebra_Cantor_Action} we will see that it also does not, up to isomorphism, depend on the choice of base point.

For general pairs \((\widehat{G}, D)\) of a separable profinite group \(\widehat{G}\) with closed core-free subgroup \(D\) we have discussed in the introduction, that they can be identified with the Ellis group and the discriminant group of a minimal equicontinuous Cantor action of some countable group \(G\). This action is not unique, however, we see as a consequence of \cref{the-Isomorphic_Hecke_Pairs}, that different choices of minimal equicontinuous Cantor actions yield isometric isomorphic Hecke algebras. This extends the concept of the Hecke algebra to general pairs \((\widehat{G}, D)\).

We start with the case that \(\widehat{G}\) is the Ellis group and \(D\) the discriminant group of a minimal equicontinuous Cantor action \((X, G, \Phi)\). The Hecke algebra \(\mathcal{H} := \mathcal{H}(\widehat{G}, D)\) consists of \(\operatorname{L}^1(\nu)\)-functions \(\varphi: D\backslash\widehat{G}/D \to \mathbb{C}\). We define the Koopman operator \(\mathcal{Q}_n: \mathcal{H}_n \to \mathcal{H}\) by
\[
	\mathcal{Q}_n\varphi := \varphi \circ q_n
\]
for all \(\varphi \in \mathcal{H}_n\). Since \(q_n\) is measure preserving, it follows that \(\mathcal{Q}_n\) is an isometry. From the commuting diagram \eqref{dia-Inverse_Limit_Double_Cosets} we obtain for \(m \geq n \geq 0\) the following commuting diagram of Banach spaces:
\begin{equation}\label{dia-Hecke_Algebras}
	\begin{tikzcd}
		&\mathcal{H}\\
		\mathcal{H}_m \arrow[ru, hook, "\mathcal{Q}_m"]\\
		\mathcal{H}_n \arrow[u, hook, "\mathcal{Q}_{n, m}"] \arrow[ruu, bend right, hook, "\mathcal{Q}_n"]
	\end{tikzcd}
\end{equation}
where all the arrows are isometries. We want to define the *-algebra structure on \(\mathcal{H}\), such that the above diagram becomes a diagram of Banach-*-algebras. In order to do this we first require the following lemma, stating that we can approximate every function in \(\mathcal{H}\) by functions in the finite dimensional subspaces \(\operatorname{ran}(\mathcal{Q}_n)\), where \(\operatorname{ran}(\mathcal{Q}_n)\) denotes the range of \(\mathcal{Q}_n\).

\begin{lemma}\label{lem-Images_Q_n_Dense}
	The subspace \(\bigcup_{n \in \mathbb{N}_0} \operatorname{ran}(\mathcal{Q}_n)\) is dense in \(\mathcal{H}\).
\end{lemma}
\begin{proof*}
	The family \(\{q_n^{-1}(U_n): U_n \subseteq G_n\backslash G/G_n\}\) is clearly closed under intersections and taking the symmetric difference for all \(n \in \mathbb{N}_0\). Thus the same is true for
	\[
		\mathcal{D} := \bigcup_{n \in \mathbb{N}_0} \{q_n^{-1}(U_n): U_n \subseteq G_n\backslash G/G_n\},
	\]
	where we use that the sets in the union are increasing. Since \(D\backslash\widehat{G}/D \cong \varprojlim_{n \in \mathbb{N}_0} G_n\backslash G/G_n\) (\cref{pro-Inverse_Limit_Double_Cosets}) we have that \(\mathcal{D}\) is a basis of the topology of \(D\backslash\widehat{G}/D\). Thus, since \(\mathcal{D}\) is countable, every open set in \(D\backslash\widehat{G}/D\) can be written as a countable union of sets in \(\mathcal{D}\), so the \(\sigma\)-algebra generated by \(\mathcal{D}\) contains all open sets in \(D\backslash\widehat{G}/D\), hence all Borel sets. Therefore \(\mathcal{D}\) generates the Borel \(\sigma\)-algebra of \(D\backslash\widehat{G}/D\) and by \citep[Satz VI.2.28]{elstrodt_mas-_2005} we obtain that
	\[
		\operatorname{span}_\mathbb{C}\{1_U: U \in \mathcal{D}\}
	\]
	is dense in \(\operatorname{L}^1(\nu)\). Now for \(n \in \mathbb{N}_0\) and \(U_n \subseteq G_n\backslash G/G_n\) we have
	\[
		1_{q_n^{-1}(U_n)} = 1_{U_n} \circ q_n = \mathcal{Q}_n 1_{U_n} \in \operatorname{ran}(\mathcal{Q}_n).
	\]
	Thus we have
	\[
		\operatorname{span}_\mathbb{C}\{1_U: U \in \mathcal{D}\} \subseteq \bigcup_{n \in \mathbb{N}_0} \operatorname{ran}(\mathcal{Q}_n),
	\]
	so the latter is also dense in \(\operatorname{L}^1(\nu)\).
\end{proof*}

Now we can continue with the definition of the *-algebra structure of \(\mathcal{H}\). In \eqref{eq-Convolution_Finite_Dimensional} we rewrote the convolution for the finite dimensional Hecke algebra \(\mathcal{H}_n\) as an integral. This definition can be directly extended to \(\mathcal{H}\). The fact that we construct this *-algebra structure through the finite dimensional Hecke algebras allows us to skip lengthy calculations for the algebraic properties such as the associativity of the convolution, however all these properties can be shown by direct computation, without the use of the finite dimensional approximation.

We remark that the Hecke algebra consists of \(\operatorname{L}^1\)-functions on \(D\backslash\widehat{G}/D\), and the formulae for the convolution and involution in the theorem below only depend on the pair \((\widehat{G}, D)\). Thus the Hecke algebra is independent of the choice of group chain \((G_n)_{n \in \mathbb{N}_0}\).

\begin{theorem}\label{the-*-Algebra_Structure_H}
	There exists a unique Banach-*-algebra structure on \(\mathcal{H} = \operatorname{L}^1(\nu)\) such that \(\mathcal{Q}_n\) is a *-algebra homomorphism for all \(n \in \mathbb{N}_0\). Its convolution is given by
	\[
		(\varphi * \psi)(D\hat{g}D) := \int_{\widehat{G}} \varphi(D\hat{f}D) \cdot \psi(D\hat{f}^{-1}\hat{g}D)\ \mathrm{d}\mu(\hat{f})
	\]
	for \(\varphi, \psi \in \mathcal{H}\) and \(\nu\)-almost all \(D\hat{g}D \in D\backslash\widehat{G}/D\) and its involution is given by
	\[
		\varphi^*(D\hat{g}D) := \overline{\varphi(D\hat{g}^{-1}D)}
	\]
	 for \(\varphi \in \mathcal{H}\) and \(\nu\)-almost all \(D\hat{g}D \in D\backslash\widehat{G}/D\). The involution is isometric.
\end{theorem}
\begin{proof*}
	We know from \cref{pro-H_n_Banach-*-algebra} that the \(\mathcal{H}_n\) are Banach-*-algebras with isometric involution. The operators \(\mathcal{Q}_n: \mathcal{H}_n \to \mathcal{H}\) are linear isometries. For \(n \in \mathbb{N}_0\) we define a convolution and involution on \(\operatorname{ran}(\mathcal{Q}_n)\) by
	\begin{equation}\label{eq-*-Algebra_Structure_ranQ_n}
		\mathcal{Q}_n\varphi_n * \mathcal{Q}_n\psi_n := \mathcal{Q}_n(\varphi_n * \psi_n) \quad\text{and}\quad (\mathcal{Q}_n\varphi_n)^* := \mathcal{Q}_n\varphi_n^*
	\end{equation}
	for all \(\varphi_n, \psi_n \in \mathcal{H}_n\). This is well-defined since \(\mathcal{Q}_n\) is injective. Since \(\mathcal{Q}_n\) is an isometry this turns \(\operatorname{ran}(\mathcal{Q}_n)\), equipped with the \(\operatorname{L}^1(\nu)\)-norm, into a Banach-*-algebra with isometric involution. For \(m \geq n \geq 0\) we have \(\operatorname{ran}(\mathcal{Q}_n) \subseteq \operatorname{ran}(\mathcal{Q}_m)\). Then for \(\varphi_n, \psi_n \in \mathcal{H}_n\) we obtain, using that \(\mathcal{Q}_{n, m}\) is a *-algebra homomorphism by \cref{pro-Q_n_m_*-Algebra_Homomorphism}, that
	\[
		\mathcal{Q}_n(\varphi_n * \psi_n) \overset{\eqref{dia-Hecke_Algebras}}{=} \mathcal{Q}_m\mathcal{Q}_{n, m}(\varphi_n * \psi_n) = \mathcal{Q}_m(\mathcal{Q}_{n, m}\varphi_n * \mathcal{Q}_{n, m}\psi_n),
	\]
	so the definitions of the convolution on \(\operatorname{ran}(\mathcal{Q}_n)\) and \(\operatorname{ran}(\mathcal{Q}_m)\) from \eqref{eq-*-Algebra_Structure_ranQ_n} agree on \(\operatorname{ran}(\mathcal{Q}_n)\). The same is true for the involution since
	\[
		\mathcal{Q}_n\varphi_n^* = \mathcal{Q}_m\mathcal{Q}_{n, m}\varphi_n^* = \mathcal{Q}_m(\mathcal{Q}_{n, m}\varphi_n)^*.
	\]
	Thus the definition in \eqref{eq-*-Algebra_Structure_ranQ_n} allows us to define a *-algebra structure on
	\[
		\mathcal{A} := \bigcup_{n \in \mathbb{N}_0} \operatorname{ran}(\mathcal{Q}_n).
	\]
	Since all the \(\operatorname{ran}(\mathcal{Q}_n)\) are Banach-*-algebras and equipped with the same norm, it follows that \(\mathcal{A}\), equipped with the \(\operatorname{L}^1(\nu)\)-norm, is a normed algebra, i.e. the norm is sub-multiplicative, and the involution on \(\mathcal{A}\) is isometric.
	
	This allows us to uniquely extend the *-algebra structure of \(\mathcal{A}\) to its completion, which by \cref{lem-Images_Q_n_Dense} can be identified with \(\operatorname{L}^1(\nu) = \mathcal{H}\), turning \(\mathcal{H}\) into a Banach-*-algebra (see \citep[Theorem II.2.III]{najmark_normed_1964}). More precisely for all \(\varphi, \psi \in \mathcal{H}\) we find sequences \((\varphi_i)_{i \in \mathbb{N}}\) and \((\psi_i)_{i \in \mathbb{N}_0}\) in \(\mathcal{A}\) that approximate \(\varphi\) and \(\psi\). Then \((\varphi_i * \psi_i)_{i \in \mathbb{N}_0}\) is a Cauchy sequence by a standard computation, so we can defined
	\[
		\varphi * \psi := \lim_{i\to\infty} \varphi_i * \psi_i.
	\]
	Another standard computation proves that this is independent of the choice of approximations \((\varphi_i)_{i \in \mathbb{N}_0}\) and \((\psi_i)_{i \in \mathbb{N}_0}\). The convolution as defined here is clearly bilinear and the \(\operatorname{L}^1(\nu)\)-norm is sub-multiplicative, so with it \(\mathcal{H}\) becomes indeed a Banach algebra. This extension of the convolution is unique, since the convolution must be continuous and \(\mathcal{A}\) is dense in \(\mathcal{H}\). Regarding the involution we note that is an isometric, hence continuous, antilinear map on \(\mathcal{A}\). Therefore the involution can be uniquely extended to \(\mathcal{H}\). This makes \(\mathcal{H}\) a Banach-*-algebra with isometric involution.
	
	This *-algebra structure is the unique one that turns all the \(\mathcal{Q}_n\) into *-algebra homomorphisms, since any such structure must fulfill equation \eqref{eq-*-Algebra_Structure_ranQ_n} and the extension from \(\mathcal{A}\) to \(\mathcal{H}\) is unique.
	
	At last we need to prove that the convolution and involution are given by the formulae in the theorem. First we prove this on \(\bigcup_{n \in \mathbb{N}_0} \operatorname{ran}(\mathcal{Q}_n)\). Fix \(n \in \mathbb{N}_0\) and \(\varphi_n, \psi_n \in \mathcal{H}_n\). Then we calculate for \(\nu\)-almost all \(D\hat{g}D \in D\backslash\widehat{G}/D\) with \(\hat{g} = (g_kC_k)_{k \in \mathbb{N}_0}\) that
	\begin{align*}
		(\mathcal{Q}_n\varphi_n * \mathcal{Q}_n\psi_n)(D\hat{g}D) &\overset{\eqref{eq-*-Algebra_Structure_ranQ_n}}{=} \mathcal{Q}_n(\varphi_n * \psi_n)(D\hat{g}D)\\
		&= (\varphi_n * \psi_n)(G_ng_nG_n)\\
		&= \int_{G/C_n} \varphi_n(G_nfG_n) \cdot \psi_n(G_nf^{-1}g_nG_n)\ \mathrm{d}\mu_n(fC_n).\\
		\intertext{Next we use the measure preserving transformation \(p_n: \widehat{G} \to G/C_n\), which maps \((f_kC_k)_{k \in \mathbb{N}_0}\) to \(f_nC_n\) (see \cref{pro-Maps_Measure_Preserving}), to obtain}
		&= \int_{\widehat{G}} \varphi_n(G_nf_nG_n) \cdot \psi_n(G_nf_n^{-1}g_nG_n)\ \mathrm{d}\mu((f_kC_k)_{k \in \mathbb{N}_0})\\
		&= \int_{\widehat{G}} \mathcal{Q}_n\varphi_n(D\hat{f}D) \cdot \mathcal{Q}_n\psi_n(D\hat{f}^{-1}\hat{g}D)\ \mathrm{d}\mu(\hat{f}).
	\end{align*}
	For the involution we calculate for \(\nu\)-almost all \(D\hat{g}D \in D\backslash\widehat{G}/D\) that
	\begin{alignat*}{3}
		(\mathcal{Q}_n\varphi_n)^*(D\hat{g}D) &\overset{\eqref{eq-*-Algebra_Structure_ranQ_n}}{=} \mathcal{Q}_n\varphi_n^*(D\hat{g}D)\\
		&= \varphi_n^*(G_ng_nG_n)\\
		&= \overline{\varphi(G_ng_n^{-1}G_n)} &= \overline{\mathcal{Q}_n\varphi_n(D\hat{g}^{-1}D)}.
	\end{alignat*}
	Thus the formulae from the theorem hold on the dense subspace \(\bigcup_{n \in \mathbb{N}_0} \operatorname{ran}(\mathcal{Q}_n)\). The convolution and involution on \(\mathcal{H}\) are the unique continuous extension of the convolution and involution on \(\bigcup_{n \in \mathbb{N}_0} \operatorname{ran}(\mathcal{Q}_n)\). Thus it suffices if we prove that
	\[
		c(\varphi, \psi)(D\hat{g}D) := \int_{\widehat{G}} \varphi(D\hat{f}D) \cdot \psi(D\hat{f}^{-1}\hat{g}D)\ \mathrm{d}\mu(\hat{f})
	\]
	defining a map \(c: \mathcal{H} \times \mathcal{H} \to \mathcal{H}\) as well as
	\[
		i(\varphi):= D\hat{g}D \mapsto \overline{\varphi(D\hat{g}^{-1}D)}
	\]
	defining a map \(i: \mathcal{H} \to \mathcal{H}\) are continuous.
	
	For the map \(c\) we also have to show that it is well-defined. We see in a moment that the integral exists, and first show that the definition is in fact independent of the choice of a coset representative \(\hat{g}\). Thus fix \(\hat{d} \in D\). It is clear that the representative \(\hat{g}\hat{d}\) of \(D\hat{g}D\) yields the same value. For the representative \(\hat{d}\hat{g}\) we obtain
	\[
		\int_{\widehat{G}} \varphi(D\hat{f}D) \cdot \psi(D\hat{f}^{-1}\hat{d}\hat{g}D)\ \mathrm{d}\mu(\hat{f}) = \int_{\widehat{G}} \varphi(D\hat{f}D) \cdot \psi(D\hat{f}^{-1}\hat{g}D)\ \mathrm{d}\mu(\hat{f})
	\]
	by using the measure preserving transformation \(\hat{f} \mapsto \hat{d}\hat{f}\). Thus \(c\) is in fact well defined as long as the integral exists. For the existence of the integral we calculate
	\begin{align*}
		&\int_{\widehat{G} \times \widehat{G}} |\varphi(D\hat{f}D) \cdot \psi(D\hat{f}^{-1}\hat{g}D)|\ \mathrm{d}(\mu \otimes \mu)(\hat{g}, \hat{f})\\
		= &\int_{\widehat{G}} |\varphi(D\hat{f}D)| \int_{\widehat{G}} |\psi(D\hat{f}^{-1}\hat{g}D)|\ \mathrm{d}\mu(\hat{g})\ \mathrm{d}\mu(\hat{f}),\\
		\intertext{where we used Tonelli's theorem. We continue with the substitution \(\hat{g} \mapsto \hat{f}\hat{g}\) in the inner integral}
		= &\int_{\widehat{G}} |\varphi(D\hat{f}D)| \int_{\widehat{G}} |\psi(D\hat{g}D)|\ \mathrm{d}\mu(\hat{g})\ \mathrm{d}\mu(\hat{f})\\
		= &\bigg( \int_{D\backslash\widehat{G}/D} |\varphi(D\hat{f}D)|\ \mathrm{d}\nu(D\hat{f}D) \bigg) \cdot \bigg( \int_{D\backslash\widehat{G}/D} |\psi(D\hat{g}D)|\ \mathrm{d}\nu(D\hat{g}D) \bigg)\\
		= & \|\varphi\|_{\operatorname{L}^1(\nu)} \cdot \|\psi\|_{\operatorname{L}^1(\nu)},
	\end{align*}
	where in the second to last step we use the substitution with the measure preserving map \(r: \widehat{G} \to D\backslash\widehat{G}/D\) (see \cref{pro-Maps_Measure_Preserving}) in both integrals. Thus Fubini's theorem implies that \(\hat{f} \mapsto \varphi(D\hat{f}D) \cdot \psi(D\hat{f}^{-1}\hat{g}D)\) is \(\mu\)-integrable for \(\mu\)-almost all \(\hat{g} \in \widehat{G}\). Thus \(c\) is in fact well-defined for \(\nu\)-almost all \(D\hat{g}D \in D\backslash\widehat{G}/D\). Further it follows from the above calculation that
	\begin{alignat*}{3}
		\|c(\varphi, \psi)\|_{\operatorname{L}^1(\nu)} &= \int_{\widehat{G}} \bigg| \int_{\widehat{G}} \varphi(D\hat{f}D) \cdot \psi(D\hat{f}^{-1}\hat{g}D)\ \mathrm{d}\mu(\hat{f}) \bigg|\ \mathrm{d}\mu(\hat{g})\\
		&\leq \int_{\widehat{G} \times \widehat{G}} |\varphi(D\hat{f}D) \cdot \psi(D\hat{f}^{-1}\hat{g}D)|\ \mathrm{d}(\mu \otimes \mu)(\hat{g}, \hat{f}) &= \|\varphi\|_{\operatorname{L}^1(\nu)} \cdot \|\psi\|_{\operatorname{L}^1(\nu)}.
	\end{alignat*}
	This shows that \(c\) is continuous.
	
	For the continuity of \(i\) we note that since \(\widehat{G}\) is compact we have that \(\hat{g} \mapsto \hat{g}^{-1}\) preserves \(\mu\). Using that \(\hat{g} \mapsto D\hat{g}D\) is measure preserving, we can rewrite the integral over \(D\backslash\widehat{G}/D\) as an integral over \(\widehat{G}\) to obtain
	\[
		\|i(\varphi)\|_{\operatorname{L}^1(\nu)} = \int_{\widehat{G}} |\overline{\varphi(D\hat{g}^{-1}D)}|\ \mathrm{d}\mu(\hat{g}) = \int_{\widehat{G}} |\varphi(D\hat{g}D)|\ \mathrm{d}\mu(\hat{g}) = \|\varphi\|_{\operatorname{L}^1(\nu)}.
	\]
	Hence also \(i\) is continuous, which concludes the proof.
\end{proof*}

In the introduction we have seen that any pair \((\widehat{G}, D)\) of a separable profinite group \(\widehat{G}\) and a closed core-free subgroup \(D\) can be identified with the Ellis group and the discriminant group of a minimal equicontinuous Cantor action. This identification is not unique. However, the following theorem implies that different choices of minimal equicontinuous Cantor actions for a given pair \((\widehat{G}, D)\) yield isomorphic Hecke algebras. The theorem states that if for two minimal equicontinuous Cantor actions \((X, G, \Phi)\) and \((Y, H, \Psi)\) the pairs \((\widehat{G}, D)\) and \((\widehat{H}, E)\), consisting of Ellis group and discriminant group, are isomorphic, then their Hecke algebras are isomorphic as well.

\begin{theorem}\label{the-Isomorphic_Hecke_Pairs}
	Let \(\widehat{G}\) and \(\widehat{H}\) be the Ellis groups of two minimal equicontinuous Cantor actions \((X, G, \Phi)\) and \((Y, H, \Psi)\) with discriminant groups \(D \leq \widehat{G}\) and \(E \leq \widehat{H}\). Assume that there exists a topological group isomorphism \(\theta: \widehat{H} \to \widehat{G}\), i.e. \(\theta\) is a group isomorphism and a homeomorphism, with \(\theta(E) = D\). Define the map \(\tilde{\theta}: E\backslash\widehat{H}/E \to D\backslash\widehat{G}/D\) by
	\[
		\tilde{\theta}(E\hat{h}E) := D\theta(\hat{h})D
	\]
	for all \(E\hat{h}E \in E\backslash\widehat{H}/E\). Then its Koopman operator \(\mathcal{K}_{\tilde{\theta}}: \mathcal{H}(\widehat{G}, D) \to \mathcal{H}(\widehat{H}, E)\) defined by
	\[
		\mathcal{K}_{\tilde{\theta}}\varphi := \varphi \circ \tilde{\theta}
	\]
	for all \(\varphi \in \mathcal{H}(\widehat{G}, D)\) is an isometric *-algebra isomorphism.
\end{theorem}
\begin{proof*}
	The fact that \(\tilde{\theta}\) is well defined follows from \(\theta(E) = D\) and \(\theta\) being a group isomorphism. This also implies that \(\tilde{\theta}\) is bijective. By definition we have
	\begin{equation}\label{eq-tilde_thete_r_H_equal_r_G_theta}
		\tilde{\theta} \circ r_H = r_G \circ \theta,
	\end{equation}
	where \(r_H: \widehat{H} \to E\backslash\widehat{H}/E\) and \(r_G: \widehat{G} \to D\backslash\widehat{G}/D\) are the surjections \(r_H(\hat{h}) = E\hat{h}E\) and \(r_G(\hat{g}) = D\hat{g}D\). This implies that \(\tilde{\theta}\) is continuous, since \(r_G \circ \theta\) is continuous and \(E\backslash\widehat{H}/E\) is equipped with the quotient topology arising from \(r_H\). Thus \(\tilde{\theta}\) is a continuous bijection between compact Hausdorff spaces, hence a homeomorphism. We denote by \(\mu_H\) and \(\mu_G\) the Haar probability measures on \(\widehat{H}\) and \(\widehat{G}\) respectively. Further let \(\nu_H := r_{H*}\mu_H\) and \(\nu_G := r_{G*}\mu_G\) be the corresponding measures on \(E\backslash\widehat{H}/E\) and \(D\backslash\widehat{G}/D\) respectively. Then \(\theta_*\mu_G\) is an invariant probability measure on \(\widehat{G}\), so by the uniqueness of the Haar measure we must have \(\mu_H = \theta_*\mu_G\), i.e. \(\theta\) is measure preserving. Thus, by \eqref{eq-tilde_thete_r_H_equal_r_G_theta} also \(\tilde{\theta}\) is measure preserving, since both \(r_G\) and \(r_H\) are. Therefore the Koopman operator \(\mathcal{K}_{\tilde{\theta}}\) defines an isometric isomorphism of Banach spaces by \citep[Proposition 7.12. and Corollary 7.21.]{eisner_operator_2015}.
	
	Now we prove that \(\mathcal{K}_{\tilde{\theta}}\) is a *-algebra homomorphism. We first show that it is compatible with the convolution. Fix \(\varphi, \psi \in \mathcal{H}(\widehat{G}, D)\) and calculate for \(\nu_H\)-almost all \(E\hat{h}E \in E\backslash\widehat{H}/E\):
	\begin{align*}
		\mathcal{K}_{\tilde{\theta}}(\varphi * \psi)(E\hat{h}E) &= (\varphi * \psi)(\tilde{\theta}(E\hat{h}E))\\
		&= (\varphi * \psi)(E\theta(\hat{h})E)\\
		&= \int_{\widehat{G}} \varphi(D\hat{f}D) \cdot \psi(D\hat{f}^{-1}\theta(\hat{h})D)\ \mathrm{d}\mu_G(\hat{f}).\\
		\intertext{Next we substitute \(\hat{f}\) by \(\theta(\hat{f}')\) using the measure preserving group isomorphism \(\theta: \widehat{H} \to \widehat{G}\) to obtain}
		&= \int_{\widehat{H}} \varphi(D\theta(\hat{f}')D) \cdot \psi(D\theta(\hat{f}')^{-1}\theta(\hat{h})D)\ \mathrm{d}\mu_H(\hat{f}')\\
		&= \int_{\widehat{H}} \varphi(D\theta(\hat{f}')D) \cdot \psi(D\theta(\hat{f}'^{-1}\hat{h})D)\ \mathrm{d}\mu_H(\hat{f}')\\
		&= \int_{\widehat{H}} \mathcal{K}_{\tilde{\theta}}\varphi(E\hat{f}'E) \cdot \mathcal{K}_{\tilde{\theta}}\psi(E\hat{f}'^{-1}\hat{h}E)\ \mathrm{d}\mu_H(\hat{f}')\\
		&= (\mathcal{K}_{\tilde{\theta}}\varphi * \mathcal{K}_{\tilde{\theta}}\psi)(E\hat{h}E).
	\end{align*}
	For the involution we calculate for all \(\varphi \in \mathcal{H}(\widehat{G}, D)\) and \(\nu_H\)-almost all \(E\hat{h}E \in E\backslash\widehat{H}/E\):
	\begin{alignat*}{3}
		(\mathcal{K}_{\tilde{\theta}}\varphi)^*(E\hat{h}E) &= \overline{\mathcal{K}_{\tilde{\theta}}\varphi(E\hat{h}^{-1}E)}\\
		&= \overline{\varphi(D\theta(\hat{h})^{-1}D)}\\
		&= \varphi^*(D\theta(\hat{h})D) &= \mathcal{K}_{\tilde{\theta}}\varphi^*(E\hat{h}E).
	\end{alignat*}
	This concludes the proof of the theorem.
\end{proof*}

This theorem can also be used for minimal equicontinuous Cantor actions. Let \((X, G, \Phi)\) be a minimal equicontinuous Cantor action with Ellis group \(\widehat{G}\) and let \(D_x\) and \(D_y\) be the discriminant groups corresponding to the base points \(x, y \in X\). Then we have the following lemma.

\begin{lemma}\label{lem-G_U_G_V_isomorphic}
	There exists a topological group isomorphism \(\theta: \widehat{G} \to \widehat{G}\) such that \(\theta(D_y) = D_x\).
\end{lemma}
\begin{proof*}
	Since \(\Phi\) is assumed to be a minimal action, it follows that \(\widehat{G}\) acts transitively on \(X\), so there exist a \(\hat{g} \in \widehat{G}\) such that \(\hat{g}x = y\). Thus \(\theta: \widehat{G} \to \widehat{G}\) defined by \(\theta(\hat{f}) := \hat{g}^{-1}\hat{f}\hat{g}\) is a homeomorphic group isomorphism with \(\theta_{\hat{g}}(D_y) = D_x\).
\end{proof*}

Combining this lemma with \cref{the-Isomorphic_Hecke_Pairs} yields the following corollary. It allows us to speak of the \textbf{Hecke algebra of the action \((X, G, \Phi)\)}, which we denote by \(\mathcal{H}(\Phi)\).

\begin{corollary}\label{cor-Hecke_Algebra_Cantor_Action}
	The Hecke algebras \(\mathcal{H}(\widehat{G}, D_x)\) and \(\mathcal{H}(\widehat{G}, D_y)\) are isometric isomorphic Banach-*-algebras.
\end{corollary}

%% file: ch-Hecke_Algebra_as_Invariant.tex
\section{The Hecke Algebra as an Invariant}\label{sec-Invariant}

In the previous section we have seen how to define the Hecke algebra for a minimal equicontinuous Cantor action. Now we consider the consequences of two minimal equicontinuous Cantor actions having isometric isomorphic Hecke algebras. In the first part we prove \cref{the-Isometric_Isomorphisms_Hecke_Algebras}, the main result of our paper. It states that any isometric *-algebra isomorphism between the Hecke algebras \(\mathcal{H}(\widehat{G}, D)\) and \(\mathcal{H}(\widehat{H}, E)\) of two separable infinite profinite groups \(\widehat{G}\) and \(\widehat{H}\) with closed core-free subgroups \(D \leq \widehat{G}\) and \(E \leq \widehat{H}\) arises as the Koopman operator of a homeomorphism \(\theta: E\backslash\widehat{H}/E \to D\backslash\widehat{G}/D\) combined with the multiplication by some character of \(\widehat{H}\), i.e. a continuous group homomorphism from \(\widehat{H}\) into the unit circle \(\mathbb{S}^1 \subseteq \mathbb{C}\). We have seen an example of such a *-algebra isomorphism in \cref{the-Isomorphic_Hecke_Pairs}.

In the second part of this section we apply this to minimal equicontinuous Cantor actions and  associated group chains. We prove in \cref{the-Group_Chains_Same_Index} that if two minimal equicontinuous Cantor actions \((X, G, \Phi)\) and \((Y, H, \Psi)\) have isometric isomorphic Hecke algebras and \((G_n)_{n \in \mathbb{N}_0}\) is a group chain associated to \((X, G, \Phi)\), then we can find a group chain \((H_n')_{n \in \mathbb{N}_0}\) associated to \((Y, H, \Psi)\) such that \(G_n\) and \(H_n'\) have the same index for all \(n \in \mathbb{N}_0\). An immediate consequence of this is that the Steinitz order of a Cantor action is invariant under isometric isomorphism of the Hecke algebras (\cref{cor-Steinitz_Order}). As a second consequence we give in \cref{the-Hecke_Algebra_Zig-zag_Map} a characterization of the Hecke algebras \(\mathcal{H}(\Phi)\) and \(\mathcal{H}(\Psi)\) being isometrically isomorphic via maps between their finite dimensional Hecke algebras \(\mathcal{H}(G, G_n)\) and \(\mathcal{H}(H, H_n)\).

\subsection{Isometric Isomorphisms between Hecke Algebras}

Let \(\widehat{G}\) and \(\widehat{H}\) be two separable infinite profinite groups with closed core-free subgroups \(D \leq \widehat{G}\) and \(E \leq \widehat{H}\). Denote by \(r_G: \widehat{G} \to D\backslash\widehat{G}/D\) and \(r_H: \widehat{H} \to E\backslash\widehat{H}/E\) the surjections given by \(r_G(\hat{g}) = D\hat{g}D\) and \(r_H(\hat{h}) = E\hat{h}E\). Further let \(\mu_G\) and \(\mu_H\) be the Haar measures on \(\widehat{G}\) and \(\widehat{H}\) respectively and denote by \(\nu_G = r_{G*}\mu_G\) and \(\nu_H = r_{H*}\mu_H\) the corresponding pushforward probability measures on \(D\backslash\widehat{G}/D\) and \(E\backslash\widehat{H}/E\). These maps and measures were discussed in more detail in \cref{sec-Cantor_Actions}.

\begin{theorem}\label{the-Isometric_Isomorphisms_Hecke_Algebras}
	Assume that \(\mathcal{T}: \mathcal{H}(\widehat{G}, D) \to \mathcal{H}(\widehat{H}, E)\) is an isometric *-algebra isomorphism. Then there exist a unique continuous group homomorphism \(\chi: \widehat{H} \to \mathbb{S}^1\) with \(E \leq \ker(\chi)\) and a unique measure preserving homeomorphism \(\theta: E\backslash\widehat{H}/E \to D\backslash\widehat{G}/D\) with \(\theta(EeE) = DeD\) such that
	\[
		\mathcal{T}\varphi(E\hat{h}E) = \chi(\hat{h}) \cdot \varphi(\theta(E\hat{h}E))
	\]
	for all \(\varphi \in \mathcal{H}(\widehat{G}, D)\) and \(\nu_H\)-almost all \(E\hat{h}E \in E\backslash\widehat{H}/E\). Moreover the Koopman operator \(\mathcal{T}': \mathcal{H}(\widehat{G}, D) \to \mathcal{H}(\widehat{H}, E)\) defined by
	\[
		\mathcal{T}'\varphi := \varphi \circ \theta
	\]
	for all \(\varphi \in \mathcal{H}(\widehat{G}, D)\) is an isometric *-algebra isomorphism as well.
\end{theorem}

We first outline the proof of this theorem and state multiple lemmata used in it. The proofs of the theorem and these lemmata then follow below. First we show in \cref{lem-T1_X=1_Y} that the map \(\hat{H} \mapsto \mathcal{T}1_{D\backslash\widehat{G}/D}(E\hat{h}E)\) is given by a continuous character \(\chi\) of \(\widehat{H}\). The continuity of \(\chi\) is due to \cref{the-Kleppner}. We then use this character to obtain a modified *-algebra isomorphism \(\mathcal{T}'\) which maps \(1_{D\backslash\widehat{G}/D}\) to \(1_{E\backslash\widehat{H}/E}\). Going onward from \cref{lem-T1_X=1_Y}, we prove for this modified *-algebra isomorphism that it is the Koopman operator of a measure preserving homeomorphism. In the end, this yields that the original *-algebra homomorphism \(\mathcal{T}\) is of the form as stated in \cref{the-Isometric_Isomorphisms_Hecke_Algebras} above.

\begin{lemma}\label{lem-T1_X=1_Y}
	Assume that \(\mathcal{T}: \mathcal{H}(\widehat{G}, D) \to \mathcal{H}(\widehat{H}, E)\) is an isometric *-algebra isomorphism. Then there exists a continuous group homomorphism \(\chi: \widehat{H} \to \mathbb{S}^1\) with \(E \leq \ker(\chi)\), such that \(\mathcal{T}1_{D\backslash\widehat{G}/D}(E\hat{h}E) = \chi(\hat{h})\) for \(\nu_H\)-almost all \(E\hat{h}E \in E\backslash\widehat{H}/E\). Then \(\mathcal{T}': \mathcal{H}(\widehat{G}, D) \to \mathcal{H}(\widehat{H}, E)\) defined by
	\[
		\mathcal{T}'\varphi(E\hat{h}E) := \overline{\chi(\hat{h})} \cdot \mathcal{T}\varphi(E\hat{h}E)
	\]
	for all \(\varphi \in \mathcal{H}(\widehat{G}, D)\) and \(\nu_H\)-almost all \(E\hat{h}E \in E\backslash\widehat{H}/E\) is an isometric *-algebra isomorphism, satisfying
	\[
		\mathcal{T}'1_{D\backslash\widehat{G}/D} = 1_{E\backslash\widehat{H}/E}.
	\]
\end{lemma}

The modified *-algebra isomorphism \(\mathcal{T}'\)  fulfills the conditions of \cref{cor-Banach-Lamperti}. Therefore \(\mathcal{T}'\) is the Koopman operator of some measure preserving map \(\theta': E\backslash\widehat{H}/E \to D\backslash\widehat{G}/D\). Our next goal is to prove that \(\theta'\) is continuous. For this we make use of approximate units and the fact that \(\mathcal{T}'\) maps approximate units to approximate units. The easiest way to construct an approximate unit in \(\mathcal{H}(\widehat{G}, D)\) is to fix a neighborhood basis \((\widetilde{U}_n)_{n \in \mathbb{N}}\) of the double coset of the identity \(DeD\), and take
\begin{equation}\label{eq-Approximate_Unit}
	\bigg( \frac{1}{\nu_G(\widetilde{U}_n)} \cdot 1_{\widetilde{U}_n} \bigg)_{n \in \mathbb{N}}.
\end{equation}
The lemma below gives a characterization of more general sequences \((\widetilde{U}_n)_{n \in \mathbb{N}}\), such that \eqref{eq-Approximate_Unit} is an approximate identity. This characterization in assertion (i) states that the proportion of \(\widetilde{U}_n\) that is contained in some neighborhood \(\widetilde{U}\) of the identity goes to one as \(n\) tends infinity. This can be interpreted as \((\widetilde{U}_n)_{n \in \mathbb{N}}\) being an approximate neighborhood basis of \(DeD\). It is in particular satisfied by every neighborhood basis of \(DeD\).

\begin{lemma}\label{lem-Approximate_Unit}
	Let \((\widetilde{U}_n)_{n \in \mathbb{N}}\) be a decreasing sequence of measurable sets in \(D\backslash\widehat{G}/D\) with positive measure. The following assertions are equivalent:
	\begin{enumerate}[(i)]
		\item for any measurable neighborhood \(\widetilde{U}\) of \(DeD\) we have
		\[
			\lim_{n \to \infty} \frac{\nu_G(\widetilde{U}_n \setminus \widetilde{U})}{\nu_G(\widetilde{U}_n)} = 0,
		\]
		\item for all \(\varphi \in \mathcal{H}(\widehat{G}, D)\) we have
		\[
			\lim_{n \to \infty} \frac{1}{\nu_G(\widetilde{U}_n)} \cdot (1_{\widetilde{U}_n} * \varphi) = \varphi
		\]
		with convergence in \(\operatorname{L}^1(\nu_G)\),
		\item for all continuous \(\varphi \in \mathcal{H}(\widehat{G}, D)\) we have
		\[
			\lim_{n \to \infty} \frac{1}{\nu_G(\widetilde{U}_n)} \cdot (1_{\widetilde{U}_n} * \varphi) = \varphi
		\]
		with uniform convergence.
	\end{enumerate}
\end{lemma}

The third assertion in the above lemma states that for an approximate unit we even have uniform convergence if we approximate continuous functions. Note that as long as the sets \(\widetilde{U}_n\) are clopen, it follows that \(1_{\widetilde{U}_n} * \varphi\) is continuous for all \(\varphi \in \mathcal{H}(\widehat{G}, D)\). Thus, in this case the uniform convergence holds if and only if \(\varphi\) is continuous. Also \(\theta'\) is continuous if and only if \(\varphi \circ \theta'\) is continuous for all continuous \(\varphi\). These two thoughts go into the lemma below, which gives a criterion for \(\theta'\) being almost everywhere equal to a continuous function.

\begin{lemma}\label{lem-Continuous_Map_from_Measurable_Map}
	Let \((V_n)_{n \in \mathbb{N}}\) be a decreasing sequence of clopen subgroups in \(\widehat{H}\) such that \(E = \bigcap_{n \in \mathbb{N}} V_n\) and set
	\[
		\widetilde{V}_n := E\backslash V_n/E \subseteq E\backslash\widehat{H}/E.
	\]
	Let \(\theta': E\backslash\widehat{H}/E \to D\backslash\widehat{G}/D\) be a measure preserving map such that for every continuous \(\varphi \in \mathcal{H}(\widehat{G}, D)\) we have
	\[
		\lim_{n \to \infty} \frac{1}{\nu_H(\widetilde{V}_n)} \cdot (1_{\widetilde{V}_n} * (\varphi \circ \theta')) = \varphi \circ \theta'
	\]
	uniformly outside of some measure-zero set. Then there exist a continuous map \(\theta: E\backslash\widehat{H}/E \to D\backslash\widehat{G}/D\) that is almost everywhere equivalent to \(\theta'\).
\end{lemma}

Now we get to the proofs. We prove the three lemmata first and then \cref{the-Isometric_Isomorphisms_Hecke_Algebras} itself.

\begin{proof*}[\cref{lem-T1_X=1_Y}]
	Define \(\chi' := \mathcal{T}1_{D\backslash\widehat{G}/D}\). Since \(\mathcal{T}\) is an isometry we have \(\|\chi'\|_{\operatorname{L}^1(\nu_H)} = 1\). From \cref{the-Banach-Lamperti} we know that \(\chi' \neq 0\) almost everywhere. Further, \(1_{D\backslash\widehat{G}/D} * 1_{D\backslash\widehat{G}/D} = 1_{D\backslash\widehat{G}/D}\) and \(1_{D\backslash\widehat{G}/D}^* = 1_{D\backslash\widehat{G}/D}\) together with \(\mathcal{T}\) being a *-algebra isomorphism implies
	\begin{equation}\label{eq-chi*chi_Equal_chi}
		\chi' * \chi' = \chi' \quad\text{and}\quad \chi'^* = \chi'.
	\end{equation}
	We now prove that \(|\chi'|\) is almost everywhere constant. For the sake of contradiction we assume that this is not the case, i.e. there exists some \(\varepsilon > 0\) such that, with \(M := \operatorname{ess\ sup} |\chi'|\), the set
	\[
		U := \{\hat{f} \in \widehat{H}: |\chi'(E\hat{f}E)| < M - \varepsilon\}
	\]
	has positive measure. Since \(|\chi'| \neq 0\) almost everywhere there is some \(\delta > 0\) such that
	\[
		V := \{\hat{f} \in \widehat{H}: |\chi'(E\hat{f}E)| \leq \delta\}
	\]
	fulfills \(\mu_H(V) < \mu_H(U)\). Note that we may have \(\mu_H(V) = 0\). Then we calculate for almost all \(\hat{h} \in \widehat{H}\):
	\begin{align*}
		|\chi'(E\hat{h}E)| &\overset{\eqref{eq-chi*chi_Equal_chi}}{=} \bigg| \int_{\widehat{H}} \chi'(EfE) \cdot \chi'(E\hat{f}^{-1}\hat{h}E)\ \mathrm{d}\mu_H(\hat{f}) \bigg|\\
		&\leq \int_{\widehat{H}} |\chi'(E\hat{f}E)| \cdot |\chi'(E\hat{f}^{-1}\hat{h}E)|\ \mathrm{d}\mu_H(\hat{f}).\\
		\intertext{Next we use that \(|\chi'(E\hat{f}^{-1}\hat{h}E)|\) is bounded by \(M\) almost everywhere and for \(\hat{f}^{-1}\hat{h} \in U\) it is even bounded by \(M - \varepsilon\), so}
		&\leq \int_{\widehat{H}} |\chi'(E\hat{f}E)| \cdot (M - \varepsilon \cdot 1_U(\hat{f}^{-1}\hat{h}))\ \mathrm{d}\mu_H(\hat{f})\\
		&= M \cdot \|\chi'\|_{\operatorname{L}^1(\nu_H)} - \varepsilon \int_{\widehat{H}} |\chi'(E\hat{f}E)| \cdot 1_U(\hat{f}^{-1}\hat{h})\ \mathrm{d}\mu_H(\hat{f})\\
		&= M - \varepsilon \int_{\hat{h}U^{-1}} |\chi'(E\hat{f}E)|\ \mathrm{d}\mu_H(\hat{f})\\
		&\leq M - \varepsilon \int_{\hat{h}U^{-1} \setminus V} |\chi'(E\hat{f}E)|\ \mathrm{d}\mu_H(f)\\
		&\leq M - \varepsilon\delta \cdot m_H(hU^{-1}\setminus V)
	\end{align*}
	since \(\hat{f} \mapsto |\chi'(E\hat{f}E)|\) is bounded from below by \(\delta\) outside of \(V\). From \(\mu_H(\hat{h}U^{-1}) = \mu_H(U^{-1}) = \mu_H(U)\) and \(\mu_H(U) > \mu_H(V)\) it follows that \(\mu_H(\hat{h}U^{-1} \setminus V) > 0\). Thus
	\[
		|\chi'(E\hat{h}E)| \leq M - \varepsilon\delta \cdot \mu_H(\hat{h}U^{-1}\setminus V) < M
	\]
	for almost all \(\hat{h} \in \widehat{H}\), which is a contradiction to \(M\) being the essential supremum of \(|\chi'|\). Therefore we have shown that \(|\chi'|\) is almost everywhere constant. Since \(\|\chi'\|_{\operatorname{L}^1(\nu_H)} = 1\), it follows that \(|\chi'| = 1\) almost everywhere.
	
	Now we show that \(\chi' \circ r_H: \widehat{H} \to \mathbb{S}^1\), i.e. the map \(\hat{h} \mapsto \chi'(E\hat{h}E)\), is almost everywhere equivalent to a continuous group homomorphism \(\chi: \widehat{H} \to \mathbb{S}^1\). For this we use \cref{the-Kleppner}. Thus we have to prove that \(\chi' \circ r_H\) is almost everywhere a group homomorphism. For \(\mu_H\)-almost all \(\hat{h} \in \widehat{H}\) we have
	\begin{alignat*}{3}
		1 &= |\chi'(E\hat{h}E)|\\
		&\overset{\eqref{eq-chi*chi_Equal_chi}}{=} |(\chi' * \chi')(E\hat{h}E)|\\
		&= \bigg| \int_{\widehat{H}} \chi'(E\hat{f}E) \cdot \chi'(E\hat{f}^{-1}\hat{h}E)\ \mathrm{d}\mu_H(\hat{f}) \bigg|\\
		&\leq \int_{\widehat{H}} |\chi'(E\hat{f}E)| \cdot |\chi'(E\hat{f}^{-1}\hat{h}E)|\ \mathrm{d}\mu_H(\hat{f}) &\ = 1.
	\end{alignat*}
	Thus it follows that
	\[
		\bigg| \int_{\widehat{H}} \chi'(E\hat{f}E) \cdot \chi'(E\hat{f}^{-1}\hat{h}E)\ \mathrm{d}\mu_H(\hat{f}) \bigg| = \int_{\widehat{H}} |\chi'(E\hat{f}E) \cdot \chi'(E\hat{f}^{-1}\hat{h}E)|\ \mathrm{d}\mu_H(\hat{f}) = 1
	\]
	must hold. Therefore \cref{lem-Constant_Absolute_Value_Triangle_Equality} implies that for almost all \(\hat{h} \in \widehat{H}\) the map \(\hat{f} \mapsto \chi'(E\hat{f}E) \cdot \chi'(E\hat{f}^{-1}\hat{h}E)\) is almost everywhere equal to some constant, with the constant depending on \(\hat{h}\). Thus the integrand on the left-hand side of
	\[
		\chi'(E\hat{h}E) \overset{\eqref{eq-chi*chi_Equal_chi}}{=} (\chi' * \chi')(E\hat{h}E) = \int_{\widehat{H}} \chi'(E\hat{f}E) \cdot \chi'(E\hat{f}^{-1}\hat{h}E)\ \mathrm{d}\mu_H(\hat{f})
	\]
	is almost everywhere independent of \(\hat{f}\) and \(\mu_H\) is a probability measure, so we have
	\[
		(\chi' \circ r_H)(\hat{h}) = \chi'(E\hat{h}E) = \chi'(E\hat{f}^{-1}E) \cdot \chi'(E\hat{f}\hat{h}E) = (\chi' \circ r_H)(\hat{f}^{-1}) \cdot (\chi' \circ r_H)(\hat{f}\hat{h})
	\]
	 for almost all \(\hat{h}, \hat{f} \in \widehat{H}\). Additionally we have with \(\chi'(E\hat{h}E) \in \mathbb{S}^1\) for almost all \(\hat{h} \in \widehat{H}\) that
	\[
		(\chi' \circ r_H)(\hat{h}^{-1}) = \chi'(E\hat{h}^{-1}E) = \overline{\chi'^*(E\hat{h}E)} \overset{\eqref{eq-chi*chi_Equal_chi}}{=} \overline{\chi'(E\hat{h}E)} = \chi'(E\hat{h}E)^{-1} = (\chi' \circ r_H)(\hat{h})^{-1}
	\]
	holds. Thus the map \(\chi' \circ r_H\) is by \cref{the-Kleppner} \(\mu\)-almost everywhere equivalent to a continuous group homomorphism \(\chi: \widehat{H} \to \mathbb{S}^1\).
	
	Next we prove that \(E \leq \ker(\chi)\), so fix \(\hat{d} \in E\). For almost all \(\hat{h} \in \widehat{H}\) we have \(\chi(\hat{h}) = \chi'(E\hat{h}E)\) and \(\chi(\hat{h}\hat{d}) = \chi'(E\hat{h}\hat{d}E)\), in particular every open set contains an element that fulfills both these equalities. Thus we can construct a sequence \((\hat{h}_n)_{n \in \mathbb{N}}\) in \(\widehat{H}\) that converges to the unit element, which has
	\[
		\chi(\hat{h}_n) = \chi'(E\hat{h}_nE) \quad\text{and}\quad \chi(\hat{h}_n\hat{d}) = \chi'(E\hat{h}_n\hat{d}E)
	\]
	for all \(n \in \mathbb{N}_0\). Then
	\[
		\chi(\hat{d}) = \lim_{n \to \infty} \chi(\hat{h}_n\hat{d}) = \lim_{n \to \infty} \chi'(E\hat{h}_n\hat{d}E) = \lim_{n \to \infty} \chi'(E\hat{h}_nE) = \lim_{n \to \infty} \chi(\hat{h}_n) = 1,
	\]
	so indeed \(E \leq \ker(\chi)\).
	
	At last we note that \(\overline{\chi}: \widehat{H} \to \mathbb{S}^1\), the pointwise complex conjugation of \(\chi\), is also a continuous group homomorphism with \(E \leq \ker(\overline{\chi})\). Thus \(\mathcal{T}'\) defined by
	\[
		\mathcal{T}'\varphi(E\hat{h}E) := \overline{\chi(\hat{h})} \cdot \mathcal{T}\varphi(E\hat{h}E)
	\]
	for all \(\varphi \in \mathcal{H}(\widehat{G}, D)\) and almost all \(E\hat{h}E \in E\backslash\widehat{H}/E\) is well defined. It is a straightforward calculation to check that \(\mathcal{T}'\) is an isometric *-algebra isomorphism. Also we have
	\[
		\mathcal{T}'1_{D\backslash\widehat{G}/D}(E\hat{h}E) = \overline{\chi(\hat{h})} \cdot \mathcal{T}1_{D\backslash\widehat{G}/D}(E\hat{h}E) = \overline{\chi(\hat{h})} \cdot \chi(\hat{h}) = 1,
	\]
	so \(\mathcal{T}'\) is as required.
\end{proof*}

\begin{proof*}[\cref{lem-Approximate_Unit}]
	We start by proving direction (i) implies (iii), so assume that
	\[
		\lim_{n \to \infty} \frac{\nu_G(\widetilde{U}_n \setminus \widetilde{U})}{\nu_G(\widetilde{U}_n)} = 0
	\]
	holds for any neighborhood \(\widetilde{U}\) of \(DeD\). Fix \(\varphi \in \mathcal{H}(\widehat{G}, D)\) continuous and \(\varepsilon > 0\). Since \(\widehat{G}\) is compact, we have that \(\hat{g} \mapsto \varphi(D\hat{g}^{-1}D)\) is uniformly continuous, so there is a neighborhood \(U'\) of \(e\) such that
	\begin{equation}\label{eq-Uniform_Continuity_varphi}
		|\varphi(D\hat{g}^{-1}D) - \varphi(D\hat{f}^{-1}D)| < \frac{\varepsilon}{2}
	\end{equation}
	whenever \(\hat{f}^{-1}\hat{g} \in U'\). Since \(\hat{g} \mapsto \varphi(D\hat{g}^{-1}D)\) is \(D\) left and right invariant, it even suffices if \(\hat{f}^{-1}\hat{g} \in U := DU'D\). Define
	\[
		\widetilde{U} := \{D\hat{g}D \in D\backslash\widehat{G}/D: \hat{g} \in U\},
	\]
	which is a neighborhood of \(DeD\). By definition \(r_G^{-1}(\tilde{U}) = U\). By (i) there is some \(N \in \mathbb{N}\) such that
	\begin{align}\label{eq_U_n_almost_contained_in_U}
		\frac{\nu_G(\widetilde{U}_n \setminus \widetilde{U})}{\nu_G(\widetilde{U}_n)} < \frac{\varepsilon}{4\|\varphi\|_\infty}
	\end{align}
	for all \(n \geq N\). Thus we have for all \(\hat{g} \in \widehat{G}\) and \(n \geq N\) that
	\begin{align*}
		&\bigg| \frac{1}{\nu_G(\widetilde{U}_n)} \cdot (1_{\widetilde{U}_n} * \varphi)(D\hat{g}D) - \varphi(D\hat{g}D) \bigg|\\
		= &\bigg| \frac{1}{\nu_G(\widetilde{U}_n)} \int_{\widehat{G}} 1_{\widetilde{U}_n}(D\hat{f}D) \cdot \varphi(D\hat{f}^{-1}\hat{g}D)\ \mathrm{d}\mu_G(\hat{f}) - \varphi(D\hat{g}D) \bigg|.
	\intertext{Next we use that \(D\hat{f}D \in \widetilde{U}_n\) if and only if \(\hat{f} \in r_G^{-1}(\widetilde{U}_n)\), so we continue with}
		= &\bigg| \frac{1}{\nu_G(\widetilde{U}_n)} \int_{r_G^{-1}(\widetilde{U}_n)} \varphi(D\hat{f}^{-1}\hat{g}D)\ \mathrm{d}\mu_G(\hat{f}) - \varphi(D\hat{g}D) \bigg|\\
		= &\bigg| \frac{1}{\nu_G(\widetilde{U}_n)} \int_{r_G^{-1}(\widetilde{U}_n)} \varphi(D\hat{f}^{-1}\hat{g}D) - \varphi(D\hat{g}D)\ \mathrm{d}\mu_G(\hat{f}) \bigg|\\
		\leq &\frac{1}{\nu_G(\widetilde{U}_n)} \int_{r_G^{-1}(\widetilde{U}_n)} \big|\varphi(D\hat{f}^{-1}\hat{g}D) - \varphi(D\hat{g}D)\big|\ \mathrm{d}\mu_G(\hat{f})\\
		\intertext{Now we use that by \eqref{eq-Uniform_Continuity_varphi} we have \(|\varphi(D\hat{f}^{-1}\hat{g}D) - \varphi(D\hat{g}D)| < \frac{\varepsilon}{2}\) whenever \(\hat{f} = (\hat{g}^{-1})^{-1}(\hat{f}^{-1}\hat{g})^{-1} \in U\) and \(|\varphi(D\hat{f}^{-1}\hat{g}D) - \varphi(D\hat{g}D)| \leq 2\|\varphi\|_\infty\) anywhere else, to obtain}
		\leq &\frac{1}{\nu_G(\widetilde{U}_n)} \int_{r_G^{-1}(\widetilde{U}_n)} \frac{\varepsilon}{2} \cdot 1_{U}(\hat{f}) + 2\|\varphi\|_\infty \cdot 1_{\widehat{G} \setminus U}(\hat{f})\ \mathrm{d}\mu_G(\hat{f}),\\
		= &\frac{\varepsilon}{2} \cdot \frac{\mu_G(r_G^{-1}(\widetilde{U}_n) \cap U)}{\nu_G(\widetilde{U}_n)} + 2\|\varphi\|_\infty \cdot \frac{\mu_G(r_G^{-1}(\widetilde{U}_n) \setminus U)}{\nu_G(\widetilde{U}_n)}.\\
		\intertext{Now we use that \(U = r_G^{-1}(\widetilde{U})\) and that \(r_G\) is measure preserving, such that}
		\leq &\frac{\varepsilon}{2} \cdot \frac{\nu_G(\widetilde{U}_n)}{\nu_G(\widetilde{U}_n)} + 2\|\varphi\|_\infty \cdot \frac{\nu_G(\widetilde{U}_n \setminus \widetilde{U})}{\nu_G(\widetilde{U}_n)}\\
		\overset{\eqref{eq_U_n_almost_contained_in_U}}{<} &\frac{\varepsilon}{2} + 2\|\varphi\|_\infty \cdot \frac{\varepsilon}{4\|\varphi\|_\infty}\\
		= &\varepsilon.
	\end{align*}
	This proves uniform convergence and thus assertion (iii) holds.
	
	Next we prove that (iii) implies (ii). Since \(\nu_G\) is a probability measure we have that uniform convergence implies \(\operatorname{L}^1\)-convergence. Also the continuous functions are dense in \(\mathcal{H}(\widehat{G}, D)\). Fix \(\varphi \in \mathcal{H}(\widehat{G}, D)\) and \(\varepsilon > 0\). Then there is a continuous function \(\varphi' \in \mathcal{H}(\widehat{G}, D)\) such that
	\begin{equation}\label{eq-varphi-varphi'}
		\|\varphi - \varphi'\|_{\operatorname{L}^1(\nu_G)} < \frac{\varepsilon}{3}
	\end{equation}
	and some \(N \in \mathbb{N}\) such that
	\begin{equation}\label{eq-Approximation_varphi'}
		\bigg\| \frac{1}{\nu_G(\widetilde{U}_n)} \cdot (1_{\widetilde{U}_n} * \varphi') - \varphi' \bigg\|_{\operatorname{L}^1(\nu_G)} < \frac{\varepsilon}{3}
	\end{equation}
	for all \(n \geq N\). Then we calculate
	\begin{align*}
		&\bigg\| \frac{1}{\nu_G(\widetilde{U}_n)} \cdot (1_{\widetilde{U}_n} * \varphi) - \frac{1}{\nu_G(\widetilde{U}_n)} \cdot (1_{\widetilde{U}_n} * \varphi') \bigg\|_{\operatorname{L}^1(\nu_G)}\\
		= &\frac{1}{\nu_G(\widetilde{U}_n)} \int_{D\backslash\widehat{G}/D} \big|(1_{\widetilde{U}_n} * \varphi)(D\hat{g}D) - (1_{\widetilde{U}_n} * \varphi')(D\hat{g}D)\big|\ \mathrm{d}\nu_G(D\hat{g}D)\\
		= &\frac{1}{\nu_G(\widetilde{U}_n)} \int_{\widehat{G}} \bigg| \int_{\widehat{G}} 1_{\widetilde{U}_n}(D\hat{f}D) \cdot \big( \varphi(D\hat{f}^{-1}\hat{g}D) - \varphi'(D\hat{f}^{-1}\hat{g}D) \big)\ \mathrm{d}\mu_G(\hat{f}) \bigg|\ \mathrm{d}\mu_G(\hat{g})\\
		\leq &\frac{1}{\nu_G(\widetilde{U}_n)} \int_{\widehat{G}} 1_{\widetilde{U}_n}(D\hat{f}D) \int_{\widehat{G}} \big|\varphi(D\hat{f}^{-1}\hat{g}D) - \varphi'(D\hat{f}^{-1}\hat{g}D)\big|\ \mathrm{d}\mu_G(\hat{g})\ \mathrm{d}\mu_G(\hat{f}).\\
		\intertext{In the inner integral we use the measure preserving transformation \(\hat{g} \mapsto \hat{f}\hat{g}\) such that}
		= &\frac{1}{\nu_G(\widetilde{U}_n)} \int_{\widehat{G}} 1_{\widetilde{U}_n}(D\hat{f}D) \int_{\widehat{G}} \big|\varphi(D\hat{g}D) - \varphi'(D\hat{g}D)\big|\ \mathrm{d}\mu_G(\hat{g})\ \mathrm{d}\mu_G(\hat{f})\\
		= &\frac{1}{\nu_G(\widetilde{U}_n)} \int_{\widehat{G}} 1_{\widetilde{U}_n}(D\hat{f}D)\ \mathrm{d}\mu_G(\hat{f}) \cdot \|\varphi - \varphi'\|_{\operatorname{L}^1(\nu_G)},
	\end{align*}
	which is bounded by \(\frac{\varepsilon}{3}\) by \eqref{eq-varphi-varphi'}. Combining this calculation with \eqref{eq-varphi-varphi'} and \eqref{eq-Approximation_varphi'} and using the triangle inequality yields
	\begin{equation}\label{eq-Approximation_varphi_epsilon}
		\bigg\| \frac{1}{\nu_G(\widetilde{U}_n)} \cdot (1_{\widetilde{U}_n} * \varphi) - \varphi \bigg\|_{\operatorname{L}^1(\nu_G)} < \varepsilon
	\end{equation}
	for all \(n \geq N\), which proves assertion (ii).
	
	At last we show that (ii) implies (i). Thus let \(\widetilde{U}\) be a measurable neighborhood of \(DeD\) and fix \(\varepsilon > 0\). Then
	\[
		r_G^{-1}(\widetilde{U}) = \{\hat{g} \in \widehat{G}: D\hat{g}D \in \widetilde{U}\}
	\]
	is a measurable neighborhood of \(e\), so it contains a clopen normal subgroup \(W'\). Then \(W := DW'\) is also a clopen subgroup and since \(r_G^{-1}(\widetilde{U})\) is \(D\) left and right invariant \(W \subseteq r_G^{-1}(\widetilde{U})\). In particular \(D\backslash W/D \subseteq \widetilde{U}\). By (ii) we have
	\[
		\lim_{n \to \infty} \frac{1}{\nu_G(\widetilde{U}_n)} \cdot (1_{\widetilde{U}_n} * 1_{D\backslash W/D}) = 1_{D \backslash W/D}
	\]
	in \(\operatorname{L}^1(\nu_G)\). Thus there is some \(N \in \mathbb{N}\) such that
	\[
		\bigg\| \frac{1}{\nu_G(\widetilde{U}_n)} \cdot (1_{\widetilde{U}_n} * 1_{D\backslash W/D}) - 1_{D \backslash W/D} \bigg\|_{\operatorname{L}^1(\nu_G)} < \varepsilon \cdot \nu_G(D\backslash W/D)
	\]
	for all \(n \geq N\). Now fix \(n \geq N\). Equation \eqref{eq-Approximation_varphi_epsilon} can only hold if
	\[
		\frac{1}{\nu_G(\widetilde{U}_n)} \cdot (1_{\widetilde{U}_n} * 1_{D\backslash W/D}) > 1 - \varepsilon
	\]
	on a subset \(\widetilde{K} \subseteq  D\backslash W/D\) with positive measure, as else
	\[
		\bigg\| \frac{1}{\nu_G(\widetilde{U}_n)} \cdot (1_{\widetilde{U}_n} * 1_{D\backslash W/D}) - 1_{D \backslash W/D} \bigg\|_{\operatorname{L}^1(\nu_G)} \geq \big\|\varepsilon \cdot 1_{D\backslash W/D}\big\|_{\operatorname{L}^1(\nu_G)} = \varepsilon \cdot \nu_G(D\backslash W/D)
	\]
	would hold. Then we calculate for almost all \(\hat{g} \in W\) with \(D\hat{g}D \in \widetilde{K}\) that
	\begin{align*}
		1 - \varepsilon &< \frac{1}{\nu_G(\widetilde{U}_n)} \cdot (1_{\widetilde{U}_n} * 1_{D\backslash W/D})(D\hat{g}D)\\
		&= \frac{1}{\nu_G(\widetilde{U}_n)} \int_{\widehat{G}} 1_{\widetilde{U}_n}(D\hat{f}D) \cdot 1_{D\backslash W/D}(D\hat{f}^{-1}\hat{g}D)\ \mathrm{d}\mu_G(\hat{f}).\
		\intertext{We use that \(D\hat{f}^{-1}\hat{g}D \in D\backslash W/D\) if and only \(\hat{f} \in W\) since \(\hat{g} \in W\), so}
		&= \frac{1}{\nu_G(\widetilde{U}_n)} \int_{\widehat{G}} 1_{r_G^{-1}(U)}(\hat{f}) \cdot 1_W(\hat{f})\ \mathrm{d}\mu_G(\hat{f})\\
		&= \frac{\mu_G(r_G^{-1}(\widetilde{U}_n) \cap W)}{\nu_G(\widetilde{U}_n)}.\\
		\intertext{Now we use that \(r_G\) is measure preserving and \(W = r_G^{-1}(D\backslash W/D)\) to obtain}
		&= \frac{\nu_G(\widetilde{U}_n \cap D\backslash W/D)}{\nu_G(\widetilde{U}_n)}\\
		&= 1 - \frac{\nu_G(\widetilde{U}_n \setminus (D\backslash W/D))}{\nu_G(\widetilde{U}_n)}.
	\end{align*}
	This implies since \(D\backslash W/D \subseteq \widetilde{U}\) that
	\[
		\frac{\nu_G(\widetilde{U}_n \setminus \widetilde{U})}{\nu_G(\widetilde{U}_n)} \leq \frac{\nu_G(\widetilde{U}_n \setminus (D\backslash W/D))}{\nu_G(\widetilde{U}_n)} < \varepsilon
	\]
	for all \(n \geq N\), which proves (i).
\end{proof*}

\begin{proof*}[\cref{lem-Continuous_Map_from_Measurable_Map}]
	Recall that \((V_n)_{n \in \mathbb{N}}\) is a decreasing sequence of clopen subgroups in \(\widehat{H}\) such that \(E = \bigcap_{n \in \mathbb{N}} V_n\) and
	\[
		\widetilde{V}_n := E\backslash V_n/E \subseteq E\backslash\widehat{H}/E.
	\]
	We are given a measure preserving map \(\theta': E\backslash\widehat{H}/E \to D\backslash\widehat{G}/D\) such that for every continuous \(\varphi \in \mathcal{H}(\widehat{G}, D)\)
	\[
		\lim_{n \to \infty} \frac{1}{\nu_H(\widetilde{V}_n)} \cdot (1_{\widetilde{V}_n} * (\varphi \circ \theta')) = \varphi \circ \theta'
	\]
	uniformly outside of some measure-zero set. We want to construct a continuous function \(\theta\) that is almost everywhere equivalent to \(\theta'\). To prove this we first show that the measure preserving map \(\zeta': \widehat{H} \to D\backslash\widehat{G}/D\) given by \(\zeta' := \theta' \circ r_H\), i.e.
	\[
		\zeta'(\hat{h}) = \theta'(E\hat{h}E)
	\]
	for all \(\hat{h} \in \widehat{H}\), is almost everywhere equivalent to some continuous \(\zeta: \widehat{H} \to D\backslash\widehat{G}/D\).
	
	As \(D\backslash\widehat{G}/D\) is totally disconnected and second countable we can fix a sequence \((\mathcal{D}_m)_{m \in \mathbb{N}}\) of finite clopen partitions of \(D\backslash\widehat{G}/D\) such that \(\mathcal{D}_{m+1}\) is a refinement of \(\mathcal{D}_m\) for all \(m \in \mathbb{N}\) and \(\bigcup_{m \in \mathbb{N}} \mathcal{D}_m\) is a topological basis of \(D\backslash\widehat{G}/D\). By assumption we find a measure-zero set \(N \subseteq \widehat{H}\) and a natural number \(n_m\) for each \(m \in \mathbb{N}\) such that
	\begin{equation}\label{eq-Approximation_1_U}
		\bigg| \frac{1}{\nu_H(\widetilde{V}_{n_m})} \cdot \big(1_{\widetilde{V}_{n_m}} * (1_{\widetilde{U}} \circ \theta')\big)(E\hat{h}E) - (1_{\widetilde{U}} \circ \theta')(E\hat{h}E) \bigg| < \frac{1}{2}
	\end{equation}
	for all \(\widetilde{U} \in \mathcal{D}_m\) and \(\hat{h} \in \widehat{H} \setminus N\). The fact that we can choose the same \(n_m\) for all \(\widetilde{U} \in \mathcal{D}_m\) comes from \(\mathcal{D}_m\) being finite. Similarly we can choose the same measure-zero set \(N\) for all \(\widetilde{U}\) since \(\bigcup_{m \in \mathbb{N}} \mathcal{D}_m\) is countable. Without loss of generality we may assume that \((n_m)_{m \in \mathbb{N}}\) is strictly increasing.
	
	We now prove that \(\zeta'\) maps almost all elements in a coset \(\hat{h}V_{n_m}\) into the same set in the partition \(\mathcal{D}_m\). Therefore fix \(m \in \mathbb{N}\) and \(\hat{h} \in \widehat{H}\). For any \(\hat{h}' \in \hat{h}V_{n_m} \setminus N\) and \(\widetilde{U} \in \mathcal{D}_m\) we calculate:
	\begin{align*}
		\frac{1}{\nu_H(\widetilde{V}_{n_m})} \cdot (1_{\widetilde{V}_{n_m}} * (1_{\widetilde{U}} \circ \theta'))(E\hat{h}'E) &= \frac{1}{\nu_H(\widetilde{V}_{n_m})} \int_{\widehat{H}} 1_{\widetilde{V}_{n_m}}(E\hat{f}E) \cdot 1_{\widetilde{U}}(\theta'(E\hat{f}^{-1}\hat{h}'E))\ \mathrm{d}\mu_H(\hat{f})\\
		&= \frac{1}{\mu_H(V_{n_m})} \int_{V_{n_m}} 1_{\widetilde{U}}(\zeta'(\hat{f}^{-1}\hat{h}'))\ \mathrm{d}\mu_H(\hat{f}),\\
		\intertext{where we used that \(r_H^{-1}(\widetilde{V}_{n_m}) = V_{n_m}\). Now we use the measure preserving transformation \(\hat{f} \mapsto \hat{h}'\hat{f}\) to obtain}
		&= \frac{1}{\mu_H(V_{n_m})} \int_{\hat{h}'^{-1}V_{n_m}} 1_{\zeta'^{-1}(\widetilde{U})}(\hat{f}^{-1})\ \mathrm{d}\mu_H(\hat{f})\\
		&= \frac{\mu_H(V_{n_m}\hat{h}' \cap \zeta'^{-1}(\widetilde{U}))}{\mu_H(V_{n_m})}.
	\end{align*}
	Plugging this into \eqref{eq-Approximation_1_U} we see that
	\[
		\frac{\mu_H(V_{n_m}\hat{h}' \cap \zeta'^{-1}(\widetilde{U}))}{\mu_H(V_{n_m})} > \frac{1}{2}
	\]
	if and only if \((1_{\tilde{U}} \circ \theta')(E\hat{h}'E) = 1\), i.e. if and only if \(\zeta'(\hat{h}') = \theta'(E\hat{h}'E) \in \tilde{U}\). For \(\hat{h}', \hat{h}'' \in V_{n_m}\hat{h}\setminus N \) and all \(\widetilde{U} \in \mathcal{D}_m\) we have
	\[
		\frac{\mu_H(V_{n_m}\hat{h}' \cap \zeta'^{-1}(\widetilde{U}))}{\mu_H(V_{n_m})} = \frac{\mu_H(V_{n_m}\hat{h}'' \cap \zeta'^{-1}(\widetilde{U}))}{\mu_H(V_{n_m})},
	\]
	so \(\zeta'(\hat{h}')\) and \(\zeta'(\hat{h}'')\) must be in the same set \(\tilde{U} \in \mathcal{D}_m\). Thus we can denote by \(\widetilde{U}_{\hat{h}}^m\) the set in \(\mathcal{D}_m\) such that \(\zeta'(\hat{h}') \in \widetilde{U}_{\hat{h}}^m\) for all \(\hat{h}' \in V_{n_m}\hat{h} \setminus N\). Note that if \(\hat{h} \in N\) we might not have \(\zeta'(h) \in \widetilde{U}_h^m\). Also note that we have \(\widetilde{U}_{\hat{h}'}^m = \widetilde{U}_{\hat{h}}^m\) for all \(\hat{h}' \in V_{n_m}\hat{h}\).
	
	Next we show that the sequence \((\widetilde{U}_{\hat{h}}^m)_{m \in \mathbb{N}}\) is strictly decreasing for each \(\hat{h} \in \widehat{H}\). Therefore fix \(\hat{h} \in \widehat{H}\) and \(m \in \mathbb{N}\). Take \(\hat{h}' \in V_{n_{m+1}}\hat{h}\setminus N\). Since \(V_{n_{m+1}} < V_{n_m}\) we also have \(\hat{h}' \in V_{n_m}\hat{h} \setminus N\). Thus \(\zeta'(\hat{h}') \in \widetilde{U}_{\hat{h}}^m\) and \(\zeta'(\hat{h}') \in \widetilde{U}_{\hat{h}}^{m+1}\) implies \(\widetilde{U}_{\hat{h}}^{m+1} \subseteq \widetilde{U}_{\tilde{h}}^m\) since \(\mathcal{D}_{m+1}\) is a refinement of \(\mathcal{D}_m\). This shows that \((\widetilde{U}_{\hat{h}}^m)_{m \in \mathbb{N}}\) is in fact strictly decreasing.
	
	As \(\bigcup_{m \in \mathbb{N}} \mathcal{D}_m\) is a basis of \(D\backslash\widehat{G}/D\) it follows that the intersection of \((\widetilde{U}_{\hat{h}}^m)_{m \in \mathbb{N}}\) consists of a single point. We define the map \(\zeta\) such that it maps \(\hat{h}\) to this point, i.e. \(\zeta: \widehat{H} \to D\backslash\widehat{G}/D\) is defined such that
	\[
		\{\zeta(\hat{h})\} = \bigcap_{m \in \mathbb{N}} \widetilde{U}_{\hat{h}}^m
	\]
	holds for all \(\hat{h} \in \widehat{H}\). For \(\hat{h} \in \widehat{H} \setminus N\) we have \(\zeta'(\hat{h}) \in \widetilde{U}_{\hat{h}}^m\) for all \(m \in \mathbb{N}\), so \(\zeta\) and \(\zeta'\) agree on \(\widehat{H} \setminus N\), in particular they agree almost everywhere.
	
	Now we prove that \(\zeta\) is continuous. Since \(\bigcup_{m \in \mathbb{N}} \mathcal{D}_m\) is a basis it suffices to show for \(m \in \mathbb{N}\) and \(\widetilde{U} \in \mathcal{D}_m\) that \(\zeta^{-1}(\widetilde{U})\) is open. This is clear if \(\zeta^{-1}(\widetilde{U})\) is empty. If it is not empty we have for all \(\hat{h} \in \zeta^{-1}(\widetilde{U})\) that \(V_{n_m}\hat{h}\) is an open neighborhood of \(\hat{h}\) and for \(\hat{h}' \in V_{n_m}\hat{h}\) that
	\[
		\zeta(\hat{h}') \in \bigcap_{n \in \mathbb{N}} \widetilde{U}_{\hat{h}'}^n \subseteq \widetilde{U}_{\hat{h}'}^m = \widetilde{U}_{\hat{h}}^m = \widetilde{U}.
	\]
	Thus \(V_{n_m}\hat{h} \subseteq \zeta^{-1}(\widetilde{U})\) for all \(\hat{h} \in \zeta^{-1}(\{\tilde{U}\})\), so \(\zeta^{-1}(\widetilde{U})\) is in fact open.
	
	At last we construct \(\theta\) from \(\zeta\). To do so, we have to show that \(\zeta\) is \(E\) left and right invariant. Recall that \(\zeta'\) was constructed the following way:
	\[
		\begin{tikzcd}
			\widehat{H} \arrow[d, ->>, "r_H"] \arrow[rd, "\zeta'"]\\
			E\backslash\widehat{H}/E \arrow[r, "\theta'"] &D\backslash\widehat{G}/D
		\end{tikzcd}
	\]
	From this construction it is clear that \(\zeta'\) is \(E\) left and right invariant. Fix \(\hat{h} \in \widehat{H}\) and \(d, d' \in E\). For almost all \(\hat{h}' \in \widehat{H}\) we have \(\zeta(\hat{h}') = \zeta'(\hat{h}')\) and \(\zeta(d\hat{h}'d') = \zeta'(d\hat{h}'d')\). In particular every open set contains an element that fulfills both these equalities. Thus we can construct a sequence \((\hat{h}_n)_{n \in \mathbb{N}}\) in \(\widehat{H}\) that converges to \(\hat{h}\) such that \(\zeta(\hat{h}_n) = \zeta'(\hat{h}_n)\) and \(\zeta(d\hat{h}_nd') = \zeta'(d\hat{h}_nd')\) for all \(n \in \mathbb{N}\). Therefore, since \(\zeta\) is continuous, we have
	\[
		\zeta(d\hat{h}d') = \lim_{n \to \infty} \zeta(d\hat{h}_nd') = \lim_{n \to \infty} \zeta'(d\hat{h}_nd') = \lim_{n \to \infty} \zeta'(\hat{h}_n) = \lim_{n \to \infty} \zeta(\hat{h}_n) = \zeta(\hat{h}),
	\]
	so \(\zeta\) is also \(E\) left and right invariant. Thus there exists a unique continuous map \(\theta: E\backslash\widehat{H}/E \to D\backslash\widehat{G}/D\) such that \(\zeta = \theta \circ r_H\). As \(\zeta\) and \(\zeta'\) are equivalent almost everywhere and \(r_H\) is measure preserving it follows that also \(\theta = \theta'\) almost everywhere.
\end{proof*}

\begin{proof*}[\cref{the-Isometric_Isomorphisms_Hecke_Algebras}]
	We are given an isometric *-algebra isomorphism \(\mathcal{T}: \mathcal{H}(\widehat{G}, D) \to \mathcal{H}(\widehat{H}, E)\). By \cref{lem-T1_X=1_Y} there exists a unique continuous group homomorphism \(\chi: \widehat{H} \to \mathbb{S}^1\) with \(E \leq \ker(\chi)\) such that \(\mathcal{T}1_{D\backslash\widehat{G}/D}(E\hat{h}E) = \chi(\hat{h})\) for \(\nu_H\)-almost all \(E\hat{h}E \in E\backslash\widehat{H}/E\). Further the same lemma yields that \(\mathcal{T}': \mathcal{H}(\widehat{G}, D) \to \mathcal{H}(\widehat{H}, E)\) defined by
	\[
		\mathcal{T}'\varphi(E\hat{h}E) := \overline{\chi(\hat{h})} \cdot \mathcal{T}\varphi(E\hat{h}E)
	\]
	for all \(\varphi \in \mathcal{H}(\widehat{G}, D)\) and \(\nu_H\)-almost all \(E\hat{h}E \in E\backslash\widehat{H}/E\) is an isometric *-algebra isomorphism, with \(\mathcal{T}'1_{D\backslash\widehat{G}/D} = 1_{E\backslash\widehat{H}/E}\). Thus by \cref{cor-Banach-Lamperti} there is a measure preserving map \(\theta': E\backslash\widehat{H}/E \to D\backslash\widehat{G}/D\) such that
	\[
		\mathcal{T}'\varphi = \varphi \circ \theta'
	\]
	for all \(\varphi \in \mathcal{H}(\widehat{G}, D)\). The map \(\theta'\) is unique up to measure zero sets. Now we need to prove that \(\theta'\) is almost everywhere equivalent to some homeomorphism \(\theta\). Once this is proven, the uniqueness of \(\theta\) follows from the almost uniqueness of \(\theta'\).
	
	Let \((V_n)_{n \in \mathbb{N}}\) be a decreasing sequence of clopen subgroups in \(\widehat{H}\) such that \(E = \bigcap_{n \in \mathbb{N}} V_n\). Such a sequence can be constructed by first taking a decreasing sequence of clopen normal subgroups \((V_n')_{n \in \mathbb{N}}\) with \(\{e\}\) as the intersection and then setting \(V_n := V_n'E\). Denote
	\[
		\widetilde{V}_n := E\backslash V_n/E \subseteq E\backslash\widehat{H}/E.
	\]
	Then \((\widetilde{V}_n)_{n \in \mathbb{N}}\) is a decreasing clopen neighborhood basis of \(EeE\). Since \(\mathcal{T}'\) is an isomorphism it follows that \(1_{\widetilde{V}_n}\) is in its range and, as \(\mathcal{T}'\) is the Koopman operator of \(\theta'\), this implies that there is a decreasing sequence \((\widetilde{U}_n)_{n \in \mathbb{N}_0}\) of Borel sets in \(D\backslash\widehat{G}/D\) such that \(\mathcal{T}'1_{\widetilde{U}_n} = 1_{\widetilde{V}_n}\), i.e.
	\[
		\widetilde{V}_n = \theta'^{-1}(\widetilde{U}_n)
	\]
	up to measure zero. Note that at this moment we do not know if the \(\tilde{U}_n\) are neighborhoods of the \(DeD\). From \cref{lem-Approximate_Unit} applied in \(\mathcal{H}(\widehat{H}, E)\) we obtain that
	\[
		\lim_{n \to \infty} \frac{1}{\nu_H(\widetilde{V}_n)} \cdot (1_{\widetilde{V}_n} * \psi) = \psi
	\]
	in \(\operatorname{L}^1(\nu_H)\) for all \(\psi \in \mathcal{H}(\widehat{H}, E)\). Thus since \(\mathcal{T}'\) is an isometric *-algebra isomorphism we have for all \(\varphi \in \mathcal{H}(\widehat{G}, D)\) that
	\begin{align}\label{eq-Approximation_varphi}
		\lim_{n \to \infty} \frac{1}{\nu_G(\widetilde{U}_n)} \cdot (1_{\widetilde{U}_n} * \varphi) &= \lim_{n \to \infty} \frac{1}{\nu_H(\theta'^{-1}(\widetilde{U}_n))} \cdot (\mathcal{T}'^{-1}1_{\widetilde{V}_n} * \mathcal{T}'^{-1}\mathcal{T}'\varphi)\\
		&= \lim_{n \to \infty} \frac{1}{\nu_H(\widetilde{V}_n)} \cdot \mathcal{T}'^{-1}(1_{\widetilde{V}_n} * \mathcal{T}'\varphi)\nonumber\\
		&= \mathcal{T}'^{-1}\mathcal{T}'\varphi\nonumber\\
		&= \varphi\nonumber
	\end{align}
	in \(\operatorname{L}^1(\nu_G)\). Hence by \cref{lem-Approximate_Unit} (this time applied in \(\mathcal{H}(\widehat{G}, D)\)) it follows that for all neighborhoods \(\widetilde{U}\) of \(DeD\)
	\begin{equation}\label{eq-U_n_Approximate_Neighborhood}
		\lim_{n \to \infty} \frac{\nu_G(\widetilde{U}_n \setminus \widetilde{U})}{\nu_G(\widetilde{U}_n)} = 0
	\end{equation}
	holds. This can be interpreted as \((\widetilde{U}_n)_{n \in \mathbb{N}}\) being approximately a neighborhood basis of \(DeD\). It follows from assertion (iii) of \cref{lem-Approximate_Unit} that for continuous \(\varphi \in \mathcal{H}(\widehat{G}, D)\) equation \eqref{eq-Approximation_varphi} also holds with uniform convergence. Since \(\mathcal{T}'\) is the Koopman operator of \(\theta'\), it preserves the supremum norm. Now, using that \(\mathcal{T}'\) is a *-algebra homomorphism, we have for all continuous \(\varphi \in \mathcal{H}(\widehat{G}, D)\) that also
	\[
		\lim_{n \to \infty} \frac{1}{\nu_H(\widetilde{V}_n)} \cdot (1_{\widetilde{V}_n} * (\varphi \circ \theta')) = \lim_{n \to \infty} \frac{1}{\nu_G(\widetilde{U}_n)} \cdot \mathcal{T}'(1_{\widetilde{U}_n} * \varphi) = \mathcal{T}'\varphi = \varphi \circ \theta'
	\]
	uniformly outside of some measure-zero set. This measure-zero set comes from the fact that the multiplicativity of \(\mathcal{T}'\) only holds almost everywhere. This is exactly the prerequisite for \cref{lem-Continuous_Map_from_Measurable_Map}, so we obtain a continuous \(\theta: E\backslash\widehat{H}/E \to D\backslash\widehat{G}/D\) that is almost everywhere equivalent to \(\theta'\). Since \(\mathcal{T}'\) is an isomorphism it is clear that \(\theta\) must be bijective and since both \(D\backslash\widehat{G}/D\) and \(E\backslash\widehat{H}/E\) are compact Hausdorff, it must be a homeomorphism.
	
	Now we prove that \(\theta(EeE) = DeD\). Since \((\widetilde{V}_n)_{n \in \mathbb{N}}\) is a neighborhood basis of \(EeE\) and \(\theta\) is a homeomorphism we know that \((\theta(\widetilde{V}_n))_{n \in \mathbb{N}}\) is a clopen neighborhood basis of \(\theta(EeE)\). Also \(\theta(\widetilde{V}_n)\) differs from \(\widetilde{U}_n\) only by a measure zero set, so from \eqref{eq-U_n_Approximate_Neighborhood} we obtain
	\[
		\lim_{n \to \infty} \frac{\nu_G(\theta(\widetilde{V}_n) \setminus \widetilde{U})}{\nu_G(\theta(\widetilde{V}_n))} = 0
	\]
	for all neighborhoods \(\tilde{U}\) of \(DeD\). If \(\theta(EeE) \neq DeD\) we would, since \(D\backslash\widehat{G}/D\) is Hausdorff, find a neighborhood \(\widetilde{U}\) of \(DeD\) that is disjoint with \(\theta(\widetilde{V}_n)\) for sufficiently large \(n \in \mathbb{N}_0\). This contradicts the convergence above, so \(\theta(EeE) = DeD\) must hold.
	
	Finally we have for the original operator \(\mathcal{T}\) that for almost all \(E\hat{h}E \in E\backslash\widehat{H}/E\)
	\[
		\mathcal{T}\varphi(E\hat{h}E) = \chi(\hat{h}) \cdot \mathcal{T}'\varphi(E\hat{h}E) = \chi(\hat{h}) \cdot \varphi(\theta(E\hat{h}E)).
	\]
\end{proof*}

\subsection{Consequences for Group Chains and Invariance of the Steinitz Order}\label{sub-Steinitz_Order}

Let \((X, G, \Phi)\) and \((Y, H, \Psi)\) be two minimal equicontinuous Cantor action. Fix group chains \((G_n)_{n \in \mathbb{N}_0}\) and \((H_n)_{n \in \mathbb{N}_0}\) associated to these actions and denote the arising Ellis groups by \(\widehat{G}\) and \(\widehat{H}\) respectively as well as their discriminant groups by \(D\) and \(E\) respectively. We label the various maps and measures from \cref{sec-Cantor_Actions,sec-Hecke_Algebras} (primarily diagrams \eqref{dia-Inverse_Limit_Double_Cosets} and \eqref{dia-Hecke_Algebras}) with sub- or super-scrips \(G\) or \(H\) to specify to which action and group chain they belong.

\begin{theorem}\label{the-Group_Chains_Same_Index}
	Let \((X, G, \Phi)\) and \((Y, H, \Psi)\) be two minimal equicontinuous Cantor actions and fix associated group chains \((G_n)_{n \in \mathbb{N}_0}\) and \((H_n)_{n \in \mathbb{N}_0}\).	Assume that there exists an isometric *-algebra isomorphism \(\mathcal{T}: \mathcal{H}(\Phi) \to \mathcal{H}(\Psi)\). Then there exists a group chain \((H_n')_{n \in \mathbb{N}_0}\) in \(H\) with
	\[
		|H/H_n'| = |G/G_n|
	\]
	for all \(n \in \mathbb{N}_0\) that is also associated to \((Y, H, \Psi)\).
	
	If we additionally assume that \(\mathcal{T}\) is a Koopman operator, then there are isometric *-algebra homomorphisms \(\mathcal{T}_n: \mathcal{H}(G, G_n) \to \mathcal{H}(H, H_n')\) such that the following diagram commutes for \(m \geq n\):
	\[
		\begin{tikzcd}
			\mathcal{H}(\Phi) \arrow[rrr, "\cong\ \mathcal{T}"] & & &\mathcal{H}(\Psi)\\
			&\mathcal{H}(G, G_m) \arrow[lu, hook', "\mathcal{Q}_m^G"] \arrow[r, hook, "\mathcal{T}_m"] &\mathcal{H}(H, H_m') \arrow[ru, hook, "\mathcal{Q}_m^{\prime H}"]\\
			&\mathcal{H}(G, G_n) \arrow[luu, hook', bend left, "\mathcal{Q}_n^G"] \arrow[u, hook, "\mathcal{Q}_{n, m}^G"] \arrow[r, hook, "\mathcal{T}_n"] &\mathcal{H}(H, H_n') \arrow[ruu, hook, bend right, "\mathcal{Q}_n^{\prime H}"] \arrow[u, hook, "\mathcal{Q}_{n, m}^{\prime H}"]
		\end{tikzcd}
	\]
	Here \(\mathcal{Q}_{n, m}^{\prime H}: \mathcal{H}(H, H_n') \to \mathcal{H}(H, H_m')\) and \(\mathcal{Q}_n^{\prime H}: \mathcal{H}(H, H_n') \to \mathcal{H}(\Psi)\) denote the Koopman operators of the measure preserving continuous maps \(q_{n, m}^{\prime H}: H_m'\backslash H/H_m' \to H_n'\backslash H/H_n'\) and \(q_n^{\prime H}: E\backslash\widehat{H}/E \to H_n'\backslash\widehat{H}/E\) (see diagrams \eqref{dia-Inverse_Limit_Double_Cosets} and \eqref{dia-Hecke_Algebras}).
\end{theorem}

An immediate consequence of this theorem is that if the Hecke algebras of two minimal equicontinuous Cantor actions are isometrically isomorphic then they have the same Steinitz order (\cref{def-Steinitz-Order}).

\begin{corollary}\label{cor-Steinitz_Order}
	Assume that \(\mathcal{H}(\Phi)\) and \(\mathcal{H}(\Psi)\) are isometrically *-isomorphic. Then the Cantor actions \((X, G, \Phi)\) and \((Y, H, \Psi)\) have the same Steinitz order.
\end{corollary}
\begin{proof*}
	This follows from combining \cref{the-Group_Chains_Same_Index} directly above with the formula of the Steinitz order given in \cref{pro-Steinitz_Order}. Here we use that this formula is independent of the choice of associated group chain.
\end{proof*}

\begin{proof*}[\cref{the-Group_Chains_Same_Index}]
	By \cref{the-Isometric_Isomorphisms_Hecke_Algebras} we find an isometric *-algebra isomorphism \(\mathcal{T}': \mathcal{H}(\Phi) \to \mathcal{H}(\Psi)\) that is the Koopman operator of a measure preserving homeomorphism \(\theta: E\backslash\widehat{H}/E \to D\backslash\widehat{G}/D\). Fix \(n \in \mathbb{N}_0\). Then \(\mathcal{T}'\mathcal{Q}_n^G: \mathcal{H}(G, G_n) \to \mathcal{H}(\Psi)\) is an isometric *-algebra homomorphism and the Koopman operator of the continuous map \(q_n^G \circ \theta: E\backslash\widehat{H}/E \to G_n\backslash G/G_n\). The set
	\[
		\widehat{H}_n' := \{\hat{h} \in \widehat{H}: (q_n^G \circ \theta)(E\hat{h}E) = G_neG_n\} = (q_n^G \circ \theta \circ r_H)^{-1}(\{G_neG_n\})
	\]
	is a clopen subset of \(\widehat{H}\) that contains \(E\). We will prove that \(\widehat{H}_n'\) is a subgroup of \(\widehat{H}\) with index \(|G/G_n|\). Since for almost all \(E\hat{h}E \in E\backslash\widehat{H}/E\) we have
	\[
		\mathcal{T}'\mathcal{Q}_n^G 1_{\{G_neG_n\}}(E\hat{h}E) = 1_{(q_n^G \circ \theta)^{-1}(\{G_neG_n\})}(E\hat{h}E) = 1_{E\backslash\widehat{H}_n'/E}(E\hat{h}E)
	\]
	and, since \(|G/G_n| \cdot 1_{\{G_neG_n\}}\) is the unit element in \(\mathcal{H}(G, G_n)\), we obtain
	\begin{alignat*}{3}
		(1_{E\backslash\widehat{H}_n'/E} * 1_{E\backslash\widehat{H}_n'/E})(E\hat{h}E) &= \mathcal{T}'\mathcal{Q}_n^G(1_{\{G_neG_n\}} * 1_{\{G_neG_n\}})(E\hat{h}E)\\
		&= \frac{1}{|G/G_n|} \cdot \mathcal{T}'\mathcal{Q}_n^G1_{\{G_neG_n\}}(E\hat{h}E) &= \frac{1}{|G/G_n|} \cdot 1_{E\backslash\widehat{H}_n'/E}.
	\end{alignat*}
	Thus we have for almost all \(\hat{h} \in \widehat{H}_n'\) that
	\begin{align*}
		1 &= 1_{E\backslash\widehat{H}_n'/E}(E\hat{h}E)\\
		&= |G/G_n| \cdot (1_{E\backslash\widehat{H}_n'/E} * 1_{E\backslash\widehat{H}_n'/E})(E\hat{h}E)\\
		&= |G/G_n| \int_{\widehat{H}} 1_{E\backslash\widehat{H}_n'/E}(E\hat{f}E) \cdot 1_{E\backslash\widehat{H}_n'/E}(E\hat{f}^{-1}\hat{h}E)\ \mathrm{d}\mu_H(\hat{f})\\
		&= |G/G_n| \int_{\widehat{H}_n'} 1_{\widehat{H}_n'}(\hat{f}^{-1}\hat{h})\ \mathrm{d}\mu_H(\hat{f}).
	\end{align*}
	From this it follows with
	\begin{equation}\label{eq-Index_H_n'}
		\mu_H(\widehat{H}_n') = \mu_H((q_n^G \circ \theta \circ r_H)^{-1}(\{G_neG_n\})) \overset{\ref{pro-Maps_Measure_Preserving}}{=} \nu_n^G(\{G_neG_n\}) \overset{\eqref{eq-nu_n}}{=} \frac{1}{|G/G_n|}
	\end{equation}
	that we must have \(\hat{f}^{-1}\hat{h} \in \widehat{H}_n'\) for almost all \(\hat{h}, \hat{f} \in \widehat{H}_n'\). Therefore, since \(\widehat{H}_n'\) is clopen, we find for all \(\hat{h}, \hat{f} \in \widehat{H}_n'\) sequences \((\hat{h}_m)_{m \in \mathbb{N}}\) and \((\hat{f}_m)_{m \in \mathbb{N}}\) that converge to \(\hat{h}\) and \(\hat{f}\) respectively with \(\hat{f}_m^{-1}\hat{h}_m \in \widehat{H}_n'\). Hence
	\[
		\hat{f}^{-1}\hat{h} = \lim_{m \to \infty} \hat{f}_m^{-1}\hat{h}_m \in \widehat{H}_n',
	\]
	so \(\widehat{H}_n'\) is indeed a clopen subgroup of \(\widehat{H}\). For its index we have since \(\mu_H\) is the Haar measure that
	\[
		|\widehat{H}/\widehat{H}_n'| = \frac{1}{\mu_H(\widehat{H}_n')} \overset{\eqref{eq-Index_H_n'}}{=} |G/G_n|.
	\]
	For \(m \geq n\) and \(\hat{h} \in \widehat{H}_m'\) we have
	\[
		(q_n^G \circ \theta)(E\hat{h}E) \overset{\eqref{dia-Inverse_Limit_Double_Cosets}}{=} q_{n, m}^G \circ (q_m^G \circ \theta)(E\hat{h}E) = q_{n, m}^G(G_meG_m) = G_neG_n,
	\]
	so \(\hat{h} \in \widehat{H}_n'\) and \((\widehat{H}_n')_{n \in \mathbb{N}_0}\) is decreasing.
	
	Now we construct the subgroups \(H_n'\) in \(H\) from these subgroups \(\widehat{H}_n'\) in \(\widehat{H}\). Recall that \(\zeta: H \to \widehat{H}\) defined by \(\zeta(h) := (hC_k)_{k \in \mathbb{N}_0}\) is a group homomorphism with dense image. Thus
	\[
		H_n' := \zeta^{-1}(\widehat{H}_n') = \zeta^{-1}(\widehat{H}_n' \cap \zeta(H))
	\]
	defines a subgroup of \(H\) with index
	\[
		|H/H_n'| = |\widehat{H}/\widehat{H}_n'| = |G/G_n|
	\]
	for all \(n \in \mathbb{N}_0\). So \((H_n')_{n \in \mathbb{N}_0}\) is a group chain in \(H\) with \(|H/H_n'| = |G/G_n|\) for all \(n \in \mathbb{N}_0\).
	
	We now prove that \((H_n')_{n \in \mathbb{N}_0}\) is associated to the minimal equicontinuous Cantor action \((Y, H, \Psi)\), i.e. we construct an adapted neighborhood basis \((U_n')_{n \in \mathbb{N}_0}\) in \(Y\), such that \(H_n'\) consists of exactly those group elements that map \(U_n'\) into itself. By \cref{the-Structure_Theorem_Equicontinuous_Actions} we may assume \(Y = \widehat{H}/E\), and by \cref{the-Inverse_Limit_Ellis_Group} \(\widehat{H}\) is given by \(\varprojlim_{n \in \mathbb{N}_0} H/C_n\), where \(C_n\) is the core of \(H_n\). The sets \((U_n')_{n \in \mathbb{N}_0} := (\widehat{H}_n'/E)_{n \in \mathbb{N}_0} = (\overline{\zeta(H_n')}/E)_{n \in \mathbb{N}_0}\) are a decreasing sequence of adapted neighborhoods of \(eE\). It is a neighborhood basis since
	\[
		\bigcap_{n \in \mathbb{N}} E\backslash\widehat{H}_n'/E = \theta^{-1}\bigg( \bigcap_{n \in \mathbb{N}} (q_n^G)^{-1}(\{G_neG_n\}) \bigg) = \theta^{-1}(\{DeD\}) = \{EeE\},
	\]
	implies \(\bigcap_{n \in \mathbb{N}} \widehat{H}_n' = E\) and thus \(\bigcap_{n \in \mathbb{N}_0} \widehat{H}_n'/E = E/E\). The group chain arising from the neighborhood basis \((\widehat{H}_n'/E)_{n \in \mathbb{N}_0}\) is \((H_n')_{n \in \mathbb{N}_0}\). Hence \((H_n')_{n \in \mathbb{N}_0}\) is associated to \((Y, H, \Psi)\).
	
	Now we assume that \(\mathcal{T}\) itself is the Koopman operator of \(\theta\) and prove the existence of the \(\mathcal{T}_n: \mathcal{H}(G, G_n) \to \mathcal{H}(H, H_n')\). The existence of such an isometric *-algebra homomorphism follows as soon as we show that
	\begin{equation}\label{eq-Inclusion_Range}
		\operatorname{ran}(\mathcal{T}\mathcal{Q}_n^G) \subseteq \operatorname{ran}(\mathcal{Q}_n^{\prime H})
	\end{equation}
	holds: since \(\mathcal{Q}_n^{\prime H}\) is injective, we can, if the above holds, define \(\mathcal{T}_n\) such that it maps \(\varphi \in \mathcal{H}(G, G_n)\) to the unique \(\psi \in \mathcal{H}(H, H_n')\) with \(\mathcal{Q}_n^{\prime H}\psi = \mathcal{T}\mathcal{Q}_n^G\varphi\). This yields an isometric *-algebra homomorphism, since the operators \(\mathcal{T}\), \(\mathcal{Q}_n^G\) and \(\mathcal{Q}_n^{\prime H}\) are all isometric *-algebra homomorphisms. By construction it fulfills \(\mathcal{T}\mathcal{Q}_n^G = \mathcal{Q}_n^{\prime H}\mathcal{T}_n\).
	
	We now prove \eqref{eq-Inclusion_Range}. Fix \(\varphi \in \mathcal{H}(G, G_n)\). We first show that the map \(\hat{h} \mapsto \mathcal{T}\mathcal{Q}_n^G\varphi(E\hat{h}E)\) is almost everywhere \(\widehat{H}_n'\) left and right invariant. The function \(|G/G_n| \cdot 1_{\{G_neG_n\}}\) is the unit element in \(\mathcal{H}(G, G_n)\), hence for almost all \(\hat{h} \in \widehat{H}\) and \(\hat{d} \in \widehat{H}_n'\)
	\begin{align*}
		\mathcal{T}\mathcal{Q}_n^G\varphi(E\hat{d}\hat{h}E) &= |G/G_n| \cdot \mathcal{T}\mathcal{Q}_n^G(1_{\{G_neG_n\}} * \varphi)(E\hat{d}\hat{h}E)\\
		&= |G/G_n| \cdot (\mathcal{T}\mathcal{Q}_n^G1_{\{G_neG_n\}} * \mathcal{T}\mathcal{Q}_m^G\varphi)(E\hat{d}\hat{h}E)\\
		&= |G/G_n| \int_{\widehat{H}} 1_{E\backslash\widehat{H}_n'/E}(E\hat{f}E) \cdot \mathcal{T}\mathcal{Q}_m^G\varphi(E\hat{f}^{-1}\hat{d}\hat{h}E)\ \mathrm{d}\mu_H(\hat{f}).\\
		\intertext{Next we use the measure preserving transformation \(\hat{f} \mapsto \hat{d}\hat{f}\) and that, since \(\hat{d} \in \widehat{H}_n'\), we have \(\hat{d}\hat{f} \in \widehat{H}_n'\) if and only if \(\hat{f} \in \widehat{H}_n'\). Thus we obtain}
		&= |G/G_n| \int_{\widehat{H}} 1_{E\backslash\widehat{H}_n'/E}(E\hat{d}\hat{f}E) \cdot \mathcal{T}\mathcal{Q}_m^G\varphi(E\hat{f}^{-1}\hat{h}E)\ \mathrm{d}\mu_H(\hat{f})\\
		&= |G/G_n| \int_{\widehat{H}} 1_{E\backslash\widehat{H}_n'/E}(E\hat{f}E) \cdot \mathcal{T}\mathcal{Q}_m^G\varphi(E\hat{f}^{-1}\hat{h}E)\ \mathrm{d}\mu_H(\hat{f})\\
		&= |G/G_n| \cdot \mathcal{T}\mathcal{Q}_n^G(1_{\{G_neG_n\}} * \varphi)(E\hat{h}E)\\
		&= \mathcal{T}\mathcal{Q}_n^G\varphi(E\hat{h}E).
	\end{align*}
	Hence for almost all \(\hat{h} \in \widehat{H}\) and \(\hat{d},\hat{d}' \in \widehat{H}_n'\)
	\begin{alignat*}{3}
		\mathcal{T}\mathcal{Q}_n^G\varphi(E\hat{d}\hat{h}\hat{d}'E) &= \mathcal{T}\mathcal{Q}_n^G\varphi(E\hat{h}\hat{d}'E)\\
		&= \overline{\mathcal{T}\mathcal{Q}_n^G\varphi^*(E\hat{d}'^{-1}\hat{h}^{-1}E)}\\
		&= \overline{\mathcal{T}\mathcal{Q}_n^G\varphi^*(E\hat{h}^{-1}E)} &= \mathcal{T}\mathcal{Q}_n^G\varphi(E\hat{h}E).
	\end{alignat*}
	Thus the map \(\hat{h} \mapsto \mathcal{T}\mathcal{Q}_n^G\varphi(E\hat{h}E)\) is almost everywhere \(\widehat{H}_n'\) left and right invariant, so it is almost surely constant on all double cosets \(\widehat{H}_n'\hat{h}\widehat{H}_n'\) and can thus be identified with a map \(\widehat{\psi}: \widehat{H}_n'\backslash\widehat{H}/\widehat{H}_n' \to \mathbb{C}\).
	
	The space \(\widehat{H}_n'\backslash\widehat{H}/\widehat{H}_n'\) is in bijection with \(H_n'\backslash H/H_n'\) through \(H_n'hH_n' \mapsto \widehat{H}_n'\zeta(h)\widehat{H}_n'\), where \(\zeta: H \to \widehat{H}\) is the group homomorphism with dense image from earlier in this proof. Therefore we have a function \(\psi: H_n'\backslash H/H_n' \to \mathbb{C}\) such that
	\begin{equation}\label{eq-psi_psi-hat}
		\psi(H_n'hH_n') = \widehat{\psi}(\widehat{H}_n'\zeta(h)\widehat{H}_n)
	\end{equation}
	for all \(h \in H\). For \(h \in H\) and \(\hat{h} \in \widehat{H}_n'\zeta(h)\widehat{H}_n'\) we have \(q_n^{\prime H}(E\hat{h}E) = H_n'hH_n'\), so
	\begin{equation}\label{eq-Qpsi=psi-hat}
		\mathcal{Q}_n^{\prime H}\psi(E\hat{h}E) = \psi(q_n^{\prime H}(\widehat{H}_n'\hat{h}\widehat{H}_n'))) = \psi(H_n'hH_n') = \widehat{\psi}(\widehat{H}_n'\zeta(h)\widehat{H}_n') = \widehat{\psi}(\widehat{H}_n'\hat{h}\widehat{H}_n').
	\end{equation}
	Thus we can conclude for almost all \(\hat{h} \in \widehat{H}\) that
	\[
		\mathcal{Q}_n^{\prime H}\psi(E\hat{h}E) \overset{\eqref{eq-Qpsi=psi-hat}}{=} \widehat{\psi}(\widehat{H}_n'\hat{h}\widehat{H}_n') \overset{\eqref{eq-psi_psi-hat}}{=} \mathcal{T}\mathcal{Q}_n^G\varphi(E\hat{h}E)
	\]
	holds. This proves \(\mathcal{T}\mathcal{Q}_n^G\varphi = \mathcal{Q}_n^{\prime H}\psi \in \operatorname{ran}(\mathcal{Q}_n^{\prime H})\).
\end{proof*}

Another application of \cref{the-Group_Chains_Same_Index} is the following characterization whether \(\mathcal{H}(\Phi)\) and \(\mathcal{H}(\Psi)\) are isometrically isomorphic. They are isometrically isomorphic if and only if we find isometric zig-zag *-algebra homomorphisms between subsequences \((\mathcal{H}(G, G_{n_k}))_{k \in \mathbb{N}_0}\) and \((\mathcal{H}(H, H_{m_k}))_{k \in \mathbb{N}_0}\) of the finite dimensional Hecke algebras. This is similar to the proof of \cref{the-Equivalence_Group_Chains} (see \citep[Theorem 3.3]{dyer_discriminant_2016}), where equivalence of group chains of conjugated actions is proven by constructing zig-zag maps between cosets spaces.

\begin{theorem}\label{the-Hecke_Algebra_Zig-zag_Map}
	Let \((X, G, \Phi)\) and \((Y, H, \Psi)\) be two minimal equicontinuous Cantor actions and fix associated group chains \((G_n)_{n \in \mathbb{N}_0}\) and \((H_n)_{n \in \mathbb{N}_0}\). Then the following assertions are equivalent:
	\begin{enumerate}[(i)]
		\item there exists an isometric *-algebra isomorphism \(\mathcal{T}: \mathcal{H}(\Phi) \to \mathcal{H}(\Psi)\),
		\item there are strictly increasing sequences \((n_k)_{k \in \mathbb{N}_0}\) and \((m_k)_{k \in \mathbb{N}_0}\) of natural numbers, together with isometric *-algebra homomorphisms
		\(\mathcal{T}_k: \mathcal{H}(G, G_{n_k}) \to \mathcal{H}(H, H_{m_k})\) and \(\mathcal{S}_k: \mathcal{H}(H, H_{m_k}) \to \mathcal{H}(G, G_{n_{k+1}})\) such that
		\[
			\mathcal{S}_k\mathcal{T}_k = \mathcal{Q}_{n_k, n_{k+1}}^G \quad\text{and}\quad \mathcal{T}_{k+1}\mathcal{S}_k = \mathcal{Q}_{m_k, m_{k+1}}^H
		\]
		for all \(k \in \mathbb{N}_0\).
	\end{enumerate}
	When following the constructions in the proof of (ii) implies (i) the diagram
	\begin{equation}\label{dia-Zig-zag}
		\begin{tikzcd}
			\mathcal{H}(\Phi) \arrow[rrr, "\cong\ \mathcal{T}"] & & &\mathcal{H}(\Psi)\\
			&\mathcal{H}(G, G_{n_{k+1}}) \arrow[lu, hook', swap, "\mathcal{Q}_{n_{k+1}}^G"] \arrow[r, hook, "\mathcal{T}_{k+1}"] &\mathcal{H}(H, H_{m_{k+1}}) \arrow[ru, hook, "\mathcal{Q}_{m_{k+1}}^H"]\\
			&\mathcal{H}(G, G_{n_k}) \arrow[luu, bend left, hook', "\mathcal{Q}_{n_k}^G"] \arrow[u, hook, swap, "\mathcal{Q}_{n_k, n_{k+1}}^G"] \arrow[r, hook, "\mathcal{T}_k"] &\mathcal{H}(H, H_{m_k}) \arrow[lu, hook', swap, "\mathcal{S}_k"] \arrow[u, swap, hook', "\mathcal{Q}_{m_k, m_{k+1}}^H"] \arrow[ruu, bend right, hook, "\mathcal{Q}_{m_k}^H"]
		\end{tikzcd}
	\end{equation}
	commutes for all \(k \in \mathbb{N}\). This is also true for the construction in the proof of (i) implies (ii) as long as \(\mathcal{T}\) is a Koopman operator.
\end{theorem}
\begin{proof*}
	First assume that (i) holds. By \cref{the-Isometric_Isomorphisms_Hecke_Algebras} we may assume that \(\mathcal{T}\) is a Koopman operator. We apply \cref{the-Group_Chains_Same_Index} to both \(\mathcal{T}\) and \(\mathcal{T}^{-1}\) to obtain group chains \((H_n')_{n \in \mathbb{N}_0}\) and \((G_n')_{n \in \mathbb{N}_0}\) in \(H\) and \(G\) respectively associated to the actions \((X, G, \Phi)\) and \((Y, H, \Psi)\) respectively. Further there are isometric *-algebra homomorphisms \(\mathcal{T}_n': \mathcal{H}(G, G_n) \to \mathcal{H}(H, H_n')\) and \(\mathcal{S}_n': \mathcal{H}(H, H_n) \to \mathcal{H}(G, G_n')\) such that
	\[
		\mathcal{Q}_n^{\prime H} \mathcal{T}_n' = \mathcal{T} \mathcal{Q}_n^G \quad\text{and}\quad \mathcal{Q}_n^{\prime G} \mathcal{S}_n' = \mathcal{T}^{-1} \mathcal{Q}_n^H
	\]
	for all \(n \in \mathbb{N}_0\).
	
	From \cref{the-Equivalence_Group_Chains} we obtain that the group chains \((G_n)_{n \in \mathbb{N}_0}\) and \((G_n')_{n \in \mathbb{N}_0}\) as well as \((H_n)_{n \in \mathbb{N}_0}\) and \((H_n')_{n \in \mathbb{N}_0}\) are equivalent, so we can construct strictly increasing sequences \((n_k)_{k \in \mathbb{N}_0}\) and \((m_k)_{k \in \mathbb{N}_0}\) such that \(H_{n_k}' \geq H_{m_k}\) and \(G_{m_k}' \geq G_{n_{k+1}}\). Then we we define \(\mathcal{T}_k: \mathcal{H}(G, G_{n_k}) \to \mathcal{H}(H, H_{m_k})\) and \(\mathcal{S}_k: \mathcal{H}(H, H_{m_k}) \to \mathcal{H}(G, G_{n_{k+1}})\) by
	\[
		\mathcal{T}_k := \widetilde{\mathcal{Q}}_k^H\mathcal{T}_{n_k}' \quad\text{and}\quad \mathcal{S}_k := \widetilde{\mathcal{Q}}_k^G\mathcal{S}_{m_k}'
	\]
	for all \(k \in \mathbb{N}\). Here \(\widetilde{\mathcal{Q}}_k^H\) and \(\widetilde{\mathcal{Q}}_k^G\) are the Koopman operators of the natural surjections \(\tilde{q}_k^H: H_{m_k}\backslash H/H_{m_k} \to H_{n_k}'\backslash H/H_{n_k}'\) and \(\tilde{q}_k^G: G_{n_{k+1}}\backslash G/G_{n_{k+1}} \to G_{m_k}'\backslash G/G_{m_k}'\). We summarize this definition in the following diagram:
	\[
		\begin{tikzcd}
			\mathcal{H}(\Phi) \arrow[rrr, "\cong\ \mathcal{T}"] & & &\mathcal{H}(\Psi)\\
			&\mathcal{H}(G, G_{n_{k+1}}) \arrow[lu, hook', swap, "\mathcal{Q}_{n_{k+1}}^G"]\\
			&\mathcal{H}(G, G_{m_k}') \arrow[u, hook, "\widetilde{\mathcal{Q}}_k^G"] \arrow[luu, hook', bend left, "\mathcal{Q}_{m_k}^{\prime G}"] &\mathcal{H}(H, H_{m_k}) \arrow[l, hook', swap, "\mathcal{S}_{m_k}'"] \arrow[lu, hook', swap, "\mathcal{S}_k"] \arrow[ruu, hook, "\mathcal{Q}_{m_k}^H"]\\
			&\mathcal{H}(G, G_{n_k}) \arrow[r, hook, "\mathcal{T}_{n_k}'"] \arrow[ru, hook, "\mathcal{T}_k"] \arrow[luuu, hook', bend left=60, "\mathcal{Q}_{n_k}^G"] &\mathcal{H}(H, H_{n_k}') \arrow[u, hook, "\widetilde{\mathcal{Q}}_k^H"] \arrow[ruuu, bend right, hook, swap, "\mathcal{Q}_{n_k}^{\prime H}"]
		\end{tikzcd}
	\]
	It commutes by the definition of \(\mathcal{T}_k\) and \(\mathcal{S}_k\), \cref{the-Group_Chains_Same_Index} and the equalities
	\[
		\mathcal{Q}_{m_k}^{\prime G} = \mathcal{Q}_{n_{k+1}}^G\widetilde{\mathcal{Q}}_k^G \quad\text{and}\quad \mathcal{Q}_{n_k}^{\prime H} = \mathcal{Q}_{m_k}^H\widetilde{\mathcal{Q}}_k^H,
	\]
	which come from \eqref{dia-Hecke_Algebras}.
	
	We now show that \(\mathcal{T}_k\) and \(\mathcal{S}_k\) satisfy the properties of assertion (ii). From this diagram above we obtain for \(k \in \mathbb{N}_0\) that
	\[
		\mathcal{Q}_{m_k}^H\mathcal{T}_k = \mathcal{T}\mathcal{Q}_{n_k}^G \quad\text{and}\quad \mathcal{Q}_{n_{k+1}}^G\mathcal{S}_k = \mathcal{T}^{-1}\mathcal{Q}_{m_k}^H.
	\]
	These two equation imply that we have
	\[
		\mathcal{Q}_{n_{k+1}}^G\mathcal{S}_k\mathcal{T}_k = \mathcal{T}^{-1}\mathcal{Q}_{m_k}^H\mathcal{T}_k = \mathcal{T}^{-1}\mathcal{T}\mathcal{Q}_{n_k}^G \overset{\eqref{dia-Hecke_Algebras}}{=} \mathcal{Q}_{n_{k+1}}^G\mathcal{Q}_{n_k, n_{k+1}}^G
	\]
	as well as
	\[
		\mathcal{Q}_{m_{k+1}}^H\mathcal{T}_{k+1}\mathcal{S}_k = \mathcal{T}\mathcal{Q}_{n_{k+1}}^G\mathcal{S}_k = \mathcal{T}\mathcal{T}^{-1}\mathcal{Q}_{m_k}^H \overset{\eqref{dia-Hecke_Algebras}}{=} \mathcal{Q}_{m_{k+1}}^H\mathcal{Q}_{m_k, m_{k+1}}^H.
	\]
	Since both \(\mathcal{Q}_{n_{k+1}}^G\) and \(\mathcal{Q}_{m_{k+1}}^H\) are injective, this implies
	\[
		\mathcal{S}_k\mathcal{T}_k = \mathcal{Q}_{n_k, n_{k+1}}^G \quad\text{and}\quad \mathcal{T}_{k+1}\mathcal{S}_k = \mathcal{Q}_{m_k, m_{k+1}}^H,
	\]
	which proves (ii).
	
	Now assume that assertion (ii) holds. We first define \(\mathcal{T}\) on the dense subspace \(\bigcup_{k \in \mathbb{N}_0} \operatorname{ran}(\mathcal{Q}_{n_k}^G)\) of \(\mathcal{H}(\Phi)\) (see \cref{lem-Images_Q_n_Dense}) and then extend it to \(\mathcal{H}(\Phi)\). A somewhat similar argument is used in the construction of the *-algebra structure of \(\mathcal{H}(\Phi)\) in \cref{the-*-Algebra_Structure_H}. Define \(\mathcal{T}: \bigcup_{k \in \mathbb{N}_0} \operatorname{ran}(\mathcal{Q}_{n_k}^G) \to \mathcal{H}(\Psi)\) by
	\begin{equation}\label{eq-definition_T}
		\mathcal{T}\mathcal{Q}_{n_k}^G\varphi_k := \mathcal{Q}_{m_k}^H\mathcal{T}_k\varphi_k
	\end{equation}
	for \(\varphi_k \in \mathcal{H}(G, G_{n_k})\). This is well defined since for \(\varphi_{k+1} \in \mathcal{H}(G, G_{n_{k+1}})\) and \(\varphi_k \in \mathcal{H}(G, G_{n_k})\) with \(\mathcal{Q}_{n_{k+1}}^G\varphi_{k+1} = \mathcal{Q}_{n_k}^G\varphi_k\) we have
	\[
		\mathcal{Q}_{n_{k+1}}^G\varphi_{k+1} = \mathcal{Q}_{n_k}^G\varphi_k \overset{\eqref{dia-Hecke_Algebras}}{=} \mathcal{Q}_{n_{k+1}}^G\mathcal{Q}_{n_k, n_{k+1}}^G\varphi_k,
	\]
	which implies \(\varphi_{k+1} = \mathcal{Q}_{n_k, n_{k+1}}^G\varphi_k\) since \(\mathcal{Q}_{n_{k+1}}^G\) is injective. Thus
	\begin{alignat*}{3}
		\mathcal{Q}_{m_{k+1}}^H\mathcal{T}_{k+1}\varphi_{k+1} &= \mathcal{Q}_{m_{k+1}}^H\mathcal{T}_{k+1}\mathcal{Q}_{n_k, n_{k+1}}^G\varphi_k\\
		&\overset{\text{(ii)}}{=} \mathcal{Q}_{m_{k+1}}^H\mathcal{T}_{k+1}\mathcal{S}_k\mathcal{T}_k\varphi_k\\
		&\overset{\text{(ii)}}{=} \mathcal{Q}_{m_{k+1}}^H\mathcal{Q}_{m_k, m_{k+1}}^H\mathcal{T}_k\varphi_k &\overset{\eqref{dia-Hecke_Algebras}}{=} \mathcal{Q}_{m_k}\mathcal{T}_k\varphi_k.
	\end{alignat*}
	This map \(\mathcal{T}\) is an isometric *-algebra homomorphism since all \(\mathcal{Q}_{n_k}^G\), \(\mathcal{Q}_{m_k}^H\) and \(\mathcal{T}_k\) are. Therefore it can be uniquely extended to \(\mathcal{H}(\Phi)\) by \cref{lem-Images_Q_n_Dense}. The diagram \eqref{dia-Zig-zag} commutes because of the relations given in (ii) and the construction of \(\mathcal{T}\).
	
	We need to prove that \(\mathcal{T}\) is an isometric isomorphism of Banach-*-algebras. We already know that it is an isometry and a *-algebra homomorphism, so it suffices to show that the dense subspace \(\bigcup_{k \in \mathbb{N}_0} \operatorname{ran}(\mathcal{Q}_{m_k}^H)\) (see \cref{lem-Images_Q_n_Dense}) is contained in its range. Fix \(k \in \mathbb{N}_0\) and \(\psi_k \in \mathcal{H}(H, H_{m_k})\). Then \(\mathcal{S}_k\psi_k \in \mathcal{H}(G, G_{n_{k+1}})\) and
	\[
		\mathcal{T}\mathcal{Q}_{n_{k+1}}^G\mathcal{S}_k\psi_k \overset{\eqref{eq-definition_T}}{=} \mathcal{Q}_{m_{k+1}}^H\mathcal{T}_{k+1}\mathcal{S}_k\psi_k \overset{\text{(ii)}}{=} \mathcal{Q}_{m_{k+1}}^H\mathcal{Q}_{m_k, m_{k+1}}^H\psi_k \overset{\eqref{dia-Hecke_Algebras}}{=} \mathcal{Q}_{m_k}^H\psi_k.
	\]
	Thus the range of \(\mathcal{T}\) contains a dense subspace. So \(\mathcal{T}\) is indeed an isometric *-algebra isomorphisms.
\end{proof*}

%% file: ch-Examples.tex
\section{Examples}\label{sec-Examples}

In this last section we consider three examples. In the first part we look at the infinite dihedral group \(H\) in the second part at the (first) Grigorchuk group \(G\) and the group \(\widetilde{G}\) constructed by Bartholdi and Grigorchuk. The actions of these groups are actions on the binary rooted tree, so all three examples have the same Steinitz order, namely \(2^\infty\).

For the dihedral group we calculate the finite dimensional Hecke algebras associated to its group chain. This is done by identifying the double cosets of the group chain and direct use of the definition of the convolution and involution. For the other two groups the double cosets of their respective group chains were already computed in \citep{bartholdi_parabolic_2002}. Both groups \(G\) and \(\widetilde{G}\) have the same double coset structure. From these structures we compute the double coset space of their Ellis groups and their discriminant groups using the inverse limit in \cref{pro-Inverse_Limit_Double_Cosets}. Similarly, the convolution and involution of the infinite dimensional Hecke algebra can be described using the approximation by finite Hecke algebras and the double coset structure. In particular they have isometric isomorphic Hecke algebras and we do not require any knowledge of the groups beyond their double coset structure to obtain the Hecke algebra of their actions.

With these examples we see that the actions of distinct groups (\(G\) and \(\widetilde{G}\)) have isometric isomorphic Hecke algebras. The actions of \(G\) and \(H\) have non-isomorphic Hecke algebras (\cref{pro-Grigorchuk_Group_Dihedral_Group}) despite having the same Steinitz order.

\subsection{Infinite Dihedral Group}\label{sub-Infinite_Dihedral_Group}

We first consider the example of the infinite dihedral group
\[
	H := \langle a, b: b^2 = abab = e \rangle.
\]
In this group we take the group chain \((H_n)_{n \in \mathbb{N}_0}\) given by the subgroups
\[
	H_n := \langle a^{2^n}, b \rangle \leq H.
\]
The index of \(H_n\) is \(2^n\) and
\[
	H/H_n = \{a^iH_n: i \in \{0, \dots, 2^n-1\}\}.
\]
We note that the minimal equicontinuous Cantor action arising from this group chain is the the iterated monodromy group of the quadratic Chebyshev polynomial (\citep[Proposition 5.6.]{grabner_fractal_2003}).

For the double cosets we note that \(a^{2^n}a^iH_n = a^iH_n\) and \(ba^iH_n = a^{-i}H_n = a^{2^n-i}H_n\), so
\[
	H_n(a^iH_n) = \{a^iH_n, a^{2^n-i}H_n\} =
	\begin{dcases}
		\{a^iH_n\} &: i \in \{0, 2^{n-1}\}\\
		\{a^iH_n, a^{2^n-i}H_n\} &: \text{otherwise}
	\end{dcases},
\]
and
\[
	H_n\backslash H/H_n = \{H_na^iH_n: i \in \{0, \dots, 2^{n-1}\}\},
\]
where each double coset consists of two cosets except if \(i = 0\) or \(i = 2^{n-1}\), in which case it consists of a single coset.

Now we calculate the convolution and the involution for characteristic functions of cosets in \(H/H_n\). It suffices to look at those functions as they form a basis of \(\mathcal{H}(H, H_n)\). This is similar to the Hecke algebras for finite dihedral groups calculated in \citep[Section III.1]{krieg_hecke_1990}.

\begin{theorem}\label{the-Hecke_Algebra_Dihedral_Group}
	For \(i, j \in \{0, \dots, 2^{n-1}\}\) we have
	\[
		1_{\{H_na^iH_n\}} * 1_{\{H_na^jH_n\}} = \frac{1}{2^n} \cdot 1_{\{H_na^{i+j}H_n, H_na^{i-j}H_n\}}
	\]
	and \(1_{\{H_na^iH_n\}}^* = 1_{\{H_na^iH_n\}}\). In particular \(\mathcal{H}(H, H_n)\) is commutative.
\end{theorem}
\begin{proof*}
	The formula for the involution follows from \(H_na^iH_n = H_na^{-i}H_n\). Fix \(i, j \in \{0, \dots, 2^{n-1}\}\). If \(i \in \{0, 2^{n-1}\}\) then
	\begin{align*}
		\big(1_{\{H_na^iH_n\}} * 1_{\{H_na^jH_n\}}\big)(H_na^kH_n) &= \frac{1}{2^n} \sum_{l=0}^{2^n-1} 1_{\{H_na^iH_n\}}(H_na^lH_n) \cdot 1_{\{H_na^jH_n\}}(H_na^{-l}a^kH_n)\\
		&= \frac{1}{2^n} \cdot 1_{\{H_na^jH_n\}}(H_na^{k-i}H_n)\\
		&= \frac{1}{2^n} \cdot 1_{\{H_na^{i-j}H_n\}}(H_na^kH_n)
	\end{align*}
	for all \(k \in \{0, \dots, 2^{n-1}\}\). This agrees with the formulation of the theorem since in this case \(H_na^{i-j}H_n = H_na^{-i+j}H_n = H_na^{i+j}H_n\). For the case \(i \notin \{0, 2^{n-1}\}\) we have for \(k \in \mathbb{Z}\):
	\begin{align*}
		&\big(1_{\{H_na^iH_n\}} * 1_{\{H_na^jH_n\}}\big)(H_na^kH_n)\\
		= &\frac{1}{2^n} \sum_{l=0}^{2^n-1} 1_{\{H_na^iH_n\}}(H_na^lH_n) \cdot 1_{\{H_na^jH_n\}}(H_na^{-l}a^kH_n)\\
		= &\frac{1}{2^n} \cdot \big( 1_{\{H_na^jH_n\}}(H_na^{k-i}H_n) + 1_{\{H_na^jH_n\}}(H_na^{k-(2^n-i)}H_n) \big)\\
		= &\frac{1}{2^n} \cdot \big( 1_{\{H_na^jH_n\}}(H_na^{k-i}H_n) + 1_{\{H_na^jH_n\}}(H_na^{k+i}H_n) \big)\\
		= &\frac{1}{2^n} \cdot 1_{\{H_na^{i+j}H_n, H_na^{i-j}H_n\}}(H_na^kH_n).
	\end{align*}
	The last step requires some explanation: For distinct \(k, k' \in \{0, \dots, 2^{n-1}\}\) we have \(H_na^{k-i}H_n \neq H_na^{k'-i}H_n\) and \(H_na^{k+i}H_n \neq H_na^{k'+i}H_n\), so either of the summands in the second to last line is non-zero for exactly one \(k \in \{0, \dots, 2^{n-1}\}\). Also the first summand is non-zero for \(k'' = i + j\) and the second summand for \(k'' = i-j\). Thus the \(k \in \{0, \dots, 2^{n-1}\}\) for which either of the summands is non-zero are those with \(H_na^kH_n = H_na^{i+j}H_n\) and \(H_na^kH_n = H_na^{i-j}H_n\) respectively. Further with \(i \notin \{0, 2^{n-1}\}\) we always have \(H_na^{i+j}H_n \neq H_na^{i-j}H_n\).
\end{proof*}

We now turn to the infinite dimensional Hecke algebra. We denote by \(\widehat{H}\) the Ellis group and by \(E\) the discriminant group associated to the group chain \((H_n)_{n \in \mathbb{N}_0}\). We know that \(\widehat{H}\) is a Cantor group. It is straightforward to check that
\[
	C_n = \operatorname{core}_H(H_n) = \langle a^{2^n} \rangle
\]
for all \(n \in \mathbb{N}_0\) and thus \(H_n/C_n = \{e, b\}\) has two elements. Therefore \(E \overset{\ref{the-Inverse_Limit_Ellis_Group}}{\cong} \varprojlim_{n \in \mathbb{N}_0} H_n/C_n\) has two elements as well. Since \(\widehat{H}\) is a cantor group and \(E\) is finite we have that \(E\backslash\widehat{H}/E\) is a Cantor set. It is difficult to give a good explicit description of the convolution in \(\mathcal{H}(\widehat{H}, E)\). However, the approximation by the finite dimensional Hecke algebras \(\mathcal{H}(H, H_n)\) as used in the proof of \cref{the-*-Algebra_Structure_H} yields with \cref{the-Hecke_Algebra_Dihedral_Group} a description of the convolution on a dense subspace. In particular, since all the \(\mathcal{H}(H, H_n)\) are commutative, it follows that \(\mathcal{H}(\widehat{H}, E)\) is commutative as well. Regarding the involution we have from \cref{the-Hecke_Algebra_Dihedral_Group} that the involution on \(\mathcal{H}(H, H_n)\) is given by pointwise complex conjugation. Thus the approximation by the \(\mathcal{H}(H, H_n)\) yields that the involution on \(\mathcal{H}(\widehat{H}, E)\) is also given by pointwise complex conjugation.

\subsection{Grigorchuk Group}\label{sub-Grigotchuk's_Groups}

In this subsection we calculate the Hecke algebra for two more groups. The first one is the (first) Grigorchuk group \(G\). It was first constructed in \citep{grigorchuk_bernsides_1980} and discussed in various other works including \citep[Section 4.]{bartholdi_parabolic_2002}. The second group \(\widetilde{G}\) was introduced in \citep[Section 5.]{bartholdi_parabolic_2002} and is called the \qq{Grigorchuk supergroup} in \citep[Section 4.3.]{bartholdi_branch_2003}. It contains \(G\) as an infinite index subgroup (\citep[Proposition 5.2.]{bartholdi_parabolic_2002}). Both groups are branch groups and have intermediate growth. However while \(G\) is a \(2\)-group, we find infinite order elements in \(\widetilde{G}\).

These groups are constructed as automorphisms of the binary tree, so they give rise to a minimal equicontinuous action on the boundary of the tree, which is a Cantor set. Group chains \((G_n)_{n \in \mathbb{N}_0}\) and \((\widetilde{G}_n)_{n \in \mathbb{N}_0}\) associated to these actions are obtained by letting \(G_n\) and \(\widetilde{G}_n\) be the stabilizers of the rightmost vertex in the \(n\)-th level of the binary tree. These are denoted by \(P_n\) and \(\widetilde{P}_n\) in \citep{bartholdi_parabolic_2002} and have index \(2^n\). We do not give an explicit definition here, since this requires formalism not introduced in this article. Also, we do not need to know the precise definition of the groups \(G\) and \(\widetilde{G}\) to obtain their Hecke algebras, it suffices to know the structure of the double coset spaces \(G_n\backslash G/G_n\) and \(\widetilde{G}_n\backslash\widetilde{G}/\widetilde{G}_n\). This was already computed in the article by Bartholdi and Grigorchuk (\citep[Lemmata 9.10. and 9.11.]{bartholdi_parabolic_2002}). Further they prove that the arising Hecke algebras \(\mathcal{H}(G, G_n)\) and \(\mathcal{H}(\widetilde{G}, \widetilde{G}_n)\) are commutative (\citep[Proposition 9.17.]{bartholdi_parabolic_2002}). For the double coset spaces they show that \(G_n\backslash G/G_n\) consists of \(n+1\) distinct double cosets. We denote representatives of those by \(g_{0, n}, \dots, g_{n, n}\), where we choose \(g_{n, n}\) to be the identity element. Following \citep[Lemma 9.10.]{bartholdi_parabolic_2002} we can choose the others such that \(G_n(g_{i, n}G_n)\) contains \(2^{n-1-i}\) left cosets for \(i \in \{0, \dots, n-1\}\). Thus
\begin{equation}\label{eq-Grigorchuk_nu_n}
	\nu_n(\{G_ng_{i, n}G_n\}) \overset{\eqref{eq-nu_n}}{=} \frac{|G_n(g_{i, n}G_n)|}{|G/G_n|} =
	\begin{dcases}
		\frac{1}{2^{i+1}} &: i < n\\
		\frac{1}{2^n} &: i = n
	\end{dcases}.
\end{equation}
The maps \(q_{n, n +1}: G_{n+1}\backslash G/G_{n+1} \to G_n\backslash G/G_n\) are then given by
\[
	q_{n, n+1}(G_{n+1}g_{i, n+1}G_{n+1}) = G_ng_{\min\{i, n\}, n}G_n.
\]
One can check that this map is the only surjective measure preserving map between those spaces that maps \(G_{n+1}g_{n+1, n+1}G_{n+1}\) to \(G_ng_{n, n}G_n\). From this formula we obtain that the Koopman operators \(\mathcal{Q}_{n, n+1}: \mathcal{H}(G, G_n) \to \mathcal{H}(G, G_{n+1})\) are given by
\begin{equation}\label{eq-Grigorchuk_Q_n_n+1}
	\mathcal{Q}_{n, n+1}1_{\{G_ng_{i, n}G_n\}} =
	\begin{dcases}
		1_{\{G_ng_{i, n+1}G_n\}} &: i \in \{0, \dots, n-1\}\\
		1_{\{G_ng_{n, n+1}G_n\}} + 1_{\{G_ng_{n+1, n+1}G_n\}} &: i = n
	\end{dcases}.
\end{equation}
The double cosets spaces \(\widetilde{G}_n\backslash\widetilde{G}/\widetilde{G}_n\) as well as the maps \(\tilde{q}_{n, n+1}: \widetilde{G}_{n+1}\backslash\widetilde{G}/\widetilde{G}_{n+1} \to \widetilde{G}_n\backslash\widetilde{G}/\widetilde{G}_n\) and their Koopman operators \(\widetilde{\mathcal{Q}}_n: \mathcal{H}(\widetilde{G}, \widetilde{G}_n) \to \mathcal{H}(\widetilde{G}, \widetilde{G}_n)\) are given analogously.

Now we can already prove that the Hecke algebras of the actions of \(G\) and \(\widetilde{G}\) are not isometrically isomorphic to the Hecke algebra of the action of the infinite dihedral group \(H\) from the previous section. Denote by \(\widehat{G}\) the Ellis group arising from the group chain \((G_n)_{n \in \mathbb{N}_0}\) of the Grigorchuk group \(G\) and by \(D\) the corresponding discriminant group.

\begin{proposition}\label{pro-Grigorchuk_Group_Dihedral_Group}
	The Hecke algebras \(\mathcal{H}(\widehat{H}, E)\) and \(\mathcal{H}(\widehat{G}, D)\) of the Cantor actions of the infinite dihedral group \(H\) introduced in \cref{sub-Infinite_Dihedral_Group} and of the Grigorchuk group \(G\) introduced above are not isometric isomorphic. The same is true for the group \(\widetilde{G}\) instead of the Grigorchuk group.
\end{proposition}
\begin{proof*}
	Assume that there exists an isometric *-algebra isomorphism \(\mathcal{T}: \mathcal{H}(\widehat{H}, E) \to \mathcal{H}(\widehat{G}, D)\). By \cref{the-Isometric_Isomorphisms_Hecke_Algebras} we may assume that \(\mathcal{T}\) is the Koopman operator of some measure preserving homeomorphism \(\theta: D\backslash\widehat{G}/D \to E\backslash\widehat{H}/E\). From \cref{the-Hecke_Algebra_Zig-zag_Map} we find some \(n > 2\) and some \(m \in \mathbb{N}\) with an isometric *-algebra homomorphism \(\mathcal{T}_n: \mathcal{H}(H, H_n) \to \mathcal{H}(G, G_m)\) such that
	\[
		\mathcal{T}\mathcal{Q}_n^H = \mathcal{Q}_m^G\mathcal{T}_n.
	\]
	Since \(\mathcal{T}\), \(\mathcal{Q}_n^H\) and \(\mathcal{Q}_m^G\) are Koopman operators of measure preserving maps, also \(\mathcal{T}_n\) is the Koopman operator of some measure preserving map \(\theta_n: G_m\backslash G/G_m \to H_n\backslash H/H_n\). We have
	\[
		\frac{1}{2} \overset{\eqref{eq-Grigorchuk_nu_n}}{=} \nu_m^G(\{G_mg_{0, m}G_m\}) \leq \nu_m^G(\theta_n^{-1}(\{\theta(G_mg_{0, m}G_m)\})) = \nu_n^H(\{\theta_n(G_ng_{0, m}G_m)\}) \overset{\eqref{eq-nu_n}}{\leq} \frac{2}{2^n},
	\]
	which is a contradiction. Here we used that for the infinite dihedral group each double coset in \(H_m\backslash H/H_m\) consists of at most two cosets. The prove for \(\widetilde{G}\) works the same way, as we do not use any properties of \(G\) which do not hold for \(\widetilde{G}\).
\end{proof*}

We want to give an explicit description of the infinite dimensional Hecke algebra \(\mathcal{H}(\widehat{G}, D)\). For the group \(\widetilde{G}\) everything can be done analogously. We first describe the double coset space \(D\backslash\widehat{G}/D\). For this we use that \(D\backslash\widehat{G}/D \cong \varprojlim_{n \in \mathbb{N}_0} G_n\backslash G/G_n\) by \cref{pro-Inverse_Limit_Double_Cosets}. We will prove that this inverse limit can be identified with the one point compactification \(\mathbb{N}_0 \cup \{\infty\}\) of \(\mathbb{N}_0\). In particular this yields an example where the double coset space \(D\backslash\widehat{G}/D\) is not a Cantor set. Since for the dihedral group \(E\backslash\widehat{H}/E\) is a Cantors set, \cref{pro-Double_Coset_Grigorchuk} together with \cref{the-Isometric_Isomorphisms_Hecke_Algebras} reconfirms \cref{pro-Grigorchuk_Group_Dihedral_Group}.

\begin{proposition}\label{pro-Double_Coset_Grigorchuk}
	Let \(G\) be the Grigorchuk group and denote by \(\widehat{G}\) and \(D\) the Ellis and discriminant group corresponding to its action on the binary tree. The space \(D\backslash\widehat{G}/D\) can be identified with \(\mathbb{N}_0 \cup \{\infty\}\). With this identification the maps \(q_n: \mathbb{N}_0 \cup \{\infty\} \to G_n\backslash G/G_n\) are given by
	\[
	q_n(i) = G_ng_{\min\{i, n\}, n}G_n
	\]
	for all \(i \in \mathbb{N}_0 \cup \{\infty\}\) and \(n \in \mathbb{N}_0\). The probability measure \(\nu\) on \(\mathbb{N}_0 \cup \{\infty\}\) is given by
	\[
	\nu(\{i\}) = \frac{1}{2^{i+1}}
	\]
	for all \(i \in \mathbb{N}_0 \cup \{\infty\}\) with convention \(2^{-\infty} = 0\). The same holds for \(\widetilde{G}\).
\end{proposition}
\begin{proof*}
	From \cref{pro-Inverse_Limit_Double_Cosets} we know
	\[
		D\backslash\widehat{G}/D \cong Z :=  \varprojlim_{n \in \mathbb{N}_0} G_n\backslash G/G_n = \{(G_kg_kG_k)_{k \in \mathbb{N}_0}: G_kg_kG_k = G_kg_{k+1}G_k\ \forall k \in \mathbb{N}_0\}.
	\]
	We define the map \(\theta: \mathbb{N}_0 \cup \{\infty\} \to Z\) by
	\[
	\theta(i) := (G_kg_{\min\{i, k\}, k}G_k)_{k \in \mathbb{N}_0}
	\]
	for all \(i \in \mathbb{N}_0 \cup \{\infty\}\). This is well defined since if \(i \leq k\) then
	\[
		G_kg_{i, k+1}G_k = q_{k, k+1}(G_{k+1}g_{i, k+1}G_{k+1}) = G_kg_{i, k}G_k
	\]
	and \(G_kg_{k+1, k+1}G_k = G_keG_k = G_kg_{k, k}G_k\) if \(i \geq k+1\) for all \(k \in \mathbb{N}_0\). The map \(\theta\) is continuous since
	\[
	\lim_{i\to\infty} \theta(i) = \lim_{i\to\infty} (G_kg_{\min\{i, k\}, k}G_k)_{k \in \mathbb{N}_0} = (G_kg_{k, k}G_k)_{k \in \mathbb{N}_0} = \theta(\infty).
	\]
	It is surjective by \citep[Corollary 1.1.6]{ribes_profinite_2000} since \(D\backslash\widehat{G}/D\) is compact Hausdorff and the map \(\theta\) is induced by the continuous surjective maps \(i \mapsto G_kg_{\min\{i, k\}, k}G_k\) for all \(k \in \mathbb{N}_0\). For distinct \(i, j \in \mathbb{N}_0 \cup \{\infty\}\) we have \(G_kg_{\min\{i, k\}, k}G_k \neq G_kg_{\min\{j, k\}, k}G_k\) for all \(k > \min\{i, j\}\), so \(\theta\) is also injective. Thus \(\theta\) is indeed a homeomorphism since \(\mathbb{N}_0 \cup \{\infty\}\) is compact Hausdorff.
	
	Identifying \(D\backslash\widehat{G}/D\) with \(\mathbb{N}_0 \cup \{\infty\}\), it is clear from the definition of \(\theta\) that the maps \(q_n: D\backslash\widehat{G}/D \to G_n\backslash G/G_n\) are given by the maps
	\[
		q_n(i) = G_ng_{\min\{i, n\}, n}G_n
	\]
	for all \(i \in \mathbb{N}_0 \cup \{\infty\}\) and \(n \in \mathbb{N}_0\). For the measure we note that for \(i \in \mathbb{N}_0\)
	\[
		\nu(\{i\}) = \nu(q_{i+1}^{-1}\{G_{i+1}g_{i, i+1}G_{i+1}\}) = \nu_{i+1}(\{G_{i+1}g_{i, i+1}G_{i+1}\}) \overset{\eqref{eq-Grigorchuk_nu_n}}{=} \frac{1}{2^{i+1}}
	\]
	and for \(i = \infty\) we have
	\[
		\nu(\{\infty\}) = 1 - \nu(\mathbb{N}_0) = 1 - \sum_{j\in\mathbb{N}_0} \frac{1}{2^{j+1}} = 0.
	\]
	All of this works identically for \(\widetilde{G}\).
\end{proof*}

Identifying \(D\backslash\widehat{G}/D\) with \(\mathbb{N}_0 \cup \{\infty\}\), we also identify \(\mathcal{H}(\widehat{G}, D)\) with \(\operatorname{L}^1(\mathbb{N}_0 \cup \{\infty\}, \nu)\). From the formula for \(q_n\) is follows that its Koopman operator \(\mathcal{Q}_n: \mathcal{H}(G, G_n) \to \mathcal{H}(\widehat{G}, D)\) is given by
\begin{equation}\label{eq-Grigorchuk_Q_n}
	\mathcal{Q}_n1_{\{G_ng_{i, n}G_n\}} =
	\begin{dcases}
		1_{\{i\}} &: i < n\\
		1_{\{i, i+1, \dots,\} \cup \{\infty\}} &: i = n
	\end{dcases}
\end{equation}
for all \(i \in \{0, \dots, n\}\) and \(n \in \mathbb{N}\). We use this to obtain explicit formulae for the convolution and involution in Hecke algebra \(\mathcal{H}(\widehat{G}, D)\). We obtain the same formulae for the group \(\widetilde{G}\), so it is clear that the actions of \(G\) and \(\widetilde{G}\) have isometric isomorphic Hecke algebras.

\begin{theorem}\label{the-Hecke_Algebra_Grigorchuk's_Group}
	Let \(G\) be the Grigorchuk group and denote by \(\widehat{G}\) and \(D\) the Ellis and discriminant group corresponding to its action on the binary tree. We identify \(D\backslash\widehat{G}/D\) with \(\mathbb{N}_0 \cup \{\infty\}\) as described in \cref{pro-Double_Coset_Grigorchuk}. Then \(\varphi^* = \overline{\varphi}\) for all \(\varphi \in \mathcal{H}(\widehat{G}, D)\) and
	\[
		(\varphi * \psi)(k) = \sum_{i=0}^{k-1} \frac{1}{2^{i+1}} \cdot \varphi(i) \cdot \psi(i) + \sum_{i = k+1}^\infty \frac{1}{2^{i+1}} \cdot (\varphi(i) \cdot \psi(k) + \varphi(k) \cdot \psi(i))
	\]
	for all \(\varphi, \psi \in \mathcal{H}(\widehat{G}, D)\) and \(k \in \mathbb{N}_0\) (note that \(\{\infty\}\) is a measure zero set). Here we use the convention, that for the case \(k = 0\) the first sum is set to zero. The same formulae hold for the group \(\widetilde{G}\).
\end{theorem}

To prove this we first require the following lemma concerning the convolution and involution for the Hecke algebra \(\mathcal{H}(G, G_n)\).

\begin{lemma}\label{lem-Hecke_Algebra_Grigorchuk}
	For \(n \in \mathbb{N}\) and \(j \in \{0, \dots, n-1\}\) we have
	\[
		1_{\{G_ng_{n-1, n}G_n\}} * 1_{\{G_ng_{j, n}G_n\}} = 1_{\{G_ng_{j, n}G_n\}} * 1_{\{G_ng_{n-1, n}G_n\}} =
		\begin{dcases}
			\frac{1}{2^n} \cdot 1_{\{G_ng_{n, n}G_n\}} &: j = n-1\\
			\frac{1}{2^n} \cdot 1_{\{G_ng_{j, n}G_n\}} &: j < n-1
		\end{dcases}
	\]
	and \(1_{\{G_ng_{n-1, n}G_n\}}^* = 1_{\{G_ng_{n-1, n}G_n\}}\).
\end{lemma}
\begin{proof*}
	Since \(g_{n, n} = e\) we have
	\[
		1_{\{G_ng_{n, n}G_n\}}^* = \overline{1_{\{G_ng_{n, n}^{-1}G_n\}}} = 1_{\{G_ng_{n, n}G_n\}}.
	\]
	From the definition of the involution we know that \(1_{\{G_ng_{n-1, n}G_n\}}^*\) is the characteristic function of some double coset. Since the involution is measure preserving (\cref{pro-H_n_Banach-*-algebra}) this double coset must have the same measure as \(G_ng_{n-1, n}G_n\). Therefore by \eqref{eq-Grigorchuk_nu_n} this double coset must be either \(G_ng_{n-1, n}G_n\) or \(G_ng_{n, n}G_n\). The latter would be a contradiction to the involution being self-inverse with the calculation above. Thus \(1_{\{G_ng_{n-1, n}G_n\}}^* = 1_{\{G_ng_{n-1, n}G_n\}}\), in particular \(G_ng_{n-1, n}^{-1}G_n = G_ng_{n-1, n}G_n\).
	
	Now we prove the formula for the convolution. Fix \(n \in \mathbb{N}\). In the following proof we repeatably use that \(2^{n-1} \cdot 1_{\{G_{n-1}g_{n-1, n-1}G_{n-1}\}}\) and \(2^n \cdot 1_{\{G_ng_{n, n}G_n\}}\) are the unit elements in their respective Hecke algebras. For \(j = n-1\) we have
	\begin{align*}
		&1_{\{G_ng_{j, n}G_n\}} * 1_{\{G_ng_{j, n}G_n\}}\\
		\overset{\eqref{eq-Grigorchuk_Q_n_n+1}}{=} &\big(\mathcal{Q}_{j, n}1_{\{G_jg_{j, j}G_j\}} - 1_{\{G_ng_{n, n}G_n\}}\big) * \big(\mathcal{Q}_{j, n}1_{\{G_jg_{j, j}G_j\}} - 1_{\{G_ng_{n, n}G_n\}}\big)\\
		= &\mathcal{Q}_{j, n}\big(1_{\{G_jg_{j, j}G_j\}} * 1_{\{G_jg_{j, j}G_j\}}\big) - \frac{2}{2^n} \cdot \mathcal{Q}_{j, n}1_{\{G_jg_{j, j}G_j\}} + \frac{1}{2^n} \cdot 1_{\{G_ng_{n, n}G_n\}}\\
		= &\frac{1}{2^j} \cdot \mathcal{Q}_{j, n}1_{\{G_jg_{j, j}G_j\}} - \frac{1}{2^j} \cdot \mathcal{Q}_{j, n}1_{\{G_jg_{j, j}G_j\}} + \frac{1}{2^n} \cdot 1_{\{G_ng_{n, n}G_n\}}\\
		= &\frac{1}{2^n} \cdot 1_{\{G_ng_{n, n}G_n\}}.
	\end{align*}
	For \(j < n - 1\) we have
	\begin{align*}
		&1_{\{G_ng_{n-1, n}G_n\}} * 1_{\{G_ng_{j, n}G_n\}}\\
		\overset{\eqref{eq-Grigorchuk_Q_n_n+1}}{=} &\big(\mathcal{Q}_{n-1, n}1_{\{G_{n-1}g_{n-1, n-1}G_{n-1}\}} - 1_{\{G_ng_{n, n}G_n\}}\big) * \mathcal{Q}_{n-1, n}1_{\{G_{n-1}g_{j, n-1}G_{n-1}\}}\\
		= &\mathcal{Q}_{n-1, n}\big(1_{\{G_{n-1}g_{n-1, n-1}G_{n-1}\}} * 1_{\{G_{n-1}g_{j, n-1}G_{n-1}\}}\big) - \big(1_{\{G_ng_{n, n}G_n\}} * \mathcal{Q}_{n-1, n}1_{\{G_{n-1}g_{j, n-1}G_{n-1}\}}\big)\\
		= &\frac{1}{2^{n-1}} \cdot \mathcal{Q}_{n-1, n}1_{\{G_{n-1}g_{j, n-1}G_{n-1}\}} - \frac{1}{2^n} \cdot \mathcal{Q}_{n-1, n}1_{\{G_{n-1}g_{j, n-1}G_{n-1}\}}\\
		= &\frac{1}{2^n} \cdot \mathcal{Q}_{n-1, n}1_{\{G_{n-1}g_{j, n-1}G_{n-1}\}}\\
		\overset{\eqref{eq-Grigorchuk_Q_n_n+1}}{=} &\frac{1}{2^n} \cdot 1_{\{G_ng_{j, n}G_n\}}.
	\end{align*}
	The calculation for \(1_{\{G_ng_{j, n}\}G_n} * 1_{\{G_ng_{n-1, n}G_n\}}\) is similar.
\end{proof*}

\begin{proof*}[\cref{the-Hecke_Algebra_Grigorchuk's_Group}]
	We use that for \(\varphi \in \mathcal{H}(\widehat{G}, D)\) we have
	\[
		\varphi = \sum_{i\in\mathbb{N}_0} \varphi(i) \cdot 1_{\{i\}}.
	\]
	Here we use \(\nu(\{\infty\}) = 0\). Since the involution is antilinear and continuous, we have that
	\[
		\varphi^* = \sum_{i\in\mathbb{N}_0} \overline{\varphi(i)} 1_{\{i\}}^* \overset{\eqref{eq-Grigorchuk_Q_n}}{=} \sum_{i\in\mathbb{N}_0} \overline{\varphi(i)} \mathcal{Q}_{i+1}1_{\{G_{i+1}g_{i, i+1}G_{i+1}\}}^* \overset{\ref{lem-Hecke_Algebra_Grigorchuk}}{=} \sum_{i\in\mathbb{N}_0} \overline{\varphi(i)} \mathcal{Q}_{i+1}1_{\{G_{i+1}g_{i, i+1}G_{i+1}\}} \overset{\eqref{eq-Grigorchuk_Q_n}}{=} \overline{\varphi}.
	\]
	For the convolution we have for \(i \in \mathbb{N}_0\) that
	\begin{align*}
		1_{\{i\}} * 1_{\{i\}} &\overset{\eqref{eq-Grigorchuk_Q_n}}{=} \mathcal{Q}_{i+1}(1_{\{G_{i+1}g_{i, i+1}G_{i+1}\}} * 1_{\{G_{i+1}g_{i, i+1}G_{i+1}\}})\\
		&\overset{\ref{lem-Hecke_Algebra_Grigorchuk}}{=} \frac{1}{2^{i+1}} \cdot \mathcal{Q}_{i+1}1_{\{G_{i+1}g_{i+1, i+1}G_{i+1}\}}\\
		&\overset{\eqref{eq-Grigorchuk_Q_n}}{=} \frac{1}{2^{i+1}} \cdot 1_{\{i+1, i+2, \dots\} \cup \{\infty\}}.
	\end{align*}	
	For \(i > j \geq 0\) we have
	\begin{alignat*}{3}
		1_{\{i\}} * 1_{\{j\}} &\overset{\eqref{eq-Grigorchuk_Q_n}}{=} \mathcal{Q}_{i+1}(1_{\{G_{i+1}g_{i, i+1}G_{i+1}\}} * 1_{\{G_{i+1}g_{j, i+1}G_{i+1}\}})\\
		&\overset{\ref{lem-Hecke_Algebra_Grigorchuk}}{=} \frac{1}{2^{i+1}} \cdot \mathcal{Q}_{i+1}1_{\{G_{i+1}g_{j, i+1}G_{i+1}\}} &\overset{\eqref{eq-Grigorchuk_Q_n}}{=} \frac{1}{2^{i+1}} \cdot 1_{\{j\}}.
	\end{alignat*}
	We obtain a similar calculation for the case \(j > i \geq 0\), such that
	\[
		1_{\{i\}} * 1_{\{j\}} = \frac{1}{2^{\max\{i, j\}+1}} \cdot 1_{\{\min\{i, j\}\}}
	\]
	for all \(i, j \in \mathbb{N}_0\) with \(i \neq j\). Therefore for \(\varphi, \psi \in \mathcal{H}(\widehat{G}, D)\) and \(k \in \mathbb{N}_0\) we have, using that the convolution is bilinear and continuous, that
	\begin{align*}
		(\varphi * \psi)(k) &= \sum_{i\in\mathbb{N}_0} \sum_{j\in\mathbb{N}_0} \varphi(i) \cdot \psi(j) \cdot (1_{\{i\}} * 1_{\{j\}})(k)\\
		&= \sum_{i=0}^{k-1} \frac{1}{2^{i+1}} \cdot \varphi(i) \cdot \psi(i) + \sum_{i = k+1}^\infty \frac{1}{2^{i+1}} \cdot (\varphi(i) \cdot \psi(k) + \varphi(k) \cdot \psi(i))
	\end{align*}
	The two sums in the last line are the distinct cases in which \((1_{\{i\}} * 1_{\{j\}})(k) \neq 0\), namely \(i = j < k\) for the first sum and \(i > j = k\) as well as \(j > i = k\) for the second sum. The same calculation works for \(\widetilde{G}\).
\end{proof*}

%% file: ch-Appendix.tex
\appendix
\section{The Banach Lamperti Theorem}\label{sec-Appendix}

In this appendix we state a corollary of the Banach-Lamperti theorem. Let \((X, \mathcal{A}, \mu)\) and \((Y, \mathcal{B}, \nu)\) be two standard probability spaces. We will show that any isometric isomorphism between their \(\operatorname{L}^1\)-spaces that maps \(1_X\) to \(1_Y\) is the Koopman operator of a measure preserving map. For this we need to introduce \textbf{regular set isomorphisms} (\citep[Section 3.]{lamperti_isometries_1958}). Those are maps \(\tilde{\theta}: \mathcal{A} \to \mathcal{B}\) such that
\begin{itemize}
	\item \(\tilde{\theta}(X\setminus U) = \tilde{\theta}(X)\setminus\tilde{\theta}(U)\) for all \(U \in \mathcal{A}\),
	\item \(\tilde{\theta}(\bigcup_{n \in \mathbb{N}} U_n) = \bigcup_{n \in \mathbb{N}} \tilde{\theta}(U_n)\) for all pairwise disjoint families \((U_n)_{n \in \mathbb{N}}\) in \(\mathcal{A}\),
	\item \(\nu(\tilde{\theta}(U)) = 0\) if and only if \(\mu(U) = 0\) for all \(U \in \mathcal{A}\).
\end{itemize}
Here the equalities in the first two properties only need to hold up to measure zero sets. If such a regular set isomorphism \(\tilde{\theta}\) is given we can define a unique linear operator \(\mathcal{S}_{\tilde{\theta}}: \operatorname{L}^1(\mu) \to \operatorname{L}^1(\nu)\) such that
\[
	\mathcal{S}_{\tilde{\theta}}1_U = 1_{\tilde{\theta}(U)}
\]
for all \(U \in \mathcal{A}\). Now we can state the Banach-Lamperti theorem. In Lamperti's paper (\citep[Theorem 3.1]{lamperti_isometries_1958}) it is only stated for \(X = Y\), this however is never used in its proof.

\begin{theorem}[Banach-Lamperti Theorem]\label{the-Banach-Lamperti}
	Assume that \(\mathcal{T}: \operatorname{L}^1(\mu) \to \operatorname{L}^1(\nu)\) is an isometric isomorphism. Then there exists a regular set isomorphism \(\tilde{\theta}: \mathcal{A} \to \mathcal{B}\) such that
	\[
		\mathcal{T}\varphi = h \cdot \mathcal{S}_{\tilde{\theta}}\varphi
	\]
	for all \(\varphi \in \operatorname{L}^1(\nu)\), with \(h := \mathcal{T}1_X\). Further \(h \neq 0\) holds \(\nu\)-almost everywhere.
\end{theorem}

If we additionally assume that \(\mathcal{T}1_X = 1_Y\) it follows that \(\tilde{\theta}\) is measure preserving and we can find a measure preserving map \(\theta: Y \to X\) such that \(\tilde{\theta}(U) = \theta^{-1}(U)\). This is the content of the following corollary.

\begin{corollary}\label{cor-Banach-Lamperti}
	Assume that \((X, \mathcal{A}, \mu)\) and \((Y, \mathcal{B}, \nu)\) are standard probability spaces and let \(\mathcal{T}: \operatorname{L}^1(\mu) \to \operatorname{L}^1(\nu)\) be an isometric isomorphism with \(\mathcal{T}1_X = 1_Y\). Then there exist a measure preserving map \(\theta: Y \to X\) such that
	\[
		\mathcal{T}\varphi = \varphi \circ \theta
	\]
	for all \(\varphi \in \operatorname{L}^1(\mu)\). The map \(\theta\) is unique up to measure zero sets.
\end{corollary}
\begin{proof*}
	By the Banach-Lamperti theorem (\cref{the-Banach-Lamperti}) and with \(\mathcal{T}1_X = 1_Y\) there exist some regular set isomorphism \(\tilde{\theta}: \mathcal{A} \to \mathcal{B}\) such \(\mathcal{T} = \mathcal{S}_{\tilde{\theta}}\). Then we have for \(U \in \mathcal{A}\)
	\[
		\nu(\tilde{\theta}(U)) = \|1_{\tilde{\theta}(U)}\|_{\operatorname{L}^1(\nu)} = \|\mathcal{T}1_U\|_{\operatorname{L}^1(\nu)} = \|1_U\|_{\operatorname{L}^1(\mu)} = \mu(U),
	\]
	so \(\tilde{\theta}\) is measure preserving.
	
	For \(U \in \mathcal{A}\) we have \(|\mathcal{T}1_U| = \mathcal{T}|1_U|\) and \(\int_Y \mathcal{T}1_U\ \mathrm{d}\nu = \int_Y 1_U\ \mathrm{d}\nu\). By linearity these equations extend to \(\operatorname{span}_{\mathbb{C}}\{1_U: U \in \mathcal{A}\}\). Approximation by those simple functions then yields \(|\mathcal{T}\varphi| = \mathcal{T}|\varphi|\) (i.e. \(\mathcal{T}\) is a lattice homomorphism) and \(\int_Y \mathcal{T}\varphi\ \mathrm{d}\nu = \int_X \varphi\ \mathrm{d}\mu\) for all \(\varphi \in \operatorname{L}^1(\mu)\). Then the existence of \(\theta\) follows from \citep[Theorem 7.20.]{eisner_operator_2015}. The uniqueness up to measure zero sets is \citep[Proposition 7.19.]{eisner_operator_2015}.
\end{proof*}

\section{Automatic Continuity of Almost Group Homomorphisms}\label{sec-Kleppner}

In this appendix we show that a measurable map \(\chi'\) from a compact group \(G\) into a compact metrizable group \(H\) that is a group homomorphism almost everywhere, is almost everywhere equivalent to a continuous group homomorphism. This is based on a theorem by Kleppner (\citep[Theorem 2.]{kleppner_measurable_1989, kleppner_correction_1991}). This theorem is proven in a more general setting, where the groups are not assume to be \(\sigma\)-compact, so their Haar measure may not be \(\sigma\)-finite. The original proof in \citep{kleppner_measurable_1989} contains a mistake, which was corrected in \citep{kleppner_correction_1991}. The corrected version of the theorem has a weaker conclusion. We show below that the incorrect step in the original proof holds in our easier setting with compact groups, so that we can still obtain the stronger conclusion from the original version, as it is found in the theorem below.

\begin{theorem}\label{the-Kleppner}
	Let \(G\) be a compact group and \(H\) a compact metrizable group with Haar probability measures \(\mu_G\) and \(\mu_H\). Let \(\chi': G \to H\) be a measurable map such that
	\[
		\chi'(gh) = \chi'(g) \cdot \chi'(h) \quad\text{and}\quad \chi'(g^{-1}) = \chi'(g)^{-1}
	\]
	for \(\mu_G\)-almost all \(g, h \in G\). Then there exists a continuous group homomorphism \(\chi: G \to H\) such that \(\chi = \chi'\) holds \(\mu_G\)-almost everywhere.
\end{theorem}
\begin{proof*}
	The map \(\chi'\) fulfills all the assumptions of \citep[Theorem 2.]{kleppner_measurable_1989}, since it is measurable, \(G\) and \(H\) are compact and \(\mu_G\) and \(\mu_H\) are finite measures. In the proof of \citep[Theorem 2.]{kleppner_measurable_1989} two group representations \(\pi\) and \(\pi_*\) of \(G\) into the bounded operators on \(\operatorname{L}^2(\mu_H)\) are constructed. It is then proven that for all \(x, y \in \operatorname{L}^2(\mu_H)\) we have
	\begin{equation}\label{eq-Kleppner_Weak_Equivalence}
		\langle \pi(g)x, y \rangle_{\operatorname{L}^2(\mu_H)} = \langle \pi_*(g)x, y \rangle_{\operatorname{L}^2(\mu_H)}
	\end{equation}
	for \(\mu_G\)-almost all \(g \in G\). The null set on which this equation does not hold may depend on \(x\) and \(y\). From this equation it is concluded that also
	\begin{equation}\label{eq-Kleppner_Strong_Equivalence}
		\pi(g) = \pi_*(g)
	\end{equation}
	holds for \(\mu_G\)-almost all \(g \in G\) holds. This is however not true in the general setting used by Kleppner, which leads to the weaker conclusion in \citep{kleppner_correction_1991}. In our setting however, it is true that \eqref{eq-Kleppner_Weak_Equivalence} implies \eqref{eq-Kleppner_Strong_Equivalence}, which is due to the fact that \(\operatorname{L}^2(\mu_H)\) is separable (\citep[Lemma B.4.9]{deitmar_principles_2009}). To prove this, we fix a countable dense subset \(S \subseteq \operatorname{L}^2(\mu_H)\). Since \(S\) is countable it follows from \eqref{eq-Kleppner_Weak_Equivalence} that we find a null set \(N \subseteq G\) such that
	\[
		\langle \pi(g)x, y \rangle_{\operatorname{L}^2(\mu_H)} = \langle \pi_*(g)x, y \rangle_{\operatorname{L}^2(\mu_H)}
	\]
	holds for all \(x, y \in S\) and \(g \in G \setminus N\). Since \(S\) is dense it follows that
	\[
		\pi(g) = \pi_*(g)
	\]
	holds for all \(g \in G \setminus N\). Thus \eqref{eq-Kleppner_Weak_Equivalence} implies \eqref{eq-Kleppner_Strong_Equivalence}, so the conclusion in \citep{kleppner_measurable_1989} is valid in our case and the theorem follows.
\end{proof*}

\section{Miscellaneous}

\begin{lemma}\label{lem-Constant_Absolute_Value_Triangle_Equality}
	Let \((X, \mathcal{A}, \mu)\) be a probability space and let \(\varphi \in \operatorname{L}^1(\mu)\) have constant absolute value and assume it fulfills
	\[
		\bigg| \int_X \varphi\ \mathrm{d}\mu \bigg| = \int_X |\varphi|\ \mathrm{d}\mu.
	\]
	Then \(\varphi\) is almost everywhere constant.
\end{lemma}
\begin{proof*}
	There exist a unique \(\lambda \in \mathbb{S}^1\) such that
	\[
		\bigg| \int_X \varphi\ \mathrm{d}\mu \bigg| = \lambda \int_X \varphi\ \mathrm{d}\mu.
	\]
	This implies that
	\[
		\int_X |\varphi|\ \mathrm{d}\mu = \bigg| \int_X \varphi\ \mathrm{d}\mu \bigg| = \lambda \int_X \varphi\ \mathrm{d}\mu = \int_X \lambda \cdot \varphi\ \mathrm{d}\mu.
	\]
	Then we also have for the real part that
	\[
		\int_X |\varphi| - \operatorname{Re}(\lambda \cdot \varphi)\ \mathrm{d}\mu = \int_X |\varphi|\ \mathrm{d}\mu - \operatorname{Re}\bigg( \int_X \lambda \cdot \varphi\ \mathrm{d}\mu \bigg) = 0,
	\]
	which together with \(|\varphi| - \operatorname{Re}(\lambda \cdot \varphi) \geq 0\) implies \(|\varphi| = \operatorname{Re}(\lambda \cdot \varphi)\) almost everywhere. This in turn implies with \(|\varphi| = |\lambda \cdot \varphi|\) that \(|\varphi| = \lambda \cdot \varphi\). Therefore
	\[
		\varphi = \frac{1}{\lambda} \cdot |\varphi|
	\]
	is almost everywhere constant.
\end{proof*}